\documentclass[preprint,12pt]{elsarticle}

\usepackage{amssymb}
\usepackage{amsmath}
\usepackage{amsthm}

\usepackage{lineno}
\usepackage{mathrsfs}
\usepackage {xfrac}
\usepackage{graphicx}
\usepackage{float}
\usepackage{tikz}
\usetikzlibrary{arrows}
\usepackage{xcolor}
\usepackage{pifont}
\usepackage{hyperref}

\newtheorem{defn}{Definition}
\newtheorem{thm}{Theorem}
\newtheorem{lem}{Lemma}
\newtheorem{exmp}{Example}

\newtheorem{cor}{Corollary}

\newtheorem{rmk}{Remark}

\begin{document}

\begin{frontmatter}



\title{Encoding argumentation frameworks with set attackers to propositional logic systems}


\author[a]{Shuai Tang\corref{cor1}} 
\ead{TangShuaiMath@outlook.com}
\author[a]{Jiachao Wu} 
\ead{wujiachao@sdnu.edu.cn}
\author[a]{Ning Zhou}
\ead{1679659758@qq.com}
\affiliation[a]{organization={School of Mathematics and Statistics, Shandong Normal University},
            addressline={No.1 Daxue Road}, 
            city={Jinan},
            postcode={250358}, 
            state={Shandong},
            country={China}}
\cortext[cor1]{Corresponding author}

\begin{abstract}
Argumentation frameworks ($AF$s) have been a useful tool for approximate reasoning. The encoding method is an important approach to formally model $AF$s under related semantics. The aim of this paper is to develop the encoding method from classical Dung's $AF$s to $AF$s with set attackers ($AFSA$s) including higher-level argumentation frames, Barringer-Gabbay-Woods higher-order $AF$s ($BHAF$s), frameworks with sets of attacking arguments and higher-order set $AF$s. 
Regarding the syntactic structure, we propose higher-order set $AF$s where the target of an attack is either an argument or an attack and the sources are sets of arguments and attacks. Regarding semantics, we propose complete semantics of $BHAF$s and higher-order set $AF$s, and translate $AFSA$s under respective complete semantics to {\L}ukasiewicz's 3-valued propositional logic system.  Moreover, we propose the equational semantics for $AFSA$s, and translate $AFSA$s to fuzzy propositional logic systems. This paper establishes relationships of model equivalence between an $AFSA$ under a given semantics and the encoded formula in a related propositional logic system. By connections of $AFSA$s and propositional logic systems, this paper provides the logical foundations for $AFSA$s associated with complete semantics and equational semantics. 
The results advance the argumentation theory by logical and equational approaches, paving the way for automated reasoning tools in AI, decision support, and multi-agent systems.
\end{abstract}



	


\end{frontmatter}



\section{Introduction}\label{sec1}
Dung's argumentation frameworks ($DAF$s) \cite{dung1995acceptability} offer a formal model for reasoning with conflicting information. Since then, argumentation frameworks ($AF$s) have seen significant development. In application, as shown in \cite{bench2007argumentation, atkinson2021argumentation}, $AF$s have been widely used in areas like artificial intelligence and law, enabling effective handling of real-world scenarios with conflicting data. Theoretically, they have been extended and refined, enhancing their ability to manage more complex reasoning.

Regarding the syntactic level of $AF$s, higher-order $AF$s ($HOAF$s) \cite{barringer2005temporal, gabbay2009semantics, baroni2009encompassing, baroni2011afra},  bipolar $AF$s \cite{karacapilidis2001computer, verheij2003deflog, cayrol2005acceptability}, higher-order bipolar $AF$s \cite{boella2010support, cohen2015approach, gottifredi2018characterizing, cayrol2018structure, cayrol2018argumentation} and frameworks with sets of attacking arguments ($SETAF$s) \cite{nielsen2006generalization, flouris2019comprehensive} have been proposed. Regarding the semantic level of $AF$s, the labelling semantics \cite{caminada2006issue} is developed from Dung's extension semantics \cite{dung1995acceptability}. The original label set $\{in, out, und\}$ in labelling semantics has been extended to the set of numbers in \cite{gabbay2016degrees}. Numerical semantics are introduced by different approaches, such as equational semantics \cite{gabbay2011introducing, gabbay2012equational} and gradual semantics \cite{cayrol2005graduality, barringer2005temporal, gabbay2015equilibrium, amgoud2022evaluation}.

$AF$s and logic systems are connected in different ways. The studies in \cite{besnard2004checking} connect $DAF$s and 2-valued propositional logic system. The encoding approach based on the quantified Boolean formula (QBF) is introduced in \cite{egly2006reasoning, arieli2013qbf}. In the paper \cite{egly2006reasoning}, the authors generalize $DAF$s as frameworks over an induced derivability relation of a deductive system and present the general encoding function for this derivability relation, such that the encoding methodology in \cite{besnard2004checking} serves as a concrete example of their general approach. In \cite{arieli2013qbf}, the authors develop the QBF-based encoding approach associated with the modified Kleene’s three-valued logic via the signed theory induced by a $DAF$. Modal logic systems and $DAF$s are connected in \cite{gabbay2009modal,caminada2009logical,grossi2010logic,villata2012logic, gabbay2015attack}. The works in \cite{gabbay2016attack, fandinno2018constructive} establish relationships between $DAF$s and intuitionistic logic systems.
$DAF$s and the fuzzy set theory are combined by the extension-based approach in \cite{janssen2008fuzzy, wu2016godel, wu2016some, wang2020dynamics, wu2020godel}. $DAF$s and the fuzzy set theory are combined by the labelling-based approach and the iterative algorithm in \cite{pereira2011changing, zhao2022efficient}.  
While \cite{dyrkolbotn2014formal} establishes a meta-logical equivalence between skeptical reasoning in argumentation semantics and {\L}ukasiewicz's three-valued logic ($\mathcal{PL}_3^L$), it does not provide a concrete encoding from argumentation frameworks to $\mathcal{PL}_3^L$ theories.
The paper \cite{strasser2014adaptive} presents a unified adaptive logic framework that encodes argumentation frameworks into logical premises and employs a dynamic proof system to uniformly model all standard semantics for both skeptical and credulous reasoning.
By encoding $DAF$s to 3-valued propositional logic systems ($\mathcal{PL}_3$s) and fuzzy propositional logic systems ($\mathcal{PL}_{[0,1]}$s), the research in \cite{tang2025encoding} connects $DAF$s under complete semantics with $\mathcal{PL}_3$s and connects $AF$s under equational semantics with $\mathcal{PL}_{[0,1]}$s. The paper \cite{dvovrak2023expressiveness} studies the logical encoding approach for $SETAF$s from the perspective of $ADF$s based on the classical binary propositional logic ($\mathcal{PL}_2$). Some logical encoding works associated with the first-order logic for higher-order and bipolar $AF$s are done in \cite{cayrol2017logical,claudette2018logical,cayrol2020logical, lagasquie2021evidential, lagasquie2021necessary, lagasquie2023handling,besnard2023generic}.

However, existing research on logical encoding of argumentation frameworks still faces significant gaps and challenges that hinder their practical application and theoretical unification. First, most existing encoding methods for $AF$s with higher-order structure or collective attackers rely on first-order logic, modal logic, or $\mathcal{PL}_2$. These approaches either introduce excessive syntactic complexity (e.g., quantifiers and complex predicates in first-order logic) that complicates integration with lightweight solvers, or lack native support for handling uncertainty and fuzzy truth-values (e.g., $\mathcal{PL}_2$ cannot directly model fuzzy acceptability degrees). Second, while existing literature has proposed structures for higher-order set argumentation frameworks, these structures suffer from a critical limitation: the elements of set attackers are restricted to arguments, and attacks are not permitted to serve as members of set attackers. From the perspective of theoretical consistency and practical scenarios, attacks, like arguments, have been recognized as entities that can be attacked; naturally, they should also be entitled to initiate collective attacks in the form of set members, just as arguments do. As a result, existing frameworks cannot accommodate scenarios where attacks form set attackers, making it difficult to generalize across diverse complex argumentation interactions. Third, the semantic generalization between multi-valued and fuzzy semantics remains underdeveloped in existing works. Most encodings either stop at three-valued semantics (without extending to continuous fuzzy semantics) or require fundamental modifications to the encoding structure when transitioning between semantic types, failing to achieve a natural and consistent extension. Additionally, practical automated reasoning in AI, decision support, and multi-agent systems demands encoding methods that are both computationally efficient (compatible with lightweight solvers like SAT or fuzzy SAT solvers) and semantically expressive (capable of capturing uncertainty and fuzzy acceptability). Existing methods often sacrifice one for the other—first-order logic-based encodings are expressive but computationally heavy, while simple propositional encodings lack the flexibility to handle complex structures and semantic gradients. Given these limitations, the motivations of this paper are derived from three core considerations, aiming to address the theoretical gaps and practical demands:
\begin{itemize}
	\item To unify $HOAF$s and $SETAF$s into a novel framework—higher-order set argumentation frameworks ($HSAF$s)—uniformly termed argumentation frameworks with set attackers ($AFSA$s), which addresses a critical limitation of existing higher-order set argumentation structures: the restriction that elements of set attackers can only be arguments, excluding attacks. Our proposed $AFSA$s explicitly allow attacks to constitute components of set attackers.
	\item To propose consistent and expressive semantics for $AFSA$s, including complete semantics (for discrete acceptability states) and numerical equational semantics (for fuzzy acceptability degrees). Existing semantics for $HOAF$s and $SETAF$s are either limited to discrete labels (in/out/undecided) or rely on ad-hoc axioms; our semantics fill this gap by providing a unified foundation that naturally connects discrete three-valued reasoning with continuous fuzzy reasoning.
	\item To develop a dedicated encoding method and formal equational approach for $AFSA$s, bridging them with propositional logic systems ($\mathcal{PLS}$s). By avoiding quantifiers and complex predicate constructions, this encoding ensures native compatibility with lightweight solvers, while the formal equational approach guarantees the existence of solutions (via Brouwer’s Fixed-Point Theorem), laying the groundwork for efficient automated reasoning in practical applications.
\end{itemize}

This work improves the encoding method from $DAF$s to $HSAF$s, lays logical foundations for the $AFSA$ by proving model equivalence between an $AFSA$ and its encoded formula, and provides a formal equational approach with guaranteed solution existence for the numerical equational semantics of $AFSA$s. In summary, this paper makes three core contributions:
\begin{itemize}
	\item \textbf{Syntactic innovation}: Propose the $HSAF$ as a higher-order extension of the $SETAF$, enabling mutual interactions between set attacks (both being attacked and attacking others).
	\item \textbf{Semantic enrichment}: Propose complete semantics for Barringer-Gabbay-Woods higher-order argumentation frameworks ($BHAF$s) and $HSAF$s, and design numerical equational semantics for $AFSA$s.
	\item \textbf{Encoding advancement}: Translate $AFSA$s under respective complete semantics to $\mathcal{PL}_3^L$, translate $AFSA$s under equational semantics to the $\mathcal{PL}_{[0,1]}$s, and explore the relationship between general real equational semantics and continuous fuzzy normal encoded equational semantics.
\end{itemize}

We structure this paper as a progression from the particular $AF$ to the general $AF$. The $BHAF$ is the generalization of the higher-level argumentation frame ($HLAF$). The $HSAF$ is the generalization of the $SETAF$ and the $BHAF$ (as illustrated in Figure \ref{fig:af_hierarchy}). We follow this progression because we want to develop the argumentation framework theory by following original basic works. This paper is arranged as follows. In Section 2, we give some basic knowledge. In Sections 3-6, we respectively encode $HLAF$s, $BHAF$s, $SETAF$s and $HSAF$s to $\mathcal{PLS}$s. Section 7 analyzes previous studies and compares them with our approach. Section 8 is the conclusion.
\begin{figure}[t]
	\centering
	\begin{tikzpicture}[node distance=3cm, auto]
		\node (HSAF) {$HSAF$};
		\node (BHAF) [below left of=HSAF] {$BHAF$};
		\node (SETAF) [below right of=HSAF] {$SETAF$};
		\node (HLAF) [below left of=BHAF] {$HLAF$};
		\draw[->] (HLAF) -- (BHAF);
		\draw[->] (BHAF) -- (HSAF);
		\draw[->] (SETAF) -- (HSAF);
	\end{tikzpicture}
	\caption{Hierarchy for $AFSA$s.}  
	\label{fig:af_hierarchy}  
\end{figure}
\section{Preliminaries}\label{pre}
In this section, we will review the basic background on $AF$s, $\mathcal{PLS}$s and related methods.

\subsection{The syntactic structures and semantics of argumentation frameworks}
\label{subsec1}
Firstly, we respectively list the definitions of  $DAF$s \cite{dung1995acceptability}, $HLAF$s \cite{gabbay2009semantics} and $SETAF$s \cite{flouris2019comprehensive} below.
\begin{defn}[\cite{dung1995acceptability}]
	A \emph{Dung's argumentation framework} ($DAF$) is a pair $(A^D, R^D)$ where $A^D$ is a finite set of arguments and $R^D \subseteq A^D \times A^D$ is the attack relation. 
\end{defn}
\begin{defn}[\cite{gabbay2009semantics}]\label{hlaf}
	A \emph{0-level argumentation frame} ($HLAF^{(0)}$) is a pair $(A^{HL(0)}, R^{HL(0)})$, where $A^{HL(0)}$ is a finite set of arguments and $R^{HL(0)}\subseteq A^{HL(0)} \times A^{HL(0)}$ is a $0$-level attack relation. An \emph{$n$-level argumentation frame} ($HLAF^{(n)}$) is a pair $(A^{HL(0)}, R^{HL(n)})$, where $n\geqslant1$, $A^{HL(n)}=A^{HL(0)}\cup A^{HL(0)} \times A^{HL(n-1)}$, the $n$-level attack relation $R^{HL(n)}\subseteq A^{HL(0)} \times A^{HL(n)}$, and $R^{HL(n)}\nsubseteq A^{HL(0)} \times A^{HL(n-1)}$. 
\end{defn}
Briefly, we use the ``$HLAF$" to represent an argumentation frame at any finite level in this paper.
Definition \ref{hlaf} is equivalent to the definition of $HLAF$s given in \cite{gabbay2009semantics}. Both definitions express the syntactic concept that higher-order attacks in an $HLAF$ are directed from arguments to attacks.
By this definition, an $HLAF^{(0)}$ coincides with a $DAF$. Specifically, $A^{HL(1)} = A^{HL(0)} \cup (A^{HL(0)} \times A^{HL(0)})$ encompasses all arguments (from $A^{HL(0)}$) and all possible 0-level attacks (i.e., attacks between arguments). Notably, all elements of $A^{HL(1)}$—both arguments and 0-level attacks—are treated as arguments for constructing the relations of the next level. Accordingly, $A^{HL(0)} \times A^{HL(1)}$ includes two types of attacks: all possible 0-level attacks (from $A^{HL(0)} \times A^{HL(0)}$) and 1-level attacks (i.e., attacks from arguments to 0-level attacks). Importantly, this set encompasses not only 0-level attacks but also higher-level attacks (specifically 1-level attacks) at this stage. Generalizing this logic, $A^{HL(2)} = A^{HL(0)} \cup (A^{HL(0)} \times A^{HL(1)})$ covers all arguments and all possible attacks up to the 2-level, with all elements of $A^{HL(2)}$ again serving as arguments for the next level. Extending this recursive process, $A^{HL(n)} = A^{HL(0)} \cup (A^{HL(0)} \times A^{HL(n-1)})$ includes all arguments and all possible attacks up to the n-level (encompassing all lower-level attacks from 0 to n-1). Furthermore, $R^{HL(n)} \subseteq A^{HL(0)} \times A^{HL(n)}$ and $R^{HL(n)} \nsubseteq A^{HL(0)} \times A^{HL(n-1)}$ ensure that the $n$-level attack relation $R^{HL(n)}$ genuinely includes $n$-level attacks (a subset of higher-level attacks) not present in the $(n-1)$-level attack relation. An updated definition of $HLAF$s is provided in Subsection \ref{sub3.1}, while the definition of $BHAF$s is detailed in Subsection \ref{sub4.1}. We use the ``$HOAF$'' to represent an $HLAF$ or a $BHAF$ in this paper.

\begin{defn}[\cite{flouris2019comprehensive}]
	A \emph{framework with sets of attacking arguments} ($SETAF$) is a pair $(A^S, R^S)$ where $A^S$ is a finite set of arguments and $R^S\subseteq (2^{A^S} \setminus\{\emptyset\}) \times A^S$ is the set attack relation. 
\end{defn}
$(S,a)\in R^S$ is usually denoted as $S\triangleright a$, where $a$ belongs to $A^S$ and $S\subseteq2^{A^S} \setminus\{\emptyset\}$ is a set attacker of $a$.
A set $B\subseteq A^S$ attacks an argument $a$ (i.e., $B$ is a set attacker of $a$) if $\exists B'\subseteq B$ such that $B'\triangleright a$.
A set $B'$ is called a minimal attacker of $a$, if $\nexists B''\subsetneq B'$ such that $B''\triangleright a$. 
In this paper, we suppose that each attacker of $a$ is a minimal attacker, which means that $S$ is a set attacker of $a$ iff (if and only if) $S\triangleright a$.

Secondly, we review key concepts related to semantics. Roughly speaking, a semantics is a method for determining the acceptability of arguments (and attacks in higher-order frameworks). A given semantics can be characterized using either extension-based or labelling-based approaches. A \emph{labelling} of an $AF$ is a function that assigns to each argument (and each attack in an $HOAF$) a label from a given set. The extension-based approach characterizes a semantics by selecting sets of arguments that satisfy certain semantic constraints. In contrast, a \emph{labelling-based semantics} for a class of $AF$s is a function that assigns to each $AF$ in that class a set of valid labellings. This paper employs exclusively the labelling-based approach. If the set of labels for a labelling $lab$ is a subset of the unit interval $[0,1]$, we call $lab$ a \emph{numerical labelling}. We only consider numerical labellings in this paper. Specifically, if the label set is $\{0, \frac{1}{2}, 1\}$ or $[0,1]$, we call $lab$ a \emph{3-valued labelling} or a \emph{$[0,1]$-valued labelling}, respectively.

\begin{defn}
For a given $HOAF$ (respectively, $SETAF$), each labelling $lab$ that satisfies a given semantics is called a \emph{model} of the $HOAF$ (respectively, $SETAF$) under this semantics. 
\end{defn}
Now let us review the complete semantics of $HLAF$s \cite{gabbay2009semantics} and $SETAF$s \cite{nielsen2006generalization, flouris2019comprehensive}.
\begin{defn}[\cite{gabbay2009semantics}]\label{defn16}\footnote{Note that the original definition of complete semantics for $HLAF$s in \cite{gabbay2009semantics} contains minor errors in the case of $\|\beta\|=\frac{1}{2}$, where this scenario conflicts with the case of $\|\beta\|=1$. Definition \ref{defn16} is a corrected version adopted in this work.}
	For an $HLAF=(A^{HL(0)}, R^{HL(n)})$, a \emph{complete labelling} for the $HLAF$ is a function $\|\cdot\| : A^{HL(0)} \cup R^{HL(n)}\rightarrow\{0, 1, \frac{1}{2}\}$ such that $\forall\beta\in A^{HL(0)}\cup R^{HL(n)}$
	\begin{equation*}
		\|\beta\| =\begin{cases}
			1 & \quad \text{ iff }\forall a ((a, \beta)\in R^{HL(n)}): \|a\|=0 \text{ or } \|(a, \beta)\|= 0\\
			0 & \quad \text{ iff } \exists a ((a, \beta)\in R^{HL(n)}): \|a\|= 1 \text{ and } \|(a, \beta)\|=1\\
			\frac{1}{2} & \quad \text{ otherwise}
		\end{cases}.
	\end{equation*}
\end{defn}
The complete semantics of $HLAF$s describes the core semantic properties: an attacking action is valid if and only if both the attacker and its associated attack are valid; an attacking action is invalid if and only if either the attacker or its associated attack is invalid; and all other cases correspond to a semi-valid attacking action.
In an $HLAF^{(0)}$ (i.e., a $DAF$), no attack is attacked by any argument, so each attack has a value of 1. When we only consider labels on arguments in the $HLAF^{(0)}$, Definition \ref{defn16} is restricted to the complete labelling semantics of the $DAF$.
\begin{defn}[\cite{flouris2019comprehensive}]
	For a given $SETAF = (A^S, R^S)$, a \emph{complete labelling} of the $SETAF$ is a function $lab^S: A^S \rightarrow \{1, 0, \frac{1}{2}\}$ such that $\forall a \in A^S$
	\begin{equation*}
		lab^S(a) = \begin{cases}
			1 & \quad \text{ iff } \forall S ((S,a)\in R^S) \exists b \in S: lab^S(b) = 0\\
			0 & \quad \text{ iff } \exists S ((S,a)\in R^S) \forall b \in S: lab^S(b) = 1\\
			\frac{1}{2} & \quad \text{ otherwise }
		\end{cases}. 
	\end{equation*}
\end{defn}
The complete semantics of $SETAF$s describes the core semantic features: a set attacking action is valid if and only if the set attacker is valid if and only if each member in the set attacker is valid; a set attacking action is invalid if and only if the set attacker is invalid if and only if there is an invalid member in the set attacker; and all other cases correspond to a semi-valid set attacking action.

Roughly speaking, an equational semantics of $HOAF$s is constructed by giving an equational system for the $HOAF$s such that the models of an $HOAF$ under the equational semantics are the solutions of the equational system applied on the $HOAF$. Equational semantics one-to-one correspond to equational system, so we do not distinguish them. The equational semantics of $SETAF$ is similar to that of $HOAF$s. We will give the formal equational approach in Section \ref{sec6}.

\subsection{Essential knowledge of propositional logic systems}
For in-depth knowledge of $\mathcal{PLS}$s, interested readers may consult \cite{klir1995fuzzy,Hajek1998,Klement2000Triangular, bergmann2008introduction, belohlavek2017fuzzy}. In this subsection, we present only the foundational concepts that are well-established in the literature.

In a $\mathcal{PLS}$, propositional connectives  $\neg, \wedge, \vee$, $\rightarrow$ and $\leftrightarrow$ respectively denote the negation,  conjunction, disjunction, implication and bi-implication connectives, where $a\leftrightarrow b:=(a\rightarrow b)\wedge(b\rightarrow a)$ for any formula $a$ and formula $b$. In this paper, we consider only the $\mathcal{PLS}$ that has a numerical assignment domain (a subset of $[0,1]$) and its corresponding numerical semantic system.
\begin{defn}
	An \emph{assignment of propositional variables} in a $\mathcal{PLS}$ is a function $\|\cdot\|: S\rightarrow L$, where  $S$ is the set of all propositional variables in the $\mathcal{PLS}$ and $L$ is a subset of $[0, 1]$.
	An \emph{assignment on a formula $\phi$} in a $\mathcal{PLS}$ is the function $\|\cdot\|: S_{\phi}\rightarrow L$, where  $S_{\phi}$ is the set of all propositional variables in $\phi$ and $L$ is a subset of $[0, 1]$.
\end{defn}
In each $\mathcal{PLS}$, there is a constant $\bot$ such that $\|\bot\|=0$ and $\|\neg\bot\|=1$.
Without causing ambiguity we also denote the evaluation of  all formulas by $\|\cdot\|$, since, for a given $\mathcal{PLS}$, an assignment $\|\cdot\|$ of propositional variables extends uniquely to the evaluation of all formulas. 
\begin{defn}
	A \emph{model} of a formula $\phi$ in a $\mathcal{PLS}$ is an assignment $\|\cdot\|$ on $\phi$ such that $\|\phi\|=1$. We denote that $\|\cdot\|$ is a model of $\phi$ in the $\mathcal{PLS}$ by $\|\cdot\| \models_{\mathcal{PLS}} \phi$.
\end{defn}

In the $\mathcal{PL}_3^L$ \cite{Lukasiewicz1970ThreeValued}, the evaluation of all formulas with connectives $\neg, \wedge$ and $\vee$ is extended as:
\begin{itemize}
	\item $\|\neg a\|=1-\|a\|$
	\item $\|a\wedge b\|=\min\{\|a\|, \|b\|\}$
	\item $\|a\vee b\|=\max\{\|a\|, \|b\|\}$
\end{itemize}
where  $a, b$ are propositional formulas. Moreover, the connective $\rightarrow$ is interpreted as $\|a\rightarrow b\|=\|a\|\Rightarrow\|b\|$ where the truth degree table of operation $\Rightarrow$ is listed as Table \ref{Rightarrow}.
\begin{table}
	\centering
	\begin{tabular}{l r r r}
		\hline
		$\Rightarrow$ & 0 & $\frac{1}{2}$ & 1 \\
		\hline
		0 & 1 & 1 & 1 \\
		$\frac{1}{2}$ & $\frac{1}{2}$ & 1 & 1 \\
		1 & 0 & $\frac{1}{2}$ & 1 \\
		\hline
	\end{tabular}
		\caption{Truth degree table of $\Rightarrow$ in $\mathcal{PL}_3^L$}
	\label{Rightarrow}
\end{table}
Since $a\leftrightarrow b:=(a\rightarrow b)\wedge(b\rightarrow a)$, we have  $\|a\|\Leftrightarrow\|b\|=\|a\leftrightarrow b\|=\min\{\|a\rightarrow b\|, \|b\rightarrow a\|\}=\min\{\|a\|\Rightarrow\|b\|, \|b\|\Rightarrow\|a\|\}$. Thus, we have the truth degree table of $\Leftrightarrow$ in the $\mathcal{PL}_3^L$ as Table \ref{Leftrightarrow}. 
\begin{table}
	\centering
	\begin{tabular}{l r r r}
		\hline
		$\Leftrightarrow$ & 0 & $\frac{1}{2}$ & 1 \\
		\hline
		0 & 1 & $\frac{1}{2}$ & 0 \\
		$\frac{1}{2}$ & $\frac{1}{2}$ & 1 & $\frac{1}{2}$ \\
		1 & 0 & $\frac{1}{2}$ & 1 \\
		\hline
	\end{tabular}
		\caption{Truth degree table of $\Leftrightarrow$ in $\mathcal{PL}_3^L$}
	\label{Leftrightarrow}
\end{table}

In $\mathcal{PL}_3^L$, the implication operator is interpreted as $\|a \rightarrow b\| = \min\{1, 1-\|a\|+\|b\|\}$. For the intermediate truth value $\frac{1}{2}$ (corresponding to ``undecided" in argumentation), this gives $\|\frac{1}{2} \rightarrow \frac{1}{2}\| = 1$, which exactly captures the core intuition of argumentation frameworks: an argument's undecided state naturally corresponds to the undecided state of its attackers. In some other logic systems, such as Kleene's 3-valued propositional logic system, the implication operator is interpreted as $\|a \rightarrow b\| = \max\{1-\|a\|, \|b\|\}$. Here, $\|\frac{1}{2} \rightarrow \frac{1}{2}\| = \frac{1}{2}$ (not 1), failing to establish the necessary correspondence between an argument's undecided state and that of its attackers—undermining the strict alignment between logical encoding and argumentation semantics.

Note that in this paper we denote ``$A$ if and only if $B$'' by ``$A \Longleftrightarrow B$'' and ``$B$ only if $A$'' (i.e., if $B$ then $A$) by ``$B \Longrightarrow A$'', where $A$ and $B$ are propositions. Next, let us review some basic knowledge of fuzzy operations in $\mathcal{PL}_{[0, 1]}s$.

\begin{defn}[\cite{klir1995fuzzy}]
	A function $N : [0, 1] \rightarrow [0, 1]$ is called a \emph{negation} if it satisfies:
	\begin{itemize}
		\item \emph{Boundary Conditions}: $N(0) = 1$ and  $N(1) = 0$.
		\item \emph{Antimonotonicity}: $m \leq n \Longrightarrow N(m) \geq N(n)$ for all $m, n\in[0, 1]$. 
	\end{itemize} 
\end{defn}
If a negation $N$ is continuous, then $N$ is called a \emph{continuous negation}.
A negation $N$ is called a \emph{standard negation} if $\forall m\in[0, 1]$: $N(m)=1-m$. 

A t-norm (short for ``triangular norm") is a fundamental mathematical operator in fuzzy logic that generalizes the classical logical conjunction (``AND") to handle truth degrees (values in the interval $[0,1]$). Unlike classical logic, t-norms quantify the ``strength" of the conjunction when $A$ and $B$ are partially true.
\begin{defn}[\cite{Hajek1998}]
	A \emph{triangular norm} (t-norm) is a function $T : [0, 1]^2 \rightarrow [0, 1]$ satisfying the following axioms for all $m, n, p, q \in [0,1]$:
	\begin{itemize}
		\item \emph{Identity}: $T (m, 1) = m$.
		\item \emph{Commutativity}: $T (m, n) = T (n, m)$.
		\item \emph{Associativity}: $T (m, T (n, p)) = T (T (m, n), p)$.
		\item \emph{Monotonicity}: If $m\leq p$ and $n \leq q$ then $T (m, n) \leq T (p, q)$.
	\end{itemize}
\end{defn}
The definition of triangular conorms (t-conorms) is dual to that of t-norms (with 0 as the identity element instead of 1). For brevity, we omit the formal definition here. We also denote $T (m, n)$ as $m\ast n$. 
Due to the commutativity and associativity properties, for the t-norm operation applied to the multivariate set $\{m_1,m_2,\dots,m_i,\dots,m_j\}$ with $m_i\in[0,1]$, we can denote it as $m_1\ast m_2\ast\cdots\ast m_i\ast\cdots\ast m_j$.
The following are three key examples of continuous t-norms \cite{Hajek1998}:
\begin{itemize}
	\item G\"{o}del t-norm (minimum): $T_G (m, n)=\min\{m, n\}$.
	\item {\L}ukasiewicz t-norm: $T_L (m, n)= \max\{0, m+n-1\}$.
	\item Product t-norm (algebraic product): $T_P (m, n) = m\cdot n$.
\end{itemize}
In terms of theory, G\"{o}del, {\L}ukasiewicz, and Product t-norms serve as the generators of continuous t-norms, forming the basis for their ordinal sum representation \cite{Klement2000Triangular}. They also act as the limit points of the Frank t-norm family, with strict t-norms derivable from the Product t-norm via monotonic transformation and nilpotent t-norms derivable from the {\L}ukasiewicz t-norm \cite{Hajek1998, Klement2000Triangular}. In terms of application, they are respectively suitable for scenarios of strict conjunction (all conditions satisfied) \cite{Hajek1998}, linear compromise (partial satisfaction) \cite{cignoli2000algebraic}, and joint modeling of independent factors (probabilistic reasoning) \cite{nguyen2000first}, collectively serving as core tools for theoretical construction and practical problem-solving in the field of fuzzy logic. In argumentation theory, G\"{o}del t-norm and Product t-norm have been widely applied \cite{gabbay2012equational, gabbay2015probabilistic}, while {\L}ukasiewicz t-norm seems to be neglected. In this paper, we incorporate all three significant t-norms.

A fuzzy implication generalizes the classical logical implication to fuzzy logic. It quantifies the truth degree of the conditional relationship between two fuzzy propositions (antecedent $A$ and consequent $B$), where $A$ and $B$ each take values in $[0,1]$.
\begin{defn}[\cite{Hajek1998}]
	A function $I : [0, 1]^2 \rightarrow[0, 1]$ is a \emph{fuzzy implication} if it satisfies that for all $m, n, p\in[0, 1]$:
	\begin{itemize}
		\item \emph{Left Antitonicity}: $m \leq n\Longrightarrow I(m, p) \geq I(n, p)$.
		\item \emph{Right Monotonicity}: $n \leq p\Longrightarrow I(m, n) \leq I(m, p)$.
		\item \emph{Boundary Conditions}: $I(1, 0) = 0$ and $I(0, 0) = I(1, 1) = 1$.
	\end{itemize}
\end{defn}
An R-implication (short for ``residual implication") is a specific type of fuzzy implication constructed from a t-norm via a mathematical operation called residuation. It is the most widely used fuzzy implication in applications (e.g., fuzzy control, expert systems) because it aligns closely with intuitive reasoning about conditional rules.
\begin{defn}[\cite{Hajek1998}]\label{def11}
	For a t-norm $T$, the \emph{residual implication} (or \emph{$R$-implication}) $I_T$ induced by $T$ is defined by:
	 \begin{equation*}
	 	I_T (m, n):= sup\{p\mid T (m, p) \leq n\}
	 \end{equation*}
	 where $m,n,p\in[0,1]$.
\end{defn}
The following are three key examples of $R$-implications associated with $T_G, T _L$ and $T_P$ \cite{Hajek1998}:
\begin{itemize}
	\item $I_G(m, n) = 
	\begin{cases}
		1 & \quad m\leq n \\
		n & \quad m>n 
	\end{cases},$
	\item $I_L(m, n) = \min\{1 - m +n, 1\},$
	\item $I_P(m, n)=\begin{cases}
		1 & \quad m\leq n \\
		\frac{n}{m} & \quad m>n
	\end{cases}$.
\end{itemize}
The restriction of $I_L$ to the 3-valued domain $\{0,1,\frac{1}{2}\}$ is the implication operator in $\mathcal{PL}_3^L$ (Table \ref{Rightarrow}). The restrictions of $I_G$ and $I_P$ to $\{0,1,\frac{1}{2}\}$ are identical, and their difference from the implication in Table \ref{Rightarrow} lies in that $\|\frac{1}{2} \rightarrow 0\| = 0$. However, as we will see in the following, we only handle the models of formulas involving bi-implication and only use the relations on the main diagonals in Table \ref{Rightarrow} or Table \ref{Leftrightarrow}. Thus, there is no difference in treating the theory of argumentation encoding in this paper. Therefore, we mainly adopt $\mathcal{PL}_3^L$ and do not elaborate on the other two 3-valued logic systems, which are obtained by substituting the implication operator in $\mathcal{PL}_3^L$ with the restricted forms of $I_G$ or $I_P$.

In this paper, evaluations of formulas with connectives $\neg, \wedge$, $\rightarrow$ and $\leftrightarrow$ in any $\mathcal{PL}_{[0, 1]}$ are extended as \cite{Hajek1998, Klement2000Triangular}:
\begin{itemize}
	\item $\|\neg a\|=N(\|a\|),$
	\item $\|a\wedge b\|=T(\|a\|, \|b\|),$
	\item $\|a\rightarrow b\|=I_T(\|a\|, \|b\|),$
\end{itemize} 
where $a, b$ are propositional formulas in the $\mathcal{PL}_{[0, 1]}$.
Since $a\leftrightarrow b$ is defined as $a\leftrightarrow b:=(a\rightarrow b)\wedge(b\rightarrow a)$, so we have
$\|a\leftrightarrow b\|=\|(a\rightarrow b)\wedge(b\rightarrow a)\|=T(\|a\rightarrow b\|, \|b\rightarrow a\|)=T(I_T(\|a\|, \|b\|), I_T(\|b\|, \|a\|))$.

 We use $\mathcal{PL}_{[0,1]}^G$ \cite{esteva2000residuated} to denote the fuzzy logic system with standard negation $N$,  G\"{o}del t-norm $T_G$ and R-implication $I_G$. Similarly, $\mathcal{PL}_{[0,1]}^P$ \cite{esteva2000residuated} denotes the fuzzy logic system with standard negation $N$,  Product t-norm $T_P$ and R-implication $I_P$, while $\mathcal{PL}_{[0,1]}^L$ \cite{Hajek1998} denotes the fuzzy logic system with standard negation $N$,  {\L}ukasiewicz t-norm $T_L$ and R-implication $I_L$. Note that in this paper we focus on using the semantic operations in a fuzzy propositional logic system, and we design the $\mathcal{PL}_{[0,1]}$ equipped with proper semantic operations.
\begin{rmk}
	In \cite{Hajek1998}, the conjunctive connective is denoted by the symbol $\&$, whereas in \cite{Klement2000Triangular}, it is denoted by $\wedge$. This paper adopts the latter notation. For formulas $a$ and $b$ in $\mathcal{PL}_{[0, 1]}$ systems, the evaluation of the formula $a \wedge b$ is given by:
	\begin{equation*}
		\|a \wedge b\|=T(\|a\|, \|b\|)=\|a\|\ast\|b\|.
	\end{equation*}
Our approach focuses on the $\mathcal{PLS}$ equipped with a numerical semantic system, such as $\mathcal{PL}_3^L$, $\mathcal{PL}_{[0,1]}^G$, $\mathcal{PL}_{[0,1]}^P$, $\mathcal{PL}_{[0,1]}^L$, Kleene's $\mathcal{PL}_3$, {\L}ukasiewicz's $n$-valued $\mathcal{PLS}$, among others. This paper primarily employs $\mathcal{PL}_3^L$, $\mathcal{PL}_{[0,1]}^G$, $\mathcal{PL}_{[0,1]}^P$, and $\mathcal{PL}_{[0,1]}^L$.
\end{rmk}

\subsection{The encoding method and the translating method}\label{sub2.3}
In this subsection, we will give definitions of the encoding method and the translating method for $AFSA$s. Similar definitions are presented for $DAF$s in \cite{tang2025encoding}.

Let the arguments and attacks of $\mathcal{HOAF}$ be in one-to-one correspondence with propositional variables in a $\mathcal{PLS}$, where $\mathcal{HOAF}$ denotes the class of all $HLAF$s or all $BHAF$s. Let $\mathcal{F_{\mathcal{PL}}}$ be the set of all formulas in the $\mathcal{PLS}$. 

\begin{defn}
An \emph{encoding} of $\mathcal{HOAF}$ is a function $ec: \mathcal{HOAF} \rightarrow \mathcal{F_{\mathcal{PL}}}$, s.t. the set of all propositional variables in the formula $ec(HOAF)$ is the set of all arguments and attacks in the $HOAF$. We call $ec(HOAF)$ the \emph{encoded formula} of the $HOAF$.
\end{defn}
Encoding $HOAF$s to a $\mathcal{PLS}$ refers to the process of mapping them to $\mathcal{F_{\mathcal{PL}}}$ via a given encoding function. The models of the formula $ec(HOAF)$ can be viewed as the models of a theory defined by the singleton set $ec(HOAF)$ in a $\mathcal{PLS}$.

\begin{defn}\label{defn13}
Given an encoding function $ec$, the \emph{encoded argumentation semantics} of $\mathcal{HOAF}$ is the semantics obtained by defining models of each $HOAF$ as models of the $ec(HOAF)$ in a given $\mathcal{PLS}$.
A \emph{fuzzy encoded argumentation semantics} of $\mathcal{HOAF}$ is an encoded argumentation semantics where the $\mathcal{PLS}$ is instantiated as a $\mathcal{PL}_{[0,1]}$.
\end{defn}

\begin{defn}
A \emph{translation} of $\mathcal{HOAF}$ under a given semantics is an encoding function $tr: \mathcal{HOAF} \rightarrow \mathcal{F_{\mathcal{PL}}}$, $HOAF \mapsto tr(HOAF)$, s.t. models of $tr(HOAF)$ in the $\mathcal{PLS}$ are in one-to-one correspondence with the models of the $HOAF$ under the given semantics. 
\end{defn}

\begin{defn}
A given semantics of  $\mathcal{HOAF}$ is called \emph{translatable} if there exists a translation of $\mathcal{HOAF}$.
\end{defn} 
Obviously, an encoded argumentation semantics of $\mathcal{HOAF}$ is translatable. 

For the case of $SETAF$s, we let arguments of $\mathcal{SETAF}$ one-to-one correspond to propositional variables in a $\mathcal{PLS}$, where $\mathcal{SETAF}$ is the set of all $SETAF$s. Then related conceptions are similar to those of $HOAF$s.

Within a $\mathcal{PLS}$, the assignment for a formula $\phi$ is conceptually parallel to the labelling for an $AF$. Thus, for the sake of ease of expression, we consider the terms ``assignment'' and ``labelling'' interchangeable. Henceforth, we usually represent the labelling function $lab$ using the notation $\|\cdot\|$. In particular, the truth degree of an argument $a$ (or an attack in $HOAF$) is denoted as $\|a\|$. 
\section{Encoding higher-level argumentation frames}\label{sec3}
In this section, we encode $HLAF$s under complete semantics to the $\mathcal{PL}_3^L$.  Then, we present three typical equational semantics for $HLAF$s and encode $HLAF$s to $\mathcal{PL}_{[0,1]}$s under these these equational semantics. It is shown that these semantics are translatable by proving the model equivalence. 
\subsection{Encoding $HLAF$s under complete semantics}\label{sub3.1}
In this subsection, we first give the normal encoding function of $HLAF$s and then give the theorem of model equivalence. In order to encode $HLAF$s, we extend $HLAF$s by adding some imaginary argument and attacks firstly. 
\begin{defn}\label{hlafdef}
	For an $HLAF=(A^{HL(0)}, R^{HL(n)})$, the extended $HLAF$ is a pair $(\mathsf{A}^{HL(0)}, \mathsf{R}^{HL(n)})$, where $\mathsf{A}^{HL(0)}=A^{HL(0)}\cup \{\bot\}$, $\mathsf{R}^{HL(n)}=R^{HL(n)}\cup \{(\bot, \beta)\mid \beta\in A^{HL(0)}\cup R^{HL(n)} \text{ and } \nexists a\in A^{HL(0)}: (a, \beta)\in R^{HL(n)}\}$, and for any assignment $\|\cdot\|$ of the extended $HLAF$ let $\|\bot\|=0$ and $\|(\bot, \beta)\|=1$ (if $(\bot, \beta)\in \mathsf{R}^{HL(n)}$).
\end{defn}
From this definition, for any argument or attack that is attacked in $(A^{HL(0)}, R^{HL(n)})$, we do nothing. Let $\beta$ denote an argument or an attack that does not have any attacker in $(A^{HL(0)}, R^{HL(n)})$. We only add an imaginary argument $\bot$ and associated attack $(\bot, \beta)$ for $\beta$. Intuitively, introducing a zero-valued attacker for an argument or an attack will not affect its validity. Let us review Definition \ref{defn16}. If we add $\bot$ and $(\bot, \beta)$ for $\beta$, then $\beta$ has only one attacker $\bot$ ($\|\bot\|=0$) and associated attacker is $(\bot, \beta)$ ($\|(\bot, \beta)\|=1$), and then we get $\|\beta\|=1$ from Definition \ref{defn16}. If we do not add $\bot$ and $(\bot, \beta)$ for $\beta$, then $\nexists y: (y, \beta) \in R^{HL(n)}$, and then we still have $\|\beta\|=1$ from Definition \ref{defn16}. Thus, this extension preserves the complete semantics of the original $HLAF$. With a slight abuse of notation, symbols $HLAF$ and $\mathcal{HLAF}$ will henceforth refer to the extended $HLAF$ and the set of all extended $HLAF$s, respectively. The meanings of the notations $A^{HL(0)}$ and $R^{HL(n)}$ remain unchanged. Formally, we present the complete semantics for the extended $HLAF$.
\begin{defn}\label{newcom}
	A \emph{complete labelling} for an $HLAF=(\mathsf{A}^{HL(0)}, \mathsf{R}^{HL(n)})$ is a function $\|\cdot\| : \mathsf{A}^{HL(0)}\cup \mathsf{R}^{HL(n)}\rightarrow\{0, 1, \frac{1}{2}\}$ such that $\|\bot\|=0$, $\|(\bot, \alpha)\|=1$ (if $(\bot, \alpha)\in \mathsf{R}^{HL(n)}$) and for each $\beta\in A^{HL(0)}\cup R^{HL(n)}$
	\begin{equation*}
		\|\beta\| =\begin{cases}
			1 & \quad \text{ iff } \forall a ((a, \beta)\in \mathsf{R}^{HL(n)})): \|a\|=0 \text{ or } \|(a, \beta)\|= 0\\
			0 & \quad \text{ iff } \exists a ((a, \beta)\in \mathsf{R}^{HL(n)})): \|a\|= 1 \text{ and } \|(a, \beta)\|=1\\
			\frac{1}{2} & \quad \text{ otherwise}
		\end{cases}.
	\end{equation*}
\end{defn}
The purpose of establishing these updated definitions is to facilitate the encoding of $HLAF$s, particularly for arguments that have no attackers. An alternative way to handle such particular cases is to define the conjunction of $\emptyset$ as $\bot$, which is an approach applicable to $HLAF$s but not to $SETAF$s. We will elaborate on the difficulties that arise in the case of $SETAF$ in Section \ref{SETAF}. Note that for any assignment $\|\cdot\|$ in any $\mathcal{PLS}$, we always let $\|\bot\|=0$ and $\|(\bot,\alpha)\|=1$ (if $(\bot,\alpha)\in\mathsf{R}^{HL(n)}$ for an $HLAF=(\mathsf{A}^{HL(0)}, \mathsf{R}^{HL(n)})$). Then based on new definitions for $HLAF$s, we present the normal encoding function of $HLAF$s below.
\begin{defn}\label{defn17}
	The \emph{normal encoding} of $HLAF$s w.r.t. the $\mathcal{PL}_3^L$ is the encoding function $ec_{HL}$ such that for a given $HLAF=(\mathsf{A}^{HL(0)}, \mathsf{R}^{HL(n)})$ we have
	\begin{equation*}
		ec_{HL}(HLAF)=\bigwedge_{\beta\in A^{HL(0)}\cup R^{HL(n)}}(\beta\leftrightarrow\bigwedge_{(a_i, \beta)\in \mathsf{R}^{HL(n)}}\neg(a_i\wedge r_{a_i}^{\beta}))
	\end{equation*}
	where $r_{a_i}^{\beta}$ denotes the attack $(a_i, \beta)$.
\end{defn}
According to Definition\ref{defn17}, in the encoded formula, any $\beta \in A^{HL(0)} \cup R^{HL(n)}$ refers to a non-imaginary argument or attack, whereas $a_i$ may be the imaginary argument $\bot$. Since the encoded formula includes the imaginary argument $\bot$ and imaginary attacks $(\bot, \beta)$, we stipulate henceforth that any labelling $\|\cdot\|$ of an $HLAF$ is defined over the extended domain $\mathsf{A}^{HL(0)}\cup\mathsf{R}^{HL(n)}$. To illustrate Definition \ref{defn17} and make it more accessible, we present the following example. 
\begin{exmp}\label{ex:2}
	Let $HLAF=(\mathsf{A}^{HL(0)}, \mathsf{R}^{HL(1)})$, $A^{HL(0)}=\{a,b\}$, $R^{HL(0)}=\{(a,b)\}$, and $R^{HL(1)}=\{(a,b), (b,(a,b))\}$. For brevity, we denote $(a,b)$ as $\alpha$ and denote $(b,(a,b))$ as $\beta$. Since $a$ and $\beta$ are not attacked by $a$ or $b$, we have $\mathsf{A}^{HL(0)}=\{a,b,\bot\}$ and $\mathsf{R}^{HL(1)}=\{\alpha,\beta,r_\bot^ a, r_\bot^ \beta\}$. Then, each element in $A^{HL(0)}\cup R^{HL(1)}$ has an attacker, and from Definition \ref{defn17}, we obtain the encoded formula $ec_{HL}(HLAF)$ of the given $HLAF$ via the normal encoding $ec_{HL}$ w.r.t. the $\mathcal{PL}_3^L$:
	\begin{equation*}
		(a\leftrightarrow \neg(\bot\wedge r_\bot^ a))\wedge(b\leftrightarrow \neg(a\wedge \alpha))\wedge(\alpha\leftrightarrow \neg(b\wedge \beta))\wedge(\beta\leftrightarrow \neg(\bot\wedge r_\bot^ \beta)).
	\end{equation*}
\end{exmp}
Based on the normal encoding of $HLAF$s, we have a theorem that for a given $HLAF$, an assignment $\|\cdot\|$ is a model of  the $HLAF$ under complete semantics iff it is a model of $ec_{HL}(HLAF)$ in the $PL_3^L$.
Denote that a labelling (alternatively called an assignment) $\|\cdot\|$ is a model of an $HLAF$ under complete semantics by $\|\cdot\| \models_{\mathcal{HL}_3^C} HLAF$. Then, we formally present the model equivalence theorem.
\begin{thm}\label{hsme}
	For an $HLAF=(\mathsf{A}^{HL(0)}, \mathsf{R}^{HL(n)})$ and an assignment $\|\cdot\|: \mathsf{A}^{HL(0)}\cup\mathsf{R}^{HL(n)}\to \{0,1,\frac{1}{2}\}$, 
	\begin{equation*}
		\|\cdot\| \models_{\mathcal{HL}_3^C} HLAF \Longleftrightarrow \|\cdot\| \models_{\mathcal{PL}_3^L} ec_{HL}(HLAF).
	\end{equation*}
\end{thm}
\begin{proof} 
	For the assignment $\|\cdot\|$, $\|\bot\|=0$ and $\|(\bot, \alpha)\|=1$ (if $(\bot, \alpha)\in \mathsf{R}^{HL(n)}$) hold in $\mathcal{PL}_3^L$ and also hold for the $HLAF$ under complete semantics by definitions.
	Thus, we need to check that for the given assignment $\|\cdot\|$ and $\forall \beta\in {A}^{HL}\cup{R}^{HL(n)}$, the value of $\beta$ (which can be 1 or 0 or $\frac{1}{2}$) satisfies complete semantics iff $\|\beta\leftrightarrow\bigwedge_{(a_i, \beta)\in \mathsf{R}^{HL(n)}}\neg(a_i\wedge r_{a_i}^{\beta})\|=1$ in the $\mathcal{PL}_3^L$. Then, we need to discuss three cases. 
	\begin{itemize}
		\item  Case 1, $\|\beta\|=1$. 
		\\$\|\beta\|=1$ by model $\|\cdot\|$ under complete semantics
		\\$\Longleftrightarrow$ for each $a_i$ such that $(a_i, \beta) \in \mathsf{R}^{HL(n)}$, we have either $\|a_i\|=0$ or $\|(a_i, \beta)\| = 0$
		\\$\Longleftrightarrow$ $\|(a_i\wedge r_{a_i}^{\beta})\|=0$ (i.e., $\|\neg(a_i\wedge r_{a_i}^{\beta})\|=1$) in the $\mathcal{PL}_3^L$
		\\$\Longleftrightarrow$ $\|\bigwedge_{(a_i, \beta)\in \mathsf{R}^{HL(n)}}\neg(a_i\wedge r_{a_i}^{\beta})\|=1$ in the $\mathcal{PL}_3^L$
		\\$\Longleftrightarrow$ $\|\beta\leftrightarrow\bigwedge_{(a_i, \beta)\in \mathsf{R}^{HL(n)}}\neg(a_i\wedge r_{a_i}^{\beta})\|=1$ in the $\mathcal{PL}_3^L$.
		
		\item  Case 2, $\|\beta\|=0$. 
		\\$\|\beta\|=0$ by model $\|\cdot\|$ under complete semantics
		\\$\Longleftrightarrow$ for some $ a_i$, $(a_i, \beta) \in \mathsf{R}^{HL(n)}$, we have $\|a_i\|= 1$ and $\|(a_i, \beta)\| =1$			
		\\$\Longleftrightarrow$ for some $ a_i$, $(a_i, \beta) \in \mathsf{R}^{HL(n)}$, we have $\|(a_i\wedge r_{a_i}^{\beta})\|=1$ (i.e., $\|\neg(a_i\wedge r_{a_i}^{\beta})\|=0$) in the $\mathcal{PL}_3^L$
		\\$\Longleftrightarrow$ $\|\bigwedge_{(a_i, \beta)\in \mathsf{R}^{HL(n)}}\neg(a_i\wedge r_{a_i}^{\beta})\|=0$ in the $\mathcal{PL}_3^L$
		\\$\Longleftrightarrow$ $\|\beta\leftrightarrow\bigwedge_{(a_i, \beta)\in \mathsf{R}^{HL(n)}}\neg(a_i\wedge r_{a_i}^{\beta})\|=1$ in the $\mathcal{PL}_3^L$.
		
		\item Case 3, $\|\beta\|=\frac{1}{2}$. 
			\\$\|\beta\|=\frac{1}{2}$ by model $\|\cdot\|$ under complete semantics
		\\$\Longleftrightarrow$ $\|\beta\|\neq1$ and $\|\beta\|\neq0$ under complete semantics
			\\$\Longleftrightarrow$ $\|\bigwedge_{(a_i, \beta)\in \mathsf{R}^{HL(n)}}\neg(a_i\wedge r_{a_i}^{\beta})\|\neq1$ and $\|\bigwedge_{(a_i, \beta)\in \mathsf{R}^{HL(n)}}\neg(a_i\wedge r_{a_i}^{\beta})\|\neq0$ in the $\mathcal{PL}_3^L$ by Case 1 and Case 2
		\\$\Longleftrightarrow$ $\|\bigwedge_{(a_i, \beta)\in \mathsf{R}^{HL(n)}}\neg(a_i\wedge r_{a_i}^{\beta})\|=\frac{1}{2}$ in the $\mathcal{PL}_3^L$
			\\$\Longleftrightarrow$ $\|\beta\leftrightarrow\bigwedge_{(a_i, \beta)\in \mathsf{R}^{HL(n)}}\neg(a_i\wedge r_{a_i}^{\beta})\|=1$ in the $\mathcal{PL}_3^L$.
	\end{itemize}
	
	From the three cases above, if the value of  any $\beta$ satisfies complete semantics, no matter whether $\|\beta\|=1$ or 0 or $\frac{1}{2}$, we have $\|\beta\leftrightarrow\bigwedge_{(a_i, \beta)\in \mathsf{R}^{HL(n)}}\neg(a_i\wedge r_{a_i}^{\beta})\|=1$ in the $\mathcal{PL}_3^L$ and vice versa by ``$\Longleftrightarrow$". 
	Thus, an assignment $\|\cdot\|$ of the $HLAF$ satisfies complete semantics iff $\|\bigwedge_{\beta\in A^{HL(0)}\cup R^{HL(n)}}(\beta\leftrightarrow\bigwedge_{(a_i, \beta)\in \mathsf{R}^{HL(n)}}\neg(a_i\wedge r_{a_i}^{\beta}))\|=1$, i.e., $\|\cdot\|$ is a model of the $ec_{HL}(HLAF)$, in the $\mathcal{PL}_3^L$.
	\end{proof} 
	We continue to use the $HLAF$ in Example \ref{ex:2} to illustrate this theorem.
\addtocounter{exmp}{-1}
\begin{exmp}[continued]\label{ex:2-cont}
	 For this given $HLAF$, we first find models under complete semantics. Since $a$ and $\beta$ are only attacked by $\bot$, $\|a\|=\|\beta\|=1$ by Definition \ref{newcom}. If $\|\alpha\|=1$ then $\|b\|=0$ by Definition \ref{newcom}, i.e., we have the first complete labelling $\{\|a\|=\|\beta\|=\|\alpha\|=\|r_\bot^ a\|=\|r_\bot^ \beta\|=1, \|b\|=\|\bot\|=0\}$. If $\|\alpha\|=0$ then $\|b\|=1$, i.e., we have the second complete labelling $\{\|a\|=\|\beta\|=\|b\|=\|r_\bot^ a\|=\|r_\bot^ \beta\|=1, \|\alpha\|=\|\bot\|=0\}$. If $\|\alpha\|=\frac{1}{2}$ then $\|b\|=\frac{1}{2}$, i.e., we have the third complete labelling $\{\|a\|=\|\beta\|=\|r_\bot^ a\|=\|r_\bot^ \beta\|=1, \|\alpha\|=\|b\|=\frac{1}{2}, \|\bot\|=0\}$.
	 
	 Next, we find models of the encoded formula $ec_{HL}(HLAF)$ in the $\mathcal{PL}_3^L$, i.e., all assignments that satisfy $\|ec_{HL}(HLAF)\|=1$. In the $\mathcal{PL}_3^L$ we have that\\
    	$\|ec_{HL}(HLAF)\|=1$\\
	 	$\Longleftrightarrow$ $\|(a\leftrightarrow \neg(\bot\wedge r_\bot^ a))\wedge(b\leftrightarrow \neg(a\wedge \alpha))\wedge(\alpha\leftrightarrow \neg(b\wedge \beta))\wedge(\beta\leftrightarrow \neg(\bot\wedge r_\bot^ \beta))\|=1$\\
	 	$\Longleftrightarrow$ 
	 	\begin{align*}
	 		1&=\|b\leftrightarrow \neg(a\wedge \alpha)\|\\
	 		&=\|\alpha\leftrightarrow \neg(b\wedge \beta)\|\\
	 		&=\|\beta\leftrightarrow \neg(\bot\wedge r_\bot^ \beta)\|\\
	 		&=\|a\leftrightarrow \neg(\bot\wedge r_\bot^ a)\|
	 	\end{align*}
	 	$\Longleftrightarrow$ 
	\begin{equation*}
		\begin{cases}
			\| \neg(\bot\wedge r_\bot^ a)\|=\|a\|\\
			\| \neg(a\wedge \alpha)\|=\|b\|\\
			\| \neg(b\wedge \beta)\|=\|\alpha\|\\
			\| \neg(\bot\wedge r_\bot^ \beta)\|=\|\beta\|
		\end{cases}.
	\end{equation*}
	The above system of equations can be simplified as follows
		\begin{equation*}
		\begin{cases}
		\|\bot\|=0\\
		\|r_\bot^ a\|=\|r_\bot^ \beta\|=1\\
		\|a\|=1\\
		\|\beta\|=1\\
		\|b\|+\|\alpha\|=1		
		\end{cases}.
	\end{equation*}
	Thus, in the $\mathcal{PL}_3^L$, we have three models $\{\|a\|=\|\beta\|=\|\alpha\|=\|r_\bot^ a\|=\|r_\bot^ \beta\|=1, \|b\|=\|\bot\|=0\}$, $\{\|a\|=\|\beta\|=\|b\|=\|r_\bot^ a\|=\|r_\bot^ \beta\|=1, \|\alpha\|=\|\bot\|=0\}$ and $\{\|a\|=\|\beta\|=\|r_\bot^ a\|=\|r_\bot^ \beta\|=1, \|\alpha\|=\|b\|=\frac{1}{2}, \|\bot\|=0\}$. They are exactly the same as models of the $HLAF$ under complete semantics.
\end{exmp}
\begin{cor}
	The $ec_{HL}$ is a translation of $\mathcal{HLAF}$ associated with $\mathcal{PL}_3^L$ under complete semantics and thus complete semantics of $\mathcal{HLAF}$ is translatable.
\end{cor}

\subsection{Three equational semantics of $HLAF$s}
In this subsection, we give three equational systems on [0,1] for $HLAF$s and then we prove the theorems of model equivalence which provide logical bases for those equational systems. The characteristics and applications of these semantics are not within the main scope of our research. Therefore, we will not conduct a detailed discussion here but only provide a few brief comments.

To simplify some proofs in this paper, we present two lemmas.
 \begin{lem}\label{remark1}
 	For any $m, n \in [0,1]$, $m=n$ iff $I_T(m,n)=I_T(n,m)=1$, where $I_T$ is an $R$-implication and $T$ is a continuous t-norm.
 \end{lem}
\begin{proof}
	 By Definition \ref{def11}, if $m=n$ then $I_T(m,n)=I_T(n,m)=1$ obviously. On the other hand, if $I_T(m,n)=sup\{p|T (m, p) \leq n\}=1$ and $T$ is  continuous, so $T(m, 1)\leq n$, i.e., $m\leq n$. Similarly, we can get $n\leq m$ from $I_T(n,m)=1$. Thus $m=n$. 	
\end{proof}

\begin{lem}\label{remark2}
	For formulas $a$ and $b$ in $\mathcal{PL}_{[0, 1]}s$, we have $\|a\|=\|b\|$ iff $\|a\leftrightarrow b\|=1$.
\end{lem}
\begin{proof}
	 If $\|a\|=\|b\|$ then $\|a\leftrightarrow b\|=T(I_T(\|a\|, \|b\|), I_T(\|b\|, \|a\|))=1$ by Lemma \ref{remark1}. On the other hand, if
	$\|a\leftrightarrow b\|=T(I_T(\|a\|, \|b\|), I_T(\|b\|, \|a\|))=1$ then $I_T(\|a\|, \|b\|)=I_T(\|b\|, \|a\|)=1$ and then $\|a\|=\| b\|$ by Lemma \ref{remark1}.
\end{proof}

In the following discussion, for a given $HLAF=(\mathsf{A}^{HL(0)}, \mathsf{R}^{HL(n)})$ and $\forall \beta\in A^{HL(0)}\cup R^{HL(n)}$, let $n_\beta = |\{a_i \mid (a_i,\beta) \in \mathsf{R}^{HL(n)}\}|$ denote the number of attackers of $\beta$. The normal encoding function of $HLAF$s w.r.t. a $\mathcal{PL}_{[0, 1]}$ is still $ec_{HL}$ such that	
\begin{equation*}
	ec_{HL}(HLAF)=\bigwedge_{\beta\in {A^{HL(0)}\cup R^{HL(n)}}}(\beta\leftrightarrow\bigwedge_{(a_i,\beta) \in \mathsf{R}^{HL(n)}}\neg(a_i\wedge r_{a_i}^{\beta})).
\end{equation*}
Note that here, with a slight abuse of notation (but without leading to confusion), we still use the symbol $ec_{HL}$ to denote the normal encoding function of $HLAF$s w.r.t. a $\mathcal{PL}_{[0,1]}$.

\subsubsection{The equational system $Eq_G^{HL}$}
For an $HLAF=(\mathsf{A}^{HL(0)}, \mathsf{R}^{HL(n)})$, we give the equational system $Eq_G^{HL}$ such that $\|\bot\|=0$, $\|r_{\bot}^\alpha\|=1$ (if $r_{\bot}^\alpha\in\mathsf{R}^{HL(n)}$), and the equation for each $\beta\in A^{HL(0)}\cup R^{HL(n)}$ is 
\begin{equation}\label{eq1}
	\|\beta\|=\min_{i=1}^{n_{\beta}}\max\{1-\|a_i\|, 1-\|r_{a_i}^{\beta}\|\}.
\end{equation}
This equational semantics $Eq_G^{HL}$ is well-founded as it models the strength of an attack as a joint requirement on both the attacker and the attack relation (via the $\max$ operator, computing the ``ineffectiveness" of a single attack), while the final value of the target is determined by its most potent threat (via the $\min$ operator). This formulation naturally extends classical argumentation semantics into a fuzzy setting: an argument is fully accepted only if all attacks are completely ineffective, and is fully defeated if any attack is fully effective, thereby providing a fine-grained representation of argument status on a continuous scale. 

Based on the normal encoding of $HLAF$s, we have a theorem that an assignment is a model of an $HLAF$ under equational semantics $Eq_G^{HL}$ iff it is a model of $ec_{HL}(HLAF)$ in the $PL_{[0, 1]}^G$. 
Denote that an assignment $\|\cdot\|$ is a model of an $HLAF$ under equational semantics $Eq_G^{HL}$ by $\|\cdot\| \models_{\mathcal{HL}_{[0,1]}^G} HLAF$. Then, we formally present the model equivalence theorem.
\begin{thm}\label{thm2}
	For an $HLAF=(\mathsf{A}^{HL(0)}, \mathsf{R}^{HL(n)})$ and an assignment $\|\cdot\|: \mathsf{A}^{HL(0)}\cup\mathsf{R}^{HL(n)}\to [0,1]$, 
	\begin{equation*}
		\|\cdot\| \models_{\mathcal{HL}_{[0,1]}^G} HLAF \Longleftrightarrow \|\cdot\| \models_{\mathcal{PL}_{[0, 1]}^G} ec_{HL}(HLAF).
	\end{equation*}
\end{thm}
\begin{proof} 
	A model $\|\cdot\|$ of $ec_{HL}(HLAF)$ in the $\mathcal{PL}_{[0, 1]}^G$
	\\$\Longleftrightarrow $ (a solution of) $\|\bigwedge_{\beta\in {A^{HL(0)}\cup R^{HL(n)}}}(\beta\leftrightarrow\bigwedge_{(a_i,\beta) \in \mathsf{R}^{HL(n)}}\neg(a_i\wedge r_{a_i}^{\beta}))\|=1$
	\\$\Longleftrightarrow $ for each $\beta\in {A^{HL(0)}\cup R^{HL(n)}}$, $\|\beta\leftrightarrow\bigwedge_{(a_i,\beta) \in \mathsf{R}^{HL(n)}}\neg(a_i\wedge r_{a_i}^{\beta})\|=1$
	\\$\Longleftrightarrow $ for each $\beta\in {A^{HL(0)}\cup R^{HL(n)}}$, $\|\beta\|=\|\bigwedge_{(a_i,\beta) \in \mathsf{R}^{HL(n)}}\neg(a_i\wedge r_{a_i}^{\beta})\|$ by Lemma \ref{remark2}
	\\$\Longleftrightarrow $ for each $\beta\in {A^{HL(0)}\cup R^{HL(n)}}$, $\|\beta\|=\min_{i=1}^{n_{\beta}}\|\neg(a_i\wedge r_{a_i}^{\beta})\|=\min_{i=1}^{n_{\beta}}(1-\|a_i\wedge r_{a_i}^{\beta}\|)=\min_{i=1}^{n_{\beta}}(1-\min\{\|a_i\|, \|r_{a_i}^{\beta}\|\})$
	\\$\Longleftrightarrow $ for each $\beta\in {A^{HL(0)}\cup R^{HL(n)}}$, $\|\beta\|=\min_{i=1}^{n_{\beta}}\max\{1-\|a_i\|, 1-\|r_{a_i}^{\beta}\|\}$	
	\\$\Longleftrightarrow$ a model $\|\cdot\|$ of the $HLAF$ under equational semantics $Eq_G^{HL}$.
\end{proof} 

\begin{cor}
The $ec_{HL}$ is a translation of $\mathcal{HLAF}$ associated with $\mathcal{PL}_{[0,1]}^G$ under the equational semantics $Eq_G^{HL}$.
\end{cor}
\begin{cor}
The equational semantics $Eq_G^{HL}$ is a fuzzy normal encoded semantics and thus it is translatable.
\end{cor}

The t-norm $T^G$ of $\mathcal{PL}_{[0,1]}^G$ characterizes the least truth degree and the t-conorm of $\mathcal{PL}_{[0,1]}^G$ characterizes the greatest truth degree. Thus the equational system $Eq_G^{HL}$ of an $HLAF$, which is model-equivalent to the encoded formula $ec_{HL}(HLAF)$ by the $\mathcal{PL}_{[0,1]}^G$, also characterizes the minimal and maximal truth degrees of all attackers and attacks of an argument. So when we mainly care about the minimality and the maximality of attackers and attacks we can choose this equational system $Eq_G^{HL}$ in practice. We continue Example \ref{ex:2} to illustrate Theorem \ref{thm2} below.
\addtocounter{exmp}{-1}
\begin{exmp}[continued]\label{ex:2-cont2}
	 According to Equation \ref{eq1}, the system of equations for the $HLAF$ is that	 
	 \begin{equation*}
	 	\begin{cases}	 		
	 		\|a\|=1\\
	 		\|b\|=\min\{\max\{1-\|a\|, 1-\|\alpha\|\}\}\\
	 		\|\alpha\|=\min\{\max\{1-\|b\|, 1-\|\beta\|\}\}\\
	 		\|\beta\|=1
	 	\end{cases},
	 \end{equation*}
	i.e.,
	\begin{equation}\label{eq2}
		\begin{cases}
			\|a\|=1\\
			\|\beta\|=1\\
			\|b\|+\|\alpha\|=1		
		\end{cases}.
	\end{equation}
	Thus, we have models $\{\|a\|=\|\beta\|=1, \|b\|=x, \|\alpha\|=1-x, \|\bot\|=0\}$ where $x\in [0,1]$.
	
	Next, we find models of $ec_{HL}(HLAF)$ in the $\mathcal{PL}_{[0, 1]}^G$, i.e., all assignments that satisfy $\|ec_{HL}(HLAF)\|=1$. In the $\mathcal{PL}_{[0, 1]}^G$ we have that\\
	$\|ec_{HL}(HLAF)\|=1$\\
	$\Longleftrightarrow$ $\|(a\leftrightarrow \neg(\bot\wedge r_\bot^ a))\wedge(b\leftrightarrow \neg(a\wedge \alpha))\wedge(\alpha\leftrightarrow \neg(b\wedge \beta))\wedge(\beta\leftrightarrow \neg(\bot\wedge r_\bot^ \beta))\|=1$\\
	$\Longleftrightarrow$ 
	\begin{equation*}
		\begin{cases}
			\| \neg(\bot\wedge r_\bot^ a)\|=\|a\|\\
			\| \neg(a\wedge \alpha)\|=\|b\|\\
			\| \neg(b\wedge \beta)\|=\|\alpha\|\\
			\| \neg(\bot\wedge r_\bot^ \beta)\|=\|\beta\|
		\end{cases}
	\end{equation*}
	$\Longleftrightarrow$ 
	\begin{equation*}
		\begin{cases}
			1-\min\{\|\bot\|, \|r_\bot^ a\|\}=\|a\|\\
			1-\min\{\|a\|, \|\alpha\|\}=\|b\|\\
			1-\min\{\|b\|, \|\beta\|\}=\|\alpha\|\\
			1-\min\{\|\bot\|, \|r_\bot^ \beta\|\}=\|\beta\|
		\end{cases}
	\end{equation*}
	$\Longleftrightarrow$ 
	\begin{equation*}
		\begin{cases}
			\|a\|=1\\
			\|\beta\|=1\\
			\|b\|+\|\alpha\|=1		
		\end{cases}
	\end{equation*}
	which is the same as Equation \ref{eq2}.
	Thus, the models of $ec_{HL}(HLAF)$ in the $\mathcal{PL}_{[0, 1]}^G$ are the same as models of the $HLAF$ under equational semantics $Eq_G^{HL}$. 
\end{exmp}

\subsubsection{The equational system $Eq_P^{HL}$}
For an $HLAF=(\mathsf{A}^{HL(0)}, \mathsf{R}^{HL(n)})$, we present the equational system $Eq_P^{HL}$ such that $\|\bot\|=0$, $\|r_{\bot}^\alpha\|=1$ (if $r_{\bot}^\alpha\in\mathsf{R}^{HL(n)}$), and the equation for each $\beta\in A^{HL(0)}\cup R^{HL(n)}$ is
\begin{equation}\label{eq3}
	\|\beta\|=\prod_{i=1}^{n_{\beta}}(1-\|a_i\|\|r_{a_i}^{\beta}\|).
\end{equation}
The equational semantics $Eq_P^{HL}$ is grounded in product logic and a probabilistic interpretation of independent events. Its rationale lies in treating each attack $(a_i, r_{a_i}^\beta)$ as an independent event that jointly diminishes the acceptability of $\beta$: $\|a_i\|\|r_{a_i}^\beta\|$ quantifies the ``joint effectiveness" of the attack, while $1 - \|a_i\|\|r_{a_i}^\beta\|$ represents its ``failure probability". Multiplying the failure probabilities of all attacks captures their cumulative attenuating effect under the assumption of independence—$\beta$ is fully accepted only if all attacks completely fail (all terms equal 1), and is entirely defeated if any attack is fully effective (the corresponding term becomes 0). This model naturally captures the synergistic yet independent impact of attacks within a continuous semantic framework.

Based on the normal encoding of $HLAF$s, we have a theorem that an assignment is a model of an $HLAF$ under equational semantics $Eq_P^{HL}$ iff it is a model of $ec_{HL}(HLAF)$ in the $PL_{[0, 1]}^P$. 
Denote that an assignment $\|\cdot\|$ is a model of an $HLAF$ under equational semantics $Eq_P^{HL}$ by $\|\cdot\| \models_{\mathcal{HL}_{[0,1]}^P} HLAF$. Then, we formally present the model equivalence theorem.
\begin{thm}\label{Thm3}
	For an $HLAF=(\mathsf{A}^{HL(0)}, \mathsf{R}^{HL(n)})$ and an assignment $\|\cdot\|: \mathsf{A}^{HL(0)}\cup\mathsf{R}^{HL(n)}\to [0,1]$, 
	\begin{equation*}
		\|\cdot\| \models_{\mathcal{HL}_{[0,1]}^P} HLAF \Longleftrightarrow \|\cdot\| \models_{\mathcal{PL}_{[0, 1]}^P} ec_{HL}(HLAF).
	\end{equation*}
\end{thm}
\begin{proof} 
	A model $\|\cdot\|$ of $ec_{HL}(HLAF)$ in the $\mathcal{PL}_{[0, 1]}^P$
	\\$\Longleftrightarrow $ $\|\bigwedge_{\beta\in {A^{HL(0)}\cup R^{HL(n)}}}(\beta\leftrightarrow\bigwedge_{(a_i,\beta) \in \mathsf{R}^{HL(n)}}\neg(a_i\wedge r_{a_i}^{\beta}))\|=1$
	\\$\Longleftrightarrow $ for each $\beta\in {A^{HL(0)}\cup R^{HL(n)}}$, $\|\beta\leftrightarrow\bigwedge_{(a_i,\beta) \in \mathsf{R}^{HL(n)}}\neg(a_i\wedge r_{a_i}^{\beta})\|=1$
	\\$\Longleftrightarrow $ for each $\beta\in {A^{HL(0)}\cup R^{HL(n)}}$, $\|\beta\|=\|\bigwedge_{(a_i,\beta) \in \mathsf{R}^{HL(n)}}\neg(a_i\wedge r_{a_i}^{\beta})\|$ by Lemma \ref{remark2}
	\\$\Longleftrightarrow $ for each $\beta\in {A^{HL(0)}\cup R^{HL(n)}}$, $\|\beta\|=\prod_{i=1}^{n_{\beta}}\|\neg(a_i\wedge r_{a_i}^{\beta})\|=\prod_{i=1}^{n_{\beta}}(1-\|a_i\wedge r_{a_i}^{\beta}\|)$
	\\$\Longleftrightarrow$ for each $\beta\in {A^{HL(0)}\cup R^{HL(n)}}$, $\|\beta\|=\prod_{i=1}^{n_{\beta}}(1-\|a_i\|\|r_{a_i}^{\beta}\|)$		
	\\$\Longleftrightarrow$ a model $\|\cdot\|$ of the $HLAF$ under equational semantics $Eq_P^{HL}$.
\end{proof} 

\begin{cor}
The $ec_{HL}$ is a translation of $\mathcal{HLAF}$ associated with $\mathcal{PL}_{[0,1]}^P$ under the equational semantics $Eq_P^{HL}$.
\end{cor}
\begin{cor}
	The equational semantics $Eq_P^{HL}$ is a fuzzy normal encoded semantics and thus it is translatable.
\end{cor}
The value obtained by applying the t-norm operation $T^P$ to a set of numbers in $\mathcal{PL}_{[0,1]}^P$ is the product of these numbers. When this feature is shown in the equational semantics $Eq_P^{HL}$ of $HLAF$s, it is the product form, which can characterize qualities and quantities of all attackers and attacks of an argument. So when we want to consider both factors of the qualities and the quantities of all attackers and attacks, we can choose this equational system $Eq_P^{HL}$ in practice. Next, we continue Example \ref{ex:2} to illustrate Theorem \ref{Thm3}.
\addtocounter{exmp}{-1}
\begin{exmp}[continued]\label{ex:2-cont3}
	According to Equation \ref{eq3}, the system of equations for the $HLAF$ is that	 
	\begin{equation*}
		\begin{cases}	 		
			\|a\|=1\\
			\|b\|=1-\|a\|\|\alpha\|\\
			\|\alpha\|=1-\|b\|\|\beta\|\\
			\|\beta\|=1
		\end{cases},
	\end{equation*}
	i.e.,
	\begin{equation}\label{eq4}
		\begin{cases}
			\|a\|=1\\
			\|\beta\|=1\\
			\|b\|+\|\alpha\|=1		
		\end{cases}.
	\end{equation}
	Thus, we have models $\{\|a\|=\|\beta\|=\|r_\bot^ a\|=\|r_\bot^ \beta\|=1, \|b\|=x, \|\alpha\|=1-x, \|\bot\|=0\}$ where $x\in [0,1]$.
	
	Next, we find models of $ec_{HL}(HLAF)$ in the $\mathcal{PL}_{[0, 1]}^P$, i.e., all assignments that satisfy $\|ec_{HL}(HLAF)\|=1$. In the $\mathcal{PL}_{[0, 1]}^P$ we have that\\
	$\|ec_{HL}(HLAF)\|=1$\\
	$\Longleftrightarrow$ $\|(a\leftrightarrow \neg(\bot\wedge r_\bot^ a))\wedge(b\leftrightarrow \neg(a\wedge \alpha))\wedge(\alpha\leftrightarrow \neg(b\wedge \beta))\wedge(\beta\leftrightarrow \neg(\bot\wedge r_\bot^ \beta))\|=1$\\
	$\Longleftrightarrow$ 
	\begin{equation*}
		\begin{cases}
			\| \neg(\bot\wedge r_\bot^ a)\|=\|a\|\\
			\| \neg(a\wedge \alpha)\|=\|b\|\\
			\| \neg(b\wedge \beta)\|=\|\alpha\|\\
			\| \neg(\bot\wedge r_\bot^ \beta)\|=\|\beta\|
		\end{cases}
	\end{equation*}
	$\Longleftrightarrow$ 
	\begin{equation*}
		\begin{cases}
			1-\|\bot\| \|r_\bot^ a\|=\|a\|\\
			1-\|a\| \|\alpha\|=\|b\|\\
			1-\|b\| \|\beta\|=\|\alpha\|\\
			1-\|\bot\| \|r_\bot^ \beta\|=\|\beta\|
		\end{cases}
	\end{equation*}
	$\Longleftrightarrow$ 
	\begin{equation*}
		\begin{cases}
			\|a\|=1\\
			\|\beta\|=1\\
			\|b\|+\|\alpha\|=1		
		\end{cases}
	\end{equation*}
	which is the same as Equation \ref{eq4}.
	Thus, the models of $ec_{HL}(HLAF)$ in the $\mathcal{PL}_{[0, 1]}^P$ are the same as models of the $HLAF$ under equational semantics $Eq_P^{HL}$. 
\end{exmp}
\subsubsection{The equational system $Eq_L^{HL}$}
For an $HLAF=(\mathsf{A}^{HL(0)}, \mathsf{R}^{HL(n)})$, we give the equational system $Eq_L^{HL}$ such that $\|\bot\|=0$, $\|r_{\bot}^\alpha\|=1$ (if $r_{\bot}^\alpha\in\mathsf{R}^{HL(n)}$), and the equation for each $\beta\in A^{HL(0)}\cup R^{HL(n)}$ is 
\begin{equation*}
	\|\beta\|=\begin{cases}
		0 & \quad \sum_{i=1}^{n_{\beta}}y_i\leqslant n_{\beta}-1\\
		\sum_{i=1}^{n_{\beta}}y_i-n_{\beta}+1 & \quad \sum_{i=1}^{n_{\beta}}y_i> n_{\beta}-1
	\end{cases}
\end{equation*}
where 
\begin{equation*}
	y_i=\begin{cases}
		1 & \quad \|a_i\|+\|r_{a_i}^{\beta}\|\leqslant 1\\
		2-\|a_i\|-\|r_{a_i}^{\beta}\| & \quad \|a_i\|+\|r_{a_i}^{\beta}\|> 1
	\end{cases}.
\end{equation*}
This equational semantics is grounded in {\L}ukasiewicz logic, and its core characteristic lies in using a threshold mechanism and a linear compensation model to capture the synergistic effect of attacks. Here, $y_i$ quantifies the degree of ineffectiveness of each attack: when the sum of the strengths of the attacker and the attack relation $|a_i|+|r_{a_i}^{\beta}|\leqslant 1$, the attack is fully ineffective ($y_i=1$); otherwise, its ineffectiveness decreases linearly as the joint strength increases. The final truth value of the target argument $\beta$ is determined by the sum of the ineffectiveness degrees of all its attacks—$\beta$ obtains a non-zero truth value only if this sum exceeds the threshold $n_\beta-1$, and its value then increases linearly with the difference. This design not only embodies the collective rationality that ``an argument is accepted only if sufficiently many of its attacks fail" but also enables a smooth, continuous progression for the argument's status from complete rejection to partial acceptance. It thereby extends the expressive power of classical semantics while maintaining computational interpretability.

Based on the normal encoding of $HLAF$s, we have a theorem that an assignment is a model of an $HLAF$ under equational semantics $Eq_L^{HL}$ iff it is a model of $ec_{HL}(HLAF)$ in the $PL_{[0, 1]}^L$. 
Denote that an assignment $\|\cdot\|$ is a model of an $HLAF$ under equational semantics $Eq_L^{HL}$ by $\|\cdot\| \models_{\mathcal{HL}_{[0,1]}^L} HLAF$. Then, we formally present the model equivalence theorem.
\begin{thm}\label{Thm4}
	For an $HLAF=(\mathsf{A}^{HL(0)}, \mathsf{R}^{HL(n)})$ and an assignment $\|\cdot\|: \mathsf{A}^{HL(0)}\cup\mathsf{R}^{HL(n)}\to [0,1]$, 
	\begin{equation*}
		\|\cdot\| \models_{\mathcal{HL}_{[0,1]}^L} HLAF \Longleftrightarrow \|\cdot\| \models_{\mathcal{PL}_{[0, 1]}^L} ec_{HL}(HLAF).
	\end{equation*}
\end{thm}

\begin{proof} 
	The following formula for the $n$-ary {\L}ukasiewicz t-norm is presented in \cite{Beliakov2007Aggregation}, as also independently shown in \cite{tang2025encoding}. For every $i \in {1,2,\dots,n}$, given $x_i \in [0,1]$:
	\begin{equation}\label{n-t-norm}
		x_1\ast x_2\ast\dots\ast x_n=
		\begin{cases}
			0 & \quad \sum_{i=1}^{n}x_i \leqslant n-1\\
			\sum_{i=1}^{n}x_i-n+1 & \quad \sum_{i=1}^{n}x_i > n-1
		\end{cases}
	\end{equation}
	where $\ast$ is the t-norm of the $\mathcal{PL}_{[0, 1]}^L$.
	
	Now let us prove this theorem with the use of the above formula \ref{n-t-norm}.
	\\A model $\|\cdot\|$ of $ec_{HL}(HLAF)$ in $\mathcal{PL}_{[0, 1]}^L$
	\\$\Longleftrightarrow $ $\|\bigwedge_{\beta\in {A^{HL(0)}\cup R^{HL(n)}}}(\beta\leftrightarrow\bigwedge_{(a_i,\beta) \in \mathsf{R}^{HL(n)}}\neg(a_i\wedge r_{a_i}^{\beta}))\|=1$
	\\$\Longleftrightarrow $ for each $\beta\in {A^{HL(0)}\cup R^{HL(n)}}$, $\|\beta\leftrightarrow\bigwedge_{(a_i,\beta) \in \mathsf{R}^{HL(n)}}\neg(a_i\wedge r_{a_i}^{\beta})\|=1$
	\\$\Longleftrightarrow $ for each $\beta$, $\|\beta\|=\|\bigwedge_{(a_i,\beta) \in \mathsf{R}^{HL(n)}}\neg(a_i\wedge r_{a_i}^{\beta})\|$ by Lemma \ref{remark2}
	\\$\Longleftrightarrow $ for each $\beta$, $\|\beta\|=\begin{cases}
		0 & \quad \sum_{i=1}^{n_{\beta}}y_i\leqslant n_{\beta}-1\\
		\sum_{i=1}^{n_{\beta}}y_i-n_{\beta}+1 & \quad \sum_{i=1}^{n_{\beta}}y_i> n_{\beta}-1
	\end{cases}$,
	where $y_i=\|\neg(a_i\wedge r_{a_i}^{\beta})\|=1-\|a_i\wedge r_{a_i}^{\beta}\|$.
	
	By the {\L}ukasiewicz t-norm $T_L (x, y)= \max\{0, x +y-1\}$, we have 
	\begin{equation*}
		\|a_i\wedge r_{a_i}^{\beta}\|=\begin{cases}
			0 & \quad \|a_i\|+\| r_{a_i}^{\beta}\|\leqslant 1\\
			\|a_i\|+\|r_{a_i}^{\beta}\|-1 & \quad \|a_i\|+\|r_{a_i}^{\beta}\| > 1
		\end{cases}.
	\end{equation*}
	Therefore, we obtain:
	\begin{equation}
		y_i=
		\begin{cases}
			1 & \quad \|a_i\|+\|r_{a_i}^{\beta}\|\leqslant 1\\
			2-\|a_i\|-\|r_{a_i}^{\beta}\| & \quad \|a_i\|+\|r_{a_i}^{\beta}\|> 1
		\end{cases}.
	\end{equation}
	
	The equations above are the same as the equational system $Eq_L^{HL}$. So $\|\cdot\|$ is a model of $ec_{HL}(HLAF)$ in the $\mathcal{PL}_{[0, 1]}^L$ iff it is
	a model of the $HLAF$ under equational semantics $Eq_L^{HL}$, i.e., we have proven this theorem.
\end{proof} 

\begin{cor}
The	$ec_{HL}$ is a translation of $\mathcal{HLAF}$ associated with $\mathcal{PL}_{[0,1]}^L$ under the equational semantics $Eq_L^{HL}$.
\end{cor}
\begin{cor}
	The equational semantics $Eq_L^{HL}$ is a fuzzy normal encoded semantics and thus it is translatable.
\end{cor}
The t-norm operation $T^L$ in $\mathcal{PL}_{[0,1]}^L$ characterizes the summation of a set of numbers with respect to a threshold. 
This characteristic of t-norm $T^L$ is reflected in $HLAF$s, where the truth value of an argument is determined by a threshold-limited sum of values of its attackers and attacks. Thereby $Eq_L^{HL}$ can also consider the qualities and quantities of attackers and attacks. 
The distinction compared to $Eq_P^{HL}$ consists in two aspects: (1) the presence/absence of threshold constraints, and (2) the choice between additive and multiplicative aggregation operators.
So if we want to characterize the truth value of an argument by the summation of the values of attackers and attacks with the threshold, we can choose this equational system $Eq_L^{HL}$ in practice. Additionally, the illustration of Theorem \ref{Thm4} using continued Example \ref{ex:2} is straightforward and will thus be omitted here.

\subsubsection{A comparative analysis of the three equational semantics}
In this part, we present a comparative analysis of the features and summarize the suitable application scenarios of the three equational semantics for $HLAF$s.
The three equational semantics---G\"{o}del-based ($Eq_G^{HL}$), Product-based ($Eq_P^{HL}$), and {\L}ukasiewicz-based ($Eq_L^{HL}$)---offer distinct approaches to modeling argument acceptability in $HLAF$s, each with unique characteristics and practical implications.

\begin{enumerate}
	\item \textbf{G\"{o}del Semantics ($Eq_G^{HL}$)}
	\begin{itemize}
		\item \textbf{Characteristics}: Employs min-max operations to determine argument strength based on the most potent threat
		\item \textbf{Strengths}: Provides strict, conservative evaluations where arguments are only fully accepted if all attacks are completely ineffective
		\item \textbf{Weaknesses}: May be overly pessimistic in scenarios with multiple moderate-strength attacks
	\end{itemize}
	
	\item \textbf{Product Semantics ($Eq_P^{HL}$)}
	\begin{itemize}
		\item \textbf{Characteristics}: Uses multiplicative operations to model cumulative weakening effects under the assumption of independent attacks
		\item \textbf{Strengths}: Naturally captures probabilistic intuition and gradual degradation from multiple simultaneous attacks
		\item \textbf{Weaknesses}: The independence assumption may not hold in highly correlated attack scenarios
	\end{itemize}
	
	\item \textbf{{\L}ukasiewicz Semantics ($Eq_L^{HL}$)}
	\begin{itemize}
		\item \textbf{Characteristics}: Implements threshold-based mechanisms with linear compensation to model collective attack effects
		\item \textbf{Strengths}: Enables smooth transitions between argument states and captures synergistic attack interactions
		\item \textbf{Weaknesses}: More computationally complex due to piecewise linear functions and threshold calculations
	\end{itemize}
\end{enumerate}

The choice among these semantics should be guided by the intrinsic characteristics of the argumentation process being modeled:
\begin{itemize}
	\item \textbf{Select G\"{o}del Semantics when}: the defeat of an argument is determined by its \emph{single most powerful attacker}. This is suitable for contexts where an argument can be completely undermined by one compelling counter-argument, and the presence of additional, weaker attackers does not further diminish its status. It models a classic, \emph{qualitative} view of argumentation.
	
	\item \textbf{Choose Product Semantics when}: multiple attackers contribute \emph{independently} to the cumulative weakening of an argument. This is appropriate when each attacker represents a distinct, separable flaw or piece of evidence, and their combined effect leads to a gradual degradation of the argument's acceptability.
	
	\item \textbf{Prefer {\L}ukasiewicz Semantics when}: attackers exhibit a \emph{synergistic or collective effect}. This is key for modeling scenarios where a \emph{critical mass} of attackers is required to defeat an argument, and the state of the argument changes smoothly with the aggregate strength of the attacking coalition, capturing a more \emph{quantitative} and collective form of reasoning.
\end{itemize}
Additionally, computational efficiency considerations may influence the choice: G\"{o}del semantics typically offers the simplest computation, followed by Product semantics, with {\L}ukasiewicz semantics being the most computationally demanding due to its threshold mechanisms.

Ultimately, the selection should align with both the epistemological assumptions of the domain (how attacks interact) and the practical constraints of the application environment.
In subsequent corresponding sections, the three equations based on these logics all exhibit the aforementioned features and are subject to the same selection considerations. Therefore, in the following chapters, we will forgo repetitive elaboration on these points.
\section{Encoding Barringer-Gabbay-Woods higher-order argumentation frameworks}\label{sec4}
In this section, we first encode $BHAF$s under complete semantics to the $\mathcal{PL}_3^L$.  Then, we give three typical equational semantics of $BHAF$s and encode $BHAF$s to $\mathcal{PL}_{[0,1]}$s under these three equational semantics. It is shown that these semantics are translatable by proving the model equivalence. 
\subsection{The $BHAF$ and its complete semantics}\label{sub4.1}
In the paper \cite{barringer2005temporal}, the authors construct a kind of $AF$ called ``attack network" where both arguments and attacks can attack an argument or an attack. We refer to it as Barringer-Gabbay-Woods higher-order argumentation frameworks ($BHAF$s) and define it in a way similar to that of $HLAF$s, where the level can be clearly highlighted.

\begin{defn}\label{bhdefn}
A \emph{0-level preparatory $BHAF$} ($pre\text{-}BHAF^{(0)}$) is a pair $(A^{BH(0)}, R^{BH(0)})$, where $R^{BH(0)}\subseteq A^{BH(0)} \times A^{BH(0)}$. An \emph{$n$-level preparatory $BHAF$} ($pre\text{-}BHAF^{(n)}$) is a pair $(A^{BH(0)}, R^{BH(n)})$, where $n\geqslant1$, $A^{BH(n)}=A^{BH(0)}\cup A^{BH(n-1)} \times A^{BH(n-1)}$, $R^{BH(n)}\subseteq A^{BH(n)} \times A^{BH(n)}$, and $R^{BH(n)}\nsubseteq A^{BH(n-1)} \times A^{BH(n-1)}$. 

An \emph{$n$-level $BHAF$} ($BHAF^{(n)}$)  is a pair $(\mathsf{A}^{BH(0)}, \mathsf{R}^{BH(n)})$, where $n\geqslant0$, $\mathsf{A}^{BH(0)}=A^{BH(0)}\cup\{\bot\}$, $\mathsf{R}^{BH(n)}=R^{BH(n)}\cup \{(\bot, \beta)\mid \beta\in A^{BH(0)}\cup R^{BH(n)} \text{ and } \nexists \alpha\in A^{BH(0)}\cup R^{BH(n)}: (\alpha, \beta)\in R^{BH(n)}\}$. 

Specifically, in a $BHAF^{(n)}$, $\mathsf{A}^{BH(0)}$ is the finite set of arguments including the imaginary argument $\bot$, and $\mathsf{R}^{BH(n)}$ is the $n$-level attack relation that may include some imaginary attacks $(\bot, \beta)$, where $\beta\in {A}^{BH(0)}\cup {R}^{BH(n)}$.
For any assignment $\|\cdot\|$ of a $BHAF^{(n)}$, let $\|\bot\|=0$ and $\|(\bot, \beta)\|=1$ (if $(\{\bot\}, \beta)\in \mathsf{R}^{BH(n)}$).
\end{defn}

We explain the intuition of this definition.
A $pre\text{-}BHAF^{(0)}$ is a $DAF$. $A^{BH(1)}=A^{BH(0)}\cup (A^{BH(0)} \times A^{BH(0)})$ includes all arguments and all possible 0-level attacks, where we can regard all arguments and all possible 0-level attacks at this level as arguments to construct the relation at the next level.
Thus, $A^{BH(1)} \times A^{BH(1)}$ includes all possible 0-level attacks, all possible attacks from arguments to 0-level attacks, all possible attacks from 0-level attacks to arguments, and all possible attacks from 0-level attacks to 0-level attacks. Note that it includes not only higher-level attacks but also 0-level attacks.
$A^{BH(2)}=A^{BH(0)}\cup (A^{BH(1)} \times A^{BH(1)})$ includes all arguments and all possible attacks under this level, and we can regard all arguments and all possible attacks as arguments to construct the next level. As this process continues, $A^{BH(n)}=A^{BH(0)}\cup (A^{BH(n-1)}\times A^{BH(n-1)})$ includes all arguments and all possible attacks under this level (including lower-level attacks). $R^{BH(n)}\subseteq A^{BH(n)} \times A^{BH(n)}$ and $R^{BH(n)}\nsubseteq A^{BH(n-1)} \times A^{BH(n-1)}$ mean that the $n$-level attack relation indeed includes some higher-level attacks and not just the attacks in the $(n-1)$-level attack relation. Following the definition of a $pre\text{-}BHAF^{(n)}$, if we add the imaginary argument $\bot$ and the imaginary attack $(\bot, \beta)$ for each $\beta$ that is not attacked in the $pre\text{-}BHAF^{(n)}$, then we obtain the $BHAF^{(n)}$. Intuitively, adding an invalid attack does not affect the truth value of the target argument (or the target attack). Note that as can be inferred from the definition above, for a $BHAF^{(n)}$, the imaginary attack $(\bot, \beta)$ neither attacks any argument or attack nor is attacked by any argument or attack. Briefly, we use the notation $BHAF$ instead of $BHAF^{(n)}$ when the level $n$ does not need to be stressed. In addition, a $BHAF$ in which every attacker is an argument is, by definition, an $HLAF$. Thus, the $BHAF$ can be seen as a generalization of the $HLAF$.

In the paper \cite{barringer2005temporal}, the authors give the definition of the valuation function and use the iterative (gradual) semantics. Next, we give the complete semantics of $BHAF$s which is generalized from the complete semantics of $HLAF$s. 

\begin{defn}\label{com-bhaf}
	A \emph{complete labelling} for a $BHAF=(\mathsf{A}^{BH(0)}, \mathsf{R}^{BH(n)})$ is a function
	$\|\cdot\| : \mathsf{A}^{BH(0)}\cup \mathsf{R}^{BH(n)}\rightarrow\{0, 1, \frac{1}{2}\}$ such that $\|\bot\|=0$, $\|(\bot, \gamma)\|=1$ (if $(\bot, \gamma)\in \mathsf{R}^{BH(n)}$) and for each $\beta\in A^{BH(0)}\cup R^{BH(n)}$
	\begin{equation*}
		\|\beta\| =\begin{cases}
			1 & \quad \text{ iff } \forall \alpha ((\alpha, \beta)\in \mathsf{R}^{BH(n)})): \|\alpha\|=0 \text{ or } \|(\alpha, \beta)\|= 0\\
			0 & \quad \text{ iff } \exists \alpha ((\alpha, \beta)\in \mathsf{R}^{BH(n)})): \|\alpha\|= 1 \text{ and } \|(\alpha, \beta)\|=1\\
			\frac{1}{2} & \quad \text{ otherwise}
		\end{cases}.
	\end{equation*}
\end{defn}

Based on the syntactic structure of $BHAF$s, the complete semantics of $BHAF$s fully inherit the features of $HLAF$s: an attacking action is considered valid only when both the attacker and its corresponding attack are valid; it is invalid if either the attacker or the attack is invalid; all other cases result in a semi-valid attacking action. The key distinction lies in the fact that, in a $BHAF$, the attacker can itself be an attack. Note that for any assignment $\|\cdot\|$ in any $\mathcal{PLS}$, we always let $\|\bot\|=0$ and $\|(\bot,\gamma)\|=1$ (if $(\bot,\gamma)\in\mathsf{R}^{BH(n)}$ for a $BHAF=(\mathsf{A}^{BH(0)}, \mathsf{R}^{BH(n)})$). 
\subsection{Encoding $BHAF$s under complete semantics}
We first present the encoding function below.
\begin{defn}
For a given $BHAF=(\mathsf{A}^{BH(0)}, \mathsf{R}^{BH(n)})$, the \emph{normal encoding function} of the $BHAF$ w.r.t. the $\mathcal{PL}_3^L$ is $ec_{BH}$ such that
\begin{equation*}
	ec_{BH}(BHAF)=\bigwedge_{\beta\in A^{BH(0)}\cup R^{BH(n)}}(\beta\leftrightarrow\bigwedge_{(\alpha_i, \beta)\in \mathsf{R}^{BH(n)}}\neg(\alpha_i\wedge r_{\alpha_i}^{\beta}))
\end{equation*}
 where $r_{\alpha_i}^{\beta}$ denotes the attack $(\alpha_i, \beta)$.
\end{defn}
The encoded formula $ec_{BH}(BHAF)$ is structurally similar to $ec_{HL}(HLAF)$. The key difference is that in $ec_{BH}(BHAF)$, the element $\alpha_i$ can be an attack, a direct consequence of the differing syntactic structures of the two frameworks. 

Denote that an assignment $\|\cdot\|$ is a model of a $BHAF$ under complete semantics by $\|\cdot\| \models_{\mathcal{BH}_3^C} BHAF$. Then, we present the model equivalence theorem.
\begin{thm}\label{bhcs}
	For a $BHAF=(\mathsf{A}^{BH(0)}, \mathsf{R}^{BH(n)})$ and an assignment $\|\cdot\|: \mathsf{A}^{BH(0)}\cup\mathsf{R}^{BH(n)}\to \{0,1,\frac{1}{2}\}$, 
	\begin{equation*}
		\|\cdot\| \models_{\mathcal{BH}_3^C} BHAF \Longleftrightarrow \|\cdot\| \models_{\mathcal{PL}_3^L} ec_{BH}(BHAF).
	\end{equation*}
\end{thm}

\begin{proof} 
	We need to check that for a given assignment $\|\cdot\|$ and $\forall \beta\in A^{BH(0)}\cup R^{BH(n)}$, any value of $\beta$ satisfies complete semantics iff
	\begin{equation}\label{eqbhbeta}
		\|\beta\leftrightarrow\bigwedge_{(\alpha_i, \beta)\in \mathsf{R}^{BH(n)}}\neg(\alpha_i\wedge r_{\alpha_i}^{\beta})\|=1
	\end{equation}
    holds in the $\mathcal{PL}_3^L$.  We need to discuss three cases.
	\begin{itemize}
		\item  Case 1, $\|\beta\|=1$. 
		\\$\|\beta\|=1$ by model $\|\cdot\|$ under complete semantics
			\\$\Longleftrightarrow$ for all $\alpha_i$ such that $(\alpha_i, \beta) \in \mathsf{R}^{BH(n)}$ we have that either $\|\alpha_i\|=0$ or $\|(\alpha_i, \beta)\|= 0$
		\\$\Longleftrightarrow$ for all $\alpha_i$ such that $(\alpha_i, \beta) \in \mathsf{R}^{BH(n)}$ we have $\|(\alpha_i\wedge r_{\alpha_i}^{\beta})\|=0$ (i.e., $\|\neg(\alpha_i\wedge r_{\alpha_i}^{\beta})\|=1$) in the $\mathcal{PL}_3^L$
		\\$\Longleftrightarrow$ $\|\bigwedge_{(\alpha_i, \beta)\in \mathsf{R}^{BH(n)}}\neg(\alpha_i\wedge r_{\alpha_i}^{\beta})\|=1$  in the $\mathcal{PL}_3^L$
		\\$\Longleftrightarrow$ Equation \ref{eqbhbeta} holds in the $\mathcal{PL}_3^L$.

		\item  Case 2, $\|\beta\|=0$. 
		\\$\|\beta\|=0$ by model $\|\cdot\|$ under complete semantics
		\\$\Longleftrightarrow$ for some $ \alpha_i$, $(\alpha_i, \beta) \in \mathsf{R}^{BH(n)}$, we have $\|\alpha_i\|= 1$ and $\|(\alpha_i, \beta)\| =1$			
		\\$\Longleftrightarrow$ for some $\alpha_i$, $(\alpha_i, \beta) \in \mathsf{R}^{BH(n)}$, we have $\|(\alpha_i\wedge r_{\alpha_i}^{\beta})\|=1$ (i.e., $\|\neg(\alpha_i\wedge r_{\alpha_i}^{\beta})\|=0$) in the $\mathcal{PL}_3^L$
		\\$\Longleftrightarrow$ $\|\bigwedge_{(\alpha_i, \beta)\in \mathsf{R}^{BH(n)}}\neg(\alpha_i\wedge r_{\alpha_i}^{\beta})\|=0$  in the $\mathcal{PL}_3^L$
		\\$\Longleftrightarrow$ Equation \ref{eqbhbeta} holds in the $\mathcal{PL}_3^L$.
		
		\item Case 3, $\|\beta\|=\frac{1}{2}$. 
		\\$\|\beta\|=\frac{1}{2}$ by model $\|\cdot\|$ under complete semantics
		\\$\Longleftrightarrow$ $\|\beta\|\neq1$ and $\|\beta\|\neq0$ under complete semantics
		\\$\Longleftrightarrow$ $\|\bigwedge_{(\alpha_i, \beta)\in \mathsf{R}^{BH(n)}}\neg(\alpha_i\wedge r_{\alpha_i}^{\beta})\|\neq1$ and $\|\bigwedge_{(\alpha_i, \beta)\in \mathsf{R}^{BH(n)}}\neg(\alpha_i\wedge r_{\alpha_i}^{\beta})\|\neq0$ in the $\mathcal{PL}_3^L$ by Case 1 and Case 2
		\\$\Longleftrightarrow$ $\|\bigwedge_{(\alpha_i, \beta)\in \mathsf{R}^{BH(n)}}\neg(\alpha_i\wedge r_{\alpha_i}^{\beta})\|=\frac{1}{2}$ in the $\mathcal{PL}_3^L$
		\\$\Longleftrightarrow$ Equation \ref{eqbhbeta} holds in the $\mathcal{PL}_3^L$.
	\end{itemize}
	
	From the three cases above, if the value of  any $\beta$ satisfies complete semantics, no matter whether $\|\beta\|=1$ or 0 or $\frac{1}{2}$, Equation \ref{eqbhbeta} holds in the $\mathcal{PL}_3^L$ and vice versa by ``$\Longleftrightarrow$". 
	So, an assignment $\|\cdot\|$ of the $BHAF$ satisfies complete semantics iff $\|\bigwedge_{\beta\in A^{BH(0)}\cup R^{BH(n)}}(\beta\leftrightarrow\bigwedge_{(\alpha_i, \beta)\in \mathsf{R}^{BH(n)}}\neg(\alpha_i\wedge r_{\alpha_i}^{\beta}))\|=1$, i.e., $\|\cdot\|$ is a model of the $ec_{BH}(BHAF)$, in the $\mathcal{PL}_3^L$.
	\end{proof}
\begin{cor}
The	$ec_{BH}$ is a translation of $\mathcal{BHAF}$ associated with $\mathcal{PL}_3^L$ under complete semantics and thus complete semantics of $\mathcal{BHAF}$ is translatable.
\end{cor}
\subsection{Three equational semantics of $BHAF$s}
In this subsection, we emphasize three equational systems on [0,1] for $BHAF$s and then we prove the theorems of model equivalence. In the following discussion, for a given $BHAF=(\mathsf{A}^{BH(0)}, \mathsf{R}^{BH(n)})$ and $\forall\beta\in A^{BH(0)}\cup R^{BH(n)}$, let $n_\beta = |\{\alpha_i \mid (\alpha_i,\beta) \in \mathsf{R}^{BH(n)}\}|$ denote the number of attackers of $\beta$. The encoding function of the $BHAF$ w.r.t. a $\mathcal{PL}_{[0, 1]}$ is still $ec_{BH}$ such that	
\begin{equation*}
	ec_{BH}(BHAF)=\bigwedge_{\beta\in A^{BH(0)}\cup R^{BH(n)}}(\beta\leftrightarrow\bigwedge_{(\alpha_i, \beta)\in \mathsf{R}^{BH(n)}}\neg(\alpha_i\wedge r_{\alpha_i}^{\beta}))
\end{equation*}
where $r_{\alpha_i}^{\beta}$ denotes the attack $(\alpha_i, \beta)$.
 
 In this subsection, based on the encoding $ec_{BH}$, we present three typical equational semantics of $BHAF$s and show three corresponding theorems of model equivalence.
The proofs of these theorems are extremely similar to proofs in the cases of $HLAF$s. In these proofs, we only need to change $a_i$ (an argument) as $\alpha_i$ (an argument or an attack) and let $(\alpha_i, \beta)$ belong to $\mathsf{R}^{BH(n)}$. So we give one proof and omit other proofs. 

\subsubsection{The equational system $Eq_G^{BH}$}
For a $BHAF=(\mathsf{A}^{BH(0)}, \mathsf{R}^{BH(n)})$, we give the equational system $Eq_G^{BH}$ such that $\|\bot\|=0$, $\|r_{\bot}^\alpha\|=1$ (if $r_{\bot}^\alpha\in\mathsf{R}^{BH(n)}$), and the equation for each $\beta\in A^{BH(0)}\cup R^{BH(n)}$ is
\begin{equation}\label{bhes}
	\|\beta\|=\min_{i=1}^{n_{\beta}}\max\{1-\|\alpha_i\|, 1-\|r_{\alpha_i}^{\beta}\|\}.
\end{equation}
This equation semantics signifies that the truth value of an argument or attack is determined by the weakest ineffectiveness among all its attackers, where the ineffectiveness of each attacker is calculated based on the truth values of both the attacker itself and the corresponding attack relation (i.e., the strength of either component can render the attack effective).

Denote that an assignment $\|\cdot\|$ is a model of a $BHAF$ under equational semantics $Eq_G^{BH}$ by $\|\cdot\| \models_{\mathcal{BH}_{[0,1]}^G} BHAF$. Then, we present the model equivalence theorem.
\begin{thm}\label{bheg}
	For a $BHAF=(\mathsf{A}^{BH(0)}, \mathsf{R}^{BH(n)})$ and an assignment $\|\cdot\|: \mathsf{A}^{BH(0)}\cup\mathsf{R}^{BH(n)}\to [0,1]$, 
	\begin{equation*}
		\|\cdot\| \models_{\mathcal{BH}_{[0,1]}^G} BHAF \Longleftrightarrow \|\cdot\| \models_{\mathcal{PL}_{[0,1]}^G} ec_{BH}(BHAF).
	\end{equation*}
\end{thm}

\begin{proof} 
	A model $\|\cdot\|$ of $ec_{BH}(BHAF)$ in the $\mathcal{PL}_{[0, 1]}^G$
	\\$\Longleftrightarrow $ (a solution of) $\|\bigwedge_{\beta\in A^{BH(0)}\cup R^{BH(n)}}(\beta\leftrightarrow\bigwedge_{(\alpha_i, \beta)\in \mathsf{R}^{BH(n)}}\neg(\alpha_i\wedge r_{\alpha_i}^{\beta}))\|=1$
	\\$\Longleftrightarrow $ for each $\beta\in A^{BH(0)}\cup R^{BH(n)}$, $\|\beta\leftrightarrow\bigwedge_{(\alpha_i, \beta)\in \mathsf{R}^{BH(n)}}\neg(\alpha_i\wedge r_{\alpha_i}^{\beta})\|=1$
	\\$\Longleftrightarrow $ for each $\beta\in A^{BH(0)}\cup R^{BH(n)}$, $\|\beta\|=\|\bigwedge_{(\alpha_i, \beta)\in \mathsf{R}^{BH(n)}}\neg(\alpha_i\wedge r_{\alpha_i}^{\beta})\|$ by Lemma \ref{remark2}
	\\$\Longleftrightarrow $ for each $\beta\in A^{BH(0)}\cup R^{BH(n)}$, $\|\beta\|=\min_{i=1}^{n_{\beta}}\|\neg(\alpha_i\wedge r_{\alpha_i}^{\beta})\|=\min_{i=1}^{n_{\beta}}(1-\|\alpha_i\wedge r_{\alpha_i}^{\beta}\|)=\min_{i=1}^{n_{\beta}}(1-\min\{\|\alpha_i\|, \|r_{\alpha_i}^{\beta}\|\})$
	\\$\Longleftrightarrow $ for each $\beta\in A^{BH(0)}\cup R^{BH(n)}$, $\|\beta\|=\min_{i=1}^{n_{\beta}}\max\{1-\|\alpha_i\|, 1-\|r_{\alpha_i}^{\beta}\|\}$	
	\\$\Longleftrightarrow$ a model $\|\cdot\|$ of the $BHAF$ under equational semantics $Eq_G^{BH}$.
\end{proof} 

\begin{cor}
The	$ec_{BH}$ is a translation of $\mathcal{BHAF}$ associated with $\mathcal{PL}_{[0,1]}^G$ under the equational semantics $Eq_G^{BH}$.
\end{cor}
\begin{cor}
	The equational semantics $Eq_G^{BH}$ is a fuzzy normal encoded semantics and thus it is translatable.
\end{cor}
\subsubsection{The equational system $Eq_P^{BH}$}
For a $BHAF=(\mathsf{A}^{BH(0)}, \mathsf{R}^{BH(n)})$, we give the equational system $Eq_P^{BH}$ such that $\|\bot\|=0$, $\|r_{\bot}^\alpha\|=1$ (if $r_{\bot}^\alpha\in\mathsf{R}^{BH(n)}$), and the equation for each $\beta\in A^{BH(0)}\cup R^{BH(n)}$ is
\begin{equation*}
	\|\beta\|=\prod_{i=1}^{n_{\beta}}(1-\|\alpha_i\|\|r_{\alpha_i}^{\beta}\|).
\end{equation*}
This equation semantics calculates the acceptability of $\beta$ as the product of the failure probabilities of all its attacks, where each attack fails with probability ($1-\|\alpha_i\|\|r_{\alpha_i}^{\beta}\|$), representing the combined strength of the attacker and its attack relation.

Denote that an assignment $\|\cdot\|$ is a model of a $BHAF$ under equational semantics $Eq_P^{BH}$ by $\|\cdot\| \models_{\mathcal{BH}_{[0,1]}^P} BHAF$. Then, we present the model equivalence theorem.
\begin{thm}\label{bhep}
	For a $BHAF=(\mathsf{A}^{BH(0)}, \mathsf{R}^{BH(n)})$ and an assignment $\|\cdot\|: \mathsf{A}^{BH(0)}\cup\mathsf{R}^{BH(n)}\to [0,1]$, 
	\begin{equation*}
		\|\cdot\| \models_{\mathcal{BH}_{[0,1]}^P} BHAF \Longleftrightarrow \|\cdot\| \models_{\mathcal{PL}_{[0,1]}^P} ec_{BH}(BHAF).
	\end{equation*}
\end{thm}
\begin{proof} 
	Similar to the proof of Theorem \ref{Thm3}.
\end{proof} 
\begin{cor}
The	$ec_{BH}$ is a translation of $\mathcal{BHAF}$ associated with $\mathcal{PL}_{[0,1]}^P$ under the equational semantics $Eq_P^{BH}$.
\end{cor}
\begin{cor}
	The equational semantics $Eq_P^{BH}$ is a fuzzy normal encoded semantics and thus it is translatable.
\end{cor}
\subsubsection{The equational system $Eq_L^{BH}$}
For a $BHAF=(\mathsf{A}^{BH(0)}, \mathsf{R}^{BH(n)})$, we give the equational system $Eq_L^{BH}$ such that $\|\bot\|=0$, $\|r_{\bot}^\alpha\|=1$ (if $r_{\bot}^\alpha\in\mathsf{R}^{BH(n)}$), and the equation for each $\beta\in A^{BH(0)}\cup R^{BH(n)}$ is 
\begin{equation*}
	\|\beta\|=\begin{cases}
		0 & \quad \sum_{i=1}^{n_{\beta}}y_i\leqslant n_{\beta}-1\\
		\sum_{i=1}^{n_{\beta}}y_i-n_{\beta}+1 & \quad \sum_{i=1}^{n_{\beta}}y_i> n_{\beta}-1
	\end{cases}
\end{equation*}
where 
\begin{equation*}
	y_i=\begin{cases}
		1 & \quad \|\alpha_i\|+\|r_{\alpha_i}^{\beta}\|\leqslant 1\\
		2-\|\alpha_i\|-\|r_{\alpha_i}^{\beta}\| & \quad \|\alpha_i\|+\|r_{\alpha_i}^{\beta}\|> 1
	\end{cases}.
\end{equation*}
This equation determines $\|\beta\|$ through a threshold mechanism: $\beta$ becomes acceptable only when the collective ineffectiveness $\sum_{i=1}^{n_{\beta}}y_i$ exceeds the threshold $n_{\beta}-1$, with its value then growing linearly with the excess, where each $y_i$ represents the ineffectiveness of attack $i$ and decreases linearly as $\|\alpha_i\|+\|r_{\alpha_i}^{\beta}\|$ increases beyond 1.

Denote that an assignment $\|\cdot\|$ is a model of a $BHAF$ under equational semantics $Eq_L^{BH}$ by $\|\cdot\| \models_{\mathcal{BH}_{[0,1]}^L} BHAF$. Then, we present the model equivalence theorem.
\begin{thm}\label{bhel}
	For a $BHAF=(\mathsf{A}^{BH(0)}, \mathsf{R}^{BH(n)})$ and an assignment $\|\cdot\|: \mathsf{A}^{BH(0)}\cup\mathsf{R}^{BH(n)}\to [0,1]$, 
	\begin{equation*}
		\|\cdot\| \models_{\mathcal{BH}_{[0,1]}^L} BHAF \Longleftrightarrow \|\cdot\| \models_{\mathcal{PL}_{[0,1]}^L} ec_{BH}(BHAF).
	\end{equation*}
\end{thm}

\begin{proof} 
		Similar to the proof of Theorem \ref{Thm4}.
\end{proof} 
\begin{cor}
The	$ec_{BH}$ is a translation of $\mathcal{BHAF}$ associated with $\mathcal{PL}_{[0,1]}^L$ under the equational semantics $Eq_L^{BH}$.
\end{cor}
\begin{cor}
	The equational semantics $Eq_L^{BH}$ is a fuzzy normal encoded semantics and thus it is translatable.
\end{cor}
The G\"{o}del semantics determines acceptability by the strongest attack, the Product semantics by the cumulative product of independent attack strengths, and the {\L}ukasiewicz semantics by whether the collective ineffectiveness of attacks surpasses a critical threshold.
\section{Encoding frameworks with sets of attacking arguments}\label{SETAF}
In this section, we first encode $SETAF$s under complete semantics to the $\mathcal{PL}_3^L$.  Then, we give three typical equational semantics of $SETAF$s and encode $SETAF$s to $\mathcal{PL}_{[0,1]}$s under these three equational semantics. It is shown that these semantics are translatable by proving the model equivalence. 

\subsection{Encoding $SETAF$s under complete semantics}
Gabbay encodes $SETAF$s to logic systems $CN$ and $CNN$ in \cite{gabbay2015attack}, where $CN$ is the classical logic with atomic strong negation and $CNN$ is two-world modal logic with world-swapping negation. We will encode $SETAF$s to the $\mathcal{PL}_3^L$ and get the theorem of model equivalence.

Similar to the approach for $HLAF$s, we first adapt the original definition and complete semantics of $SETAF$s to facilitate their encoding, and then introduce the corresponding normal encoding function.
\begin{defn}\label{setdefn}
	For an $SETAF=(A^S, R^S)$, the extended $SETAF$ is a pair $(\mathsf{A}^{S}, \mathsf{R}^{S})$, where $\mathsf{A}^{S}=A^{S}\cup \{\bot\}$, $\mathsf{R}^{S}=R^{S}\cup \{(\{\bot\}, b)\mid b\in A^{S} \text{ and } \nexists B\subseteq A^{S}: (B, b)\in R^{S}\}$, and for any assignment $\|\cdot\|$ of the extended $SETAF$ let $\|\bot\|=0$.
\end{defn}
The extended $SETAF$ is constructed by adding an imaginary singleton attacker $\{\bot\}$ for each unattacked argument in the original framework. Intuitively, adding an invalid attacker for an unattacked argument will not affect the validity of this argument. Henceforth, all references to an $SETAF$ or the $\mathcal{SETAF}$ shall denote an extended $SETAF$ and the set of all extended $SETAF$s, respectively. The meanings of the notations $A^{S}$ and $R^{S}$ remain unchanged. We denote $(B, a)\in \mathsf{R}^{S}$ as $r_B^a\in \mathsf{R}^{S}$. Then we update the complete semantics of $SETAF$s.

\begin{defn}
	For a given $SETAF = (\mathsf{A}^{S}, \mathsf{R}^{S})$, a \emph{complete labelling} of the $SETAF$ is a function $\|\cdot\|:  \mathsf{A}^{S}\rightarrow\{0, 1, \frac{1}{2}\}$ such that $\|\bot\|=0$ and for each $a\in A^{S}$
	\begin{equation*}
		\|a\| = \begin{cases}
			1 & \quad \text{ iff } \forall S_i (r_{S_i}^a\in \mathsf{R}^{S}) \exists b \in S_i: \|b\| = 0\\
			0 & \quad \text{ iff } \exists S_i (r_{S_i}^a\in \mathsf{R}^{S}) \forall b \in S_i: \|b\| = 1\\
			\frac{1}{2} & \quad \text{ otherwise }
		\end{cases}. 
	\end{equation*}
\end{defn}
Obviously, this definition preserves the original complete semantics. Then we present the normal encoding function of the $SETAF$.
\begin{defn}
	For a given $SETAF=(\mathsf{A}^{S}, \mathsf{R}^{S})$, the \emph{normal encoding function} of the $SETAF$ w.r.t. the $\mathcal{PL}_3^L$ is $ec_{S}$ such that
	\begin{equation*}
		ec_{S}(SETAF)=\bigwedge_{a\in A^S}(a\leftrightarrow\bigwedge_{r_{S_i}^a\in \mathsf{R}^{S}}(\neg \bigwedge_{b_{ij}\in S_i}b_{ij})).
	\end{equation*}
\end{defn}
This encoding captures the collective attack semantics of $SETAF$s by logically representing that an argument is accepted if and only if for every set attacking it, at least one member in that set is rejected. From the encoded formula, each $a \in A^S$ represents a non-imaginary argument, while $b_{ij}$ may be the imaginary argument $\bot$. 
\begin{rmk}
	The classical approach of defining the conjunction over an empty set as $\bot$ leads to a contradiction. Consider an unattacked argument $a$ and the formula $\bigwedge_{r_{S_i}^a\in \mathsf{R}^{S}} (\neg \bigwedge_{b_{ij} \in S_i} b_{ij})$. Since there are no attackers, the outer conjunction ranges over an empty set, yielding $\bot$. However, each inner conjunction $\bigwedge_{b_{ij} \in S_i} b_{ij}$ also ranges over an empty set within its negated context, which would likewise be $\bot$, making the entire expression $\bigwedge_{r_{S_i}^a\in \mathsf{R}^{S}} \neg \bot$. This is not equivalent to $\bot$, resulting in a contradiction. It is for this reason that we adopt the approach of adding an imaginary argument $\bot$ and imaginary attacks but not defining the conjunction over an empty set as $\bot$.
\end{rmk}  

Denote that an assignment $\|\cdot\|$ is a model of an $SETAF$ under complete semantics by $\|\cdot\| \models_{\mathcal{S}_3^C} SETAF$. Then, we present the model equivalence theorem.
\begin{thm}\label{setme}
	For an $SETAF=(\mathsf{A}^{S}, \mathsf{R}^{S})$ and an assignment $\|\cdot\|: \mathsf{A}^{S}\to \{0,1,\frac{1}{2}\}$, 
	\begin{equation*}
		\|\cdot\| \models_{\mathcal{S}_3^C} SETAF \Longleftrightarrow \|\cdot\| \models_{\mathcal{PL}_3^L} ec_{S}(SETAF).
	\end{equation*}
\end{thm}

\begin{proof} 
We need to check that for a given assignment $\|\cdot\|$ and $\forall a\in A^S$, any value of $a$ satisfies complete semantics iff $\|a\leftrightarrow\bigwedge_{r_{S_i}^a\in \mathsf{R}^{S}}(\neg \bigwedge_{b_{ij}\in S_i}b_{ij})\|=1$ in the $\mathcal{PL}_3^L$.  We need to discuss three cases.
	\begin{itemize}
		\item Case 1, $\|a\|=1$. 
		\\$\|a\|=1$ by model $\|\cdot\|$ under complete semantics
		\\$\Longleftrightarrow$ by complete semantics, $\forall r_{S_i}^a\in \mathsf{R}^{S}$, $\exists b_{ij} \in S_i, \|b_{ij}\|= 0$ 
		\\$\Longleftrightarrow$  $\forall r_{S_i}^a\in \mathsf{R}^{S}$, $\|\bigwedge_{b_{ij}\in S_i}b_{ij}\|= 0$ (i.e., $\|\neg\bigwedge_{b_{ij}\in S_i}b_{ij}\|= 1$) in the $\mathcal{PL}_3^L$
		\\$\Longleftrightarrow$ $\|\bigwedge_{r_{S_i}^a\in \mathsf{R}^{S}}(\neg \bigwedge_{b_{ij}\in S_i}b_{ij})\|= 1$ in the $\mathcal{PL}_3^L$
		\\$\Longleftrightarrow$ $\|a\leftrightarrow\bigwedge_{r_{S_i}^a\in \mathsf{R}^{S}}(\neg \bigwedge_{b_{ij}\in S_i}b_{ij})\|= 1$ in the $\mathcal{PL}_3^L$
		
		\item Case 2, $\|a\|=0$. 
		\\$\|a\|=0$ by model $\|\cdot\|$ under complete semantics
		\\$\Longleftrightarrow$ $\exists r_{S_i}^a\in \mathsf{R}^{S}$, $\forall$ $b_{ij} \in S_i$, $\|b_{ij}\|= 1$ by complete semantics
		\\$\Longleftrightarrow$ $\exists r_{S_i}^a\in \mathsf{R}^{S}$, $\|\bigwedge_{b_{ij}\in S_i}b_{ij}\|=1$ (i.e., $\|\neg\bigwedge_{b_{ij}\in S_i}b_{ij}\|=0$) in the $\mathcal{PL}_3^L$
		\\$\Longleftrightarrow$ $\|\bigwedge_{r_{S_i}^a\in \mathsf{R}^{S}}(\neg \bigwedge_{b_{ij}\in S_i}b_{ij})\|= 0$ in the $\mathcal{PL}_3^L$
		\\$\Longleftrightarrow$ $\|a\leftrightarrow\bigwedge_{r_{S_i}^a\in \mathsf{R}^{S}}(\neg \bigwedge_{b_{ij}\in S_i}b_{ij})\|= 1$ in the $\mathcal{PL}_3^L$
	
		\item Case 3, $\|a\|=\frac{1}{2}$.
		\\$\|a\|=\frac{1}{2}$ by model $\|\cdot\|$ under complete semantics
		\\$\Longleftrightarrow$ $\|a\|\neq1$ and $\|a\|\neq0$ under complete semantics
		\\$\Longleftrightarrow$ $\|\bigwedge_{r_{S_i}^a\in \mathsf{R}^{S}}(\neg \bigwedge_{b_{ij}\in S_i}b_{ij})\|\neq1$ and $\|\bigwedge_{r_{S_i}^a\in \mathsf{R}^{S}}(\neg \bigwedge_{b_{ij}\in S_i}b_{ij})\|\neq0$ in the $\mathcal{PL}_3^L$ by Case 1 and Case 2
		\\$\Longleftrightarrow$ $\|\bigwedge_{r_{S_i}^a\in \mathsf{R}^{S}}(\neg \bigwedge_{b_{ij}\in S_i}b_{ij})\|=\frac{1}{2}$ in the $\mathcal{PL}_3^L$
		\\$\Longleftrightarrow$ $\|a\leftrightarrow\bigwedge_{r_{S_i}^a\in \mathsf{R}^{S}}(\neg \bigwedge_{b_{ij}\in S_i}b_{ij})\|=1$ in the $\mathcal{PL}_3^L$.
	\end{itemize}

	From the three cases above, the value of  each $a\in A^S$ satisfies complete semantics iff $\|a\leftrightarrow\bigwedge_{r_{S_i}^a\in \mathsf{R}^{S}}(\neg \bigwedge_{b_{ij}\in S_i}b_{ij})\|=1$ in the $\mathcal{PL}_3^L$. 
	So, an assignment $\|\cdot\|$ of the $SETAF$ satisfies complete semantics iff $\|\bigwedge_{a\in A^S}(a\leftrightarrow\bigwedge_{r_{S_i}^a\in \mathsf{R}^{S}}(\neg \bigwedge_{b_{ij}\in S_i}b_{ij}))\|=\|ec_{S}(SETAF)\|=1$, i.e., $\|\cdot\|$ is a model of the $ec_{S}(SETAF)$, in the $\mathcal{PL}_3^L$.
\end{proof} 
\begin{cor}
The	$ec_{S}$ is a translation of $\mathcal{SETAF}$ associated with $\mathcal{PL}_3^L$ under complete semantics and thus complete semantics of $\mathcal{SETAF}$ is translatable.
\end{cor}

\subsection{Three equational semantics of $SETAF$s}
In this subsection, we emphasize three equational systems on [0,1] for $SETAF$s and then we prove the theorems of model equivalence. The encoding function of $SETAF$s w.r.t. a $\mathcal{PL}_{[0, 1]}$ is still $ec_{S}$ such that 
\begin{equation*}
	ec_{S}(SETAF)=\bigwedge_{a\in A^S}(a\leftrightarrow\bigwedge_{r_{S_i}^a\in \mathsf{R}^{S}}(\neg \bigwedge_{b_{ij}\in S_i}b_{ij})).
\end{equation*}

 In the following discussion, for a given $SETAF=(\mathsf{A}^{S}, \mathsf{R}^{S})$ and $a\in A^S$, let $n_a = |\{S_i \mid r_{S_i}^a\in \mathsf{R}^{S}\}|$ denote the number of attackers of $a$. For  $r_{S_i}^a\in \mathsf{R}^{S}$ and $b_{ij}\in S_i$, let $k_i= |S_i|=|\{b_{ij} \mid b_{ij} \in S_i\}|$ denote the cardinality of the set attacker $S_i$.

\subsubsection{The equational system $Eq_G^{S}$}
For an $SETAF=(\mathsf{A}^{S}, \mathsf{R}^{S})$, we give the equational system $Eq_G^{S}$ such that $\|\bot\|=0$ and the equation for each $a\in A^S$ is
\begin{equation*}
	\|a\|=\min_{i=1}^{n_a}\max_{j=1}^{k_i}(1-\|b_{ij}\|).
\end{equation*}
This equation $\|a\| = \min_{i=1}^{n_a} \max_{j=1}^{k_i} (1 - \|b_{ij}\|)$ captures the collective attack semantics of $SETAF$s through a two-level evaluation. For each set attacker $S_i$, the $\max$ operator determines the set's effectiveness by identifying the \emph{weakest} attacker in the set—specifically, the one with the highest ineffectiveness value $(1 - \|b_{ij}\|)$. This reflects that a set attack fails if \emph{any} of its members is ineffective. The outer $\min$ operator then selects the \emph{strongest} set attack—the one with the lowest ineffectiveness—establishing that argument $a$'s acceptability is constrained by its most potent collective threat. Thus, $a$ is fully acceptable only if all set attacks contain at least one completely ineffective member.

Denote that an assignment $\|\cdot\|$ is a model of an $SETAF$ under equational semantics $Eq_G^{S}$ by $\|\cdot\| \models_{\mathcal{S}_{[0,1]}^G} SETAF$. Then, we present the model equivalence theorem.
\begin{thm}\label{sgme}
	For an $SETAF=(\mathsf{A}^{S}, \mathsf{R}^{S})$ and an assignment $\|\cdot\|: \mathsf{A}^{S}\to [0,1]$, 
	\begin{equation*}
		\|\cdot\| \models_{\mathcal{S}_{[0,1]}^G} SETAF \Longleftrightarrow \|\cdot\| \models_{\mathcal{PL}_{[0,1]}^G} ec_{S}(SETAF).
	\end{equation*}
\end{thm}

\begin{proof} 
		A model $\|\cdot\|$ of $ec_{S}(SETAF)$ in the $\mathcal{PL}_{[0, 1]}^G$
	\\$\Longleftrightarrow $ (a solution of) $\|\bigwedge_{a\in A^S}(a\leftrightarrow\bigwedge_{r_{S_i}^a\in \mathsf{R}^{S}}(\neg \bigwedge_{b_{ij}\in S_i}b_{ij}))\|=1$
	\\$\Longleftrightarrow $ for each $a\in A^S$, $\|a\leftrightarrow\bigwedge_{r_{S_i}^a\in \mathsf{R}^{S}}(\neg \bigwedge_{b_{ij}\in S_i}b_{ij})\|=1$
	\\$\Longleftrightarrow $ for each $a\in A^S$, $\|a\|=\|\bigwedge_{r_{S_i}^a\in \mathsf{R}^{S}}(\neg \bigwedge_{b_{ij}\in S_i}b_{ij})\|$ by Lemma \ref{remark2}
	\\$\Longleftrightarrow $ for each $a\in A^S$, $\|a\|=\min_{i=1}^{n_a}\|\neg \bigwedge_{b_{ij}\in S_i}b_{ij}\|=\min_{i=1}^{n_a}(1-\|\bigwedge_{b_{ij}\in S_i}b_{ij}\|)$ $=\min_{i=1}^{n_a}(1-\min_{j=1}^{k_i}\|b_{ij}\|)=\min_{i=1}^{n_a}\max_{j=1}^{k_i}(1-\|b_{ij}\|)$
	\\$\Longleftrightarrow$ a model $\|\cdot\|$ of the $SETAF$ under equational semantics $Eq_G^{S}$.
\end{proof} 
\begin{cor}
The	$ec_{S}$ is a translation of $\mathcal{SETAF}$ associated with $\mathcal{PL}_{[0,1]}^G$ under the equational semantics $Eq_G^{S}$.
\end{cor}
\begin{cor}
	The equational semantics $Eq_G^{S}$ is a fuzzy normal encoded semantics and thus it is translatable.
\end{cor}

\subsubsection{The equational system $Eq_P^{S}$}
For an $SETAF=(\mathsf{A}^{S}, \mathsf{R}^{S})$, we give the equational system $Eq_P^{S}$ such that $\|\bot\|=0$ and the equation for each $a\in A^S$ is
\begin{equation*}
	\|a\|=\prod_{i=1}^{n_a}(1-\prod_{j=1}^{k_i}\|b_{ij}\|).
\end{equation*}
The equational semantics $Eq_P^{S}$ for $SETAF$s evaluates argument acceptability through set attackers. Each term $1 - \prod_{j=1}^{k_i} \|b_{ij}\|$ represents the failure probability of set attacker $S_i$.  
This probability equals 1 if and only if at least one argument $b_{ij} \in S_i$ has $\|b_{ij}\| = 0$.  
Thus, $\|a\| = 1$ only when all attacks completely fail, meaning every set attacker contains at least one fully rejected argument.

Denote that an assignment $\|\cdot\|$ is a model of an $SETAF$ under equational semantics $Eq_P^{S}$ by $\|\cdot\| \models_{\mathcal{S}_{[0,1]}^P} SETAF$. Then, we present the model equivalence theorem.
\begin{thm}\label{spme}
	For an $SETAF=(\mathsf{A}^{S}, \mathsf{R}^{S})$ and an assignment $\|\cdot\|: \mathsf{A}^{S}\to [0,1]$, 
	\begin{equation*}
		\|\cdot\| \models_{\mathcal{S}_{[0,1]}^P} SETAF \Longleftrightarrow \|\cdot\| \models_{\mathcal{PL}_{[0,1]}^P} ec_{S}(SETAF).
	\end{equation*}
\end{thm}

\begin{proof} 
	A model $\|\cdot\|$ of $ec_{S}(SETAF)$ in the $\mathcal{PL}_{[0, 1]}^P$
	\\$\Longleftrightarrow $ (a solution of) $\|\bigwedge_{a\in A^S}(a\leftrightarrow\bigwedge_{r_{S_i}^a\in \mathsf{R}^{S}}(\neg \bigwedge_{b_{ij}\in S_i}b_{ij}))\|=1$
	\\$\Longleftrightarrow $ for each $a\in A^S$, $\|a\leftrightarrow\bigwedge_{r_{S_i}^a\in \mathsf{R}^{S}}(\neg \bigwedge_{b_{ij}\in S_i}b_{ij})\|=1$
	\\$\Longleftrightarrow $ for each $a\in A^S$, $\|a\|=\|\bigwedge_{r_{S_i}^a\in \mathsf{R}^{S}}(\neg \bigwedge_{b_{ij}\in S_i}b_{ij})\|$ by Lemma \ref{remark2}
	\\$\Longleftrightarrow $ for each $a\in A^S$, $\|a\|=\prod_{i=1}^{n_a}\|\neg \bigwedge_{b_{ij}\in S_i}b_{ij}\|=\prod_{i=1}^{n_a}(1-\|\bigwedge_{b_{ij}\in S_i}b_{ij}\|)=\prod_{i=1}^{n_a}(1-\prod_{j=1}^{k_i}\|b_{ij}\|)$
	\\$\Longleftrightarrow$ a model $\|\cdot\|$ of the $SETAF$ under equational semantics $Eq_P^{S}$.
\end{proof} 
\begin{cor}
The	$ec_{S}$ is a translation of $\mathcal{SETAF}$ associated with $\mathcal{PL}_{[0,1]}^P$ under the equational semantics $Eq_P^{S}$.
\end{cor}
\begin{cor}
The equational semantics $Eq_P^{S}$ is a fuzzy normal encoded semantics and thus it is translatable.
\end{cor}
\subsubsection{The equational system $Eq_L^{S}$}
For an $SETAF=(\mathsf{A}^{S}, \mathsf{R}^{S})$, we give the equational system $Eq_L^{S}$ such that $\|\bot\|=0$ and the equation for each $a\in A^S$ is
\begin{equation*}
	\|a\|=\begin{cases}
		0 & \quad \sum_{i=1}^{n_a}y_i\leqslant n_a-1\\
		\sum_{i=1}^{n_a}y_i-n_a+1 & \quad \sum_{i=1}^{n_a}y_i> n_a-1
	\end{cases}
\end{equation*}
where
\begin{equation*}
	y_i=\begin{cases}
		1 & \quad \sum_{j=1}^{k_i}\|b_{ij}\|\leqslant k_i-1\\
		k_i-\sum_{j=1}^{k_i}\|b_{ij}\| & \quad \sum_{j=1}^{k_i}\|b_{ij}\|> k_i-1
	\end{cases}.
\end{equation*}
The equational semantics $Eq_L^{S}$ for $SETAF$s evaluates argument acceptability using a threshold mechanism. For each argument $a$, we compute the collective ineffectiveness $y_i$ of its set attackers $S_i$, where $y_i = 1$ if the sum of truth values in $S_i$ is at most $k_i-1$ and $y_i$ decreases linearly otherwise. The final acceptability $\|a\|$ is determined by comparing the total ineffectiveness $\sum y_i$ against threshold $n_a-1$: $\|a\| = 0$ if the threshold is not exceeded, and increases linearly with the excess otherwise.

Denote that an assignment $\|\cdot\|$ is a model of an $SETAF$ under equational semantics $Eq_L^{S}$ by $\|\cdot\| \models_{\mathcal{S}_{[0,1]}^L} SETAF$. Then, we present the model equivalence theorem.
\begin{thm}\label{slme}
	For an $SETAF=(\mathsf{A}^{S}, \mathsf{R}^{S})$ and an assignment $\|\cdot\|: \mathsf{A}^{S}\to [0,1]$, 
	\begin{equation*}
		\|\cdot\| \models_{\mathcal{S}_{[0,1]}^L} SETAF \Longleftrightarrow \|\cdot\| \models_{\mathcal{PL}_{[0,1]}^L} ec_{S}(SETAF).
	\end{equation*}
\end{thm}
\begin{proof} 
A model $\|\cdot\|$ of $ec_{S}(SETAF)$ in the $\mathcal{PL}_{[0, 1]}^L$
\\$\Longleftrightarrow $ (a solution of)
	$\|\bigwedge_{a\in A^S}(a\leftrightarrow\bigwedge_{r_{S_i}^a\in \mathsf{R}^{S}}(\neg \bigwedge_{b_{ij}\in S_i}b_{ij}))\|=1$
	\\$\Longleftrightarrow $ for each $a\in A^S$, $\|a\leftrightarrow\bigwedge_{r_{S_i}^a\in \mathsf{R}^{S}}(\neg \bigwedge_{b_{ij}\in S_i}b_{ij})\|=1$
	\\$\Longleftrightarrow $ for each $a\in A^S$, $\|a\|=\|\bigwedge_{r_{S_i}^a\in \mathsf{R}^{S}}(\neg \bigwedge_{b_{ij}\in S_i}b_{ij})\|$ by Lemma \ref{remark2}
	\\$\Longleftrightarrow $ for each $a\in A^S$, by the formula \ref{n-t-norm},	
	\begin{equation}\|a\|=
		\begin{cases}
			0 & \quad \sum_{i=1}^{n_a}y_i \leq n_a-1\\
			\sum_{i=1}^{n_a}y_i-n_a+1 & \quad \sum_{i=1}^{n_a}y_i > n_a-1
		\end{cases}
	\end{equation}
	where $y_i=\|\neg \bigwedge_{b_{ij}\in S_i}b_{ij}\|=1-\|\bigwedge_{b_{ij}\in S_i}b_{ij}\|$.
	By the formula \ref{n-t-norm} again, we have
	\begin{equation*}
		\|\bigwedge_{b_{ij}\in S_i}b_{ij}\|=
		\begin{cases}
			0 & \quad \sum_{j=1}^{k_i}\|b_{ij}\| \leq k_i-1\\
			\sum_{j=1}^{k_i}\|b_{ij}\|-k_i+1 & \quad \sum_{j=1}^{k_i}\|b_{ij}\| > k_i-1
		\end{cases},
	\end{equation*} 
	i.e., 
	\begin{equation*}
		y_i=
		\begin{cases}
			1 & \quad \sum_{j=1}^{k_i}\|b_{ij}\| \leq k_i-1\\
			k_i-\sum_{j=1}^{k_i}\|b_{ij}\| & \quad \sum_{j=1}^{k_i}\|b_{ij}\| > k_i-1
		\end{cases}.
	\end{equation*}
	The above equations are the same as the equational system $Eq_L^S$. So $\|\cdot\|$ is a model of $ec_{S}(SETAF)$ in the $\mathcal{PL}_{[0, 1]}^L$ iff it is a model of the $SETAF$ under equational semantics $Eq_L^{S}$, i.e., we have proven this theorem.	
\end{proof} 
\begin{cor}
The	$ec_{S}$ is a translation of $\mathcal{SETAF}$ associated with $\mathcal{PL}_{[0,1]}^L$ under the equational semantics $Eq_L^{S}$.
\end{cor}
\begin{cor}
The equational semantics $Eq_L^{S}$ is a fuzzy normal encoded semantics and thus it is translatable.
\end{cor}

\section{Encoding higher-order set argumentation frameworks}\label{sec6}
In $RAFN$s \cite{cayrol2018structure}, the source of an attack is a single argument, whereas the source of a support is a nonempty set of arguments. The target of either can be an argument, an attack, or a support.
In $REBAF$s \cite{cayrol2018argumentation}, the source of an attack or a support is a nonempty set of arguments, while the target may be an argument, an attack, or a support. Additionally, the concept of prima-facie elements is introduced in $REBAF$s.
Both $RAFN$s and $REBAF$s do not allow a source to include attacks or supports. However, since an attack can itself be construed as an argument and is subject to being attacked, it must also be able to participate in a set-attack. Next, we present two concrete argumentation scenarios to illustrate the necessity of allowing existing attacks to be elements of a set-attack's source.

The first scenario involves Party X (claimant) and Party Y (defendant) in a dispute over the legal validity of their signed sales contract (Party X claims validity, while Y intends to refute it via a set-attack). The relevant arguments are as follows: 
\begin{itemize}
	\item Argument $P$ (core claim by X): The sales contract is legally valid (supporting grounds: meets formal requirements like mutual signatures and no statutory invalidity causes).
	\item Argument $A$ (by Y): The contract lacks a key formal requirement—X's signature was forged (directly challenging $P$'s ``formal validity").
	\item Argument $B$ (by X): Even with a forged signature, Y orally agreed to terms and accepted goods, constituting ``de facto performance" (a statutory exception to formal requirements, defending $P$ against $A$).
	\item Argument $C$ (by Y): Y's acceptance resulted from misunderstanding (mistaking goods for a third party's), not recognition of the forged contract, so $B$ cannot justify validity.
\end{itemize}
The basic attack $\alpha$ from $C$ to $B$ occurs because $C$ negates $B$'s core premise, thus invalidating $B$ as a defense for $P$. The set-attack by $\{A, \alpha\}$ on $P$ is logically necessary since neither component alone achieves a full refutation:
\begin{itemize}
	\item $A$ alone: $X$ can use $B$ to defend $P$, making $A$ insufficient.
	\item $\alpha$ alone: only invalidates $B$ without directly challenging $P$.
	\item $A$ and $\alpha$ combined: $A$ undermines $P$'s validity while $\alpha$ removes its defense ($B$), thus fully refuting $P$.
\end{itemize}

Then we present the second argumentation scenario. This is a classic argumentation involving penguins. Consider two agents, M and N: M first puts forward an opinion $\beta$ that all kinds of birds can fly. Subsequently, N thinks of penguins and presents an argument $\gamma$ that penguins are birds. At this stage, $\gamma$ does not attack $\beta$, as it remains unknown whether penguins can fly. The dialogue proceeds: to reinforce his argument $\beta$, M puts forward another opinion $\delta$ that penguins have wings, so they can fly. At this juncture, no attacks have been established among the three arguments ($\beta$, $\gamma$, and $\delta$). Next, N concludes the argumentation by presenting a critical argument $\epsilon$: an animal with wings cannot fly if its wings are too small compared to its body, failing to provide sufficient propulsion for flight. We now analyze the attacks in this scenario: $\epsilon$ does not directly attack $\beta$ or $\gamma$, but directly attacks $\delta$---we denote this attack as $\eta$. $\eta$ amounts to claiming that penguins cannot fly. In the absence of $\gamma$, $\eta$ does not attack $\beta$, as penguins might not be classified as birds. Instead, it is the collective set of $\gamma$ and $\eta$ that attacks $\beta$, and we denote this set attack as $\varphi$. Thus, this scenario demonstrates the need to introduce attacks into set attackers (i.e., higher-order collective attacks), which motivates the proposal of $HSAF$s.

This section presents the core contributions of this work by introducing and formally developing the $HSAF$—a unified structure that generalizes both higher-order attacks and collective attacks. In an $HSAF$, targets can be either arguments or attacks, while attackers may be sets comprising both arguments and attacks, thus enabling richer and more expressive interaction patterns. The main objectives of this chapter are:
\begin{itemize}
	\item To formally define $HSAF$s and extend the complete semantics from $BHAF$s and $SETAF$s to this generalized $AF$.
	\item To establish logical encodings of $HSAF$s under complete semantics into $\mathcal{PL}_3^L$, and to prove model equivalence.
	\item To propose three types of equational semantics for $HSAF$s—based on $\mathcal{PL}_{[0,1]}^G$, $\mathcal{PL}_{[0,1]}^P$, and $\mathcal{PL}_{[0,1]}^L$—and encode them to corresponding $\mathcal{PL}_{[0,1]}$s, with accompanying translatability results.
	\item To introduce a formal equational approach for $HSAF$s, relating continuous fuzzy normal encoded equational systems to general real equational systems, and establishing solution existence via fixed-point theorems.
	\item To analyze semantic correspondences between fuzzy normal encoded semantics and complete semantics of $HSAF$s.
	\item To unify the family of $AFSA$s by providing systematic transformations from each $AFSA$ type into semantically equivalent $SETAF$s and $HSAF$s.
\end{itemize}
Through these contributions, this chapter not only expands the syntactic and semantic scope of argumentation frameworks but also offers a unified logical and equational foundation for reasoning about complex attack structures.
 
For $HSAF$s, related concepts of the encoding method are very similar to those of $HOAF$s and $SETAF$s in Subsection \ref{sub2.3}, so we will omit the detailed introduction.
\subsection{The $HSAF$ and its complete semantics }\label{sub6.1}

First, we give the definition of $HSAF$s. In this definition, we follow the approach used in defining $BHAF$s by introducing an imaginary argument $\bot$ and its associated attacks. Additionally, the level of an $HSAF$ is emphasized.
\begin{defn}\label{hsafd}
A \emph{0-level preparatory higher-order set argumentation framework} ($pre\text{-}HSAF^{(0)}$) is a pair $(A^{HS(0)}, R^{HS(0)})$, where $R^{HS(0)}\subseteq (2^{A^{HS(0)}} \setminus\{\emptyset\}) \times A^{HS(0)}$. An \emph{$n$-level preparatory higher-order set argumentation framework} ($pre\text{-}HSAF^{(n)}$) is a pair 
 $(A^{HS(0)}, R^{HS(n)})$, where $n\geqslant1$, $A^{HS(n)}=A^{HS(0)}\cup ((2^{A^{HS(n-1)}} \setminus\{\emptyset\}) \times A^{HS(n-1)})$, $R^{HS(n)}\subseteq(2^{A^{HS(n)}} \setminus\{\emptyset\}) \times A^{HS(n)}$, and $R^{HS(n)}\nsubseteq(2^{A^{HS(n-1)}} \setminus\{\emptyset\}) \times A^{HS(n-1)}$. 
 
An \emph{$n$-level higher-order set argumentation framework} ($HSAF^{(n)}$)  is a pair $(\mathsf{A}^{HS(0)}, \mathsf{R}^{HS(n)})$, where $n\geqslant0$, $\mathsf{A}^{HS(0)}=A^{HS(0)}\cup\{\bot\}$, $\mathsf{R}^{HS(n)}=R^{HS(n)}\cup \{(\{\bot\}, \beta)\mid \beta\in A^{HS(0)}\cup R^{HS(n)} \text{ and } \nexists \alpha\in A^{HS(0)}\cup R^{HS(n)}: (\alpha, \beta)\in R^{HS(n)}\}$. 

Specifically, in an $HSAF^{(n)}$, $\mathsf{A}^{HS(0)}$ is the finite set of arguments including the imaginary argument $\bot$, and $\mathsf{R}^{HS(n)}$ is the $n$-level attack relation that may include some imaginary attacks $(\{\bot\}, \beta)$, where $\beta\in {A}^{HS(0)}\cup {R}^{HS(n)}$.
For any assignment $\|\cdot\|$ of an $HSAF^{(n)}$, let $\|\bot\|=0$ and $\|(\{\bot\}, \beta)\|=1$ (if $(\{\bot\}, \beta)\in \mathsf{R}^{HS(n)}$).
\end{defn}

Following a similar understanding of the definition of $BHAF^{(n)}$s, we can understand that of $HSAF^{(n)}$s, where each attacker is a set in $HSAF^{(n)}$s.
From this definition, the difference between an $HSAF^{(n)}$ and a $pre\text{-}HSAF^{(n)}$ lies in the fact that we assign to each argument or attack that is not attacked in the $pre\text{-}HSAF^{(n)}$ an invalid attacker $\{\bot\}$. We hold that adding a thoroughly invalid attacker to an argument or attack does not affect it in any way. Thus, the distinction between a $pre\text{-}HSAF^{(n)}$ and a $HSAF^{(n)}$ resides in form rather than essence. Briefly, we use the ``$HSAF$" to represent a higher-order set argumentation framework at any finite level.
The particular case of the $HSAF$ is that only subsets of $A^{HS(0)}$ can attack others. Note that an $HSAF$ where each attacker is a set of arguments is syntactically equivalent to a $REBAF$ that has only the attack relation.
Indeed, an $SETAF$ is an $HSAF^{(0)}$ when we ignore the imaginary attacks of set attacks in the $HSAF^{(0)}$. Under the convention of not distinguishing an element from its singleton, an $HLAF$ can be characterized as an $HSAF^{(n)}$ where every set attacker is a singleton of an argument, while a $BHAF$ corresponds to an $HSAF^{(n)}$ in which every set attacker is a singleton of either an argument or an attack. The hierarchy for $AFSA$s is shown in Figure \ref{fig:af_hierarchy}.
\begin{rmk}
	Another approach is to define a preparatory higher-order set argumentation framework ($pre\text{-}HSAF$) as $pre\text{-}HSAF = (A^{HS}, R^{HS})$, where $A^{HS}$ is a set of arguments and $R^{HS} \subseteq (2^{A^{HS} \cup R^{HS}} \setminus\{\emptyset\}) \times A^{HS}$. However, such a definition fails to characterize the level of a $pre\text{-}HSAF$.
\end{rmk}

We give some notations here. If $(S_i, \beta)\in\mathsf{R}^{HS(n)}$, then $S_i$ is called a set attacker of $\beta$, i.e., the set $S_i$ attacks $\beta$. Denoting $(S_i, \beta)$ as $r_{S_i}^{\beta}$, $S_i$ is also called the source of $r_{S_i}^{\beta}$ and $\beta$ is called the target of $r_{S_i}^{\beta}$. We denote that $S_i=\mathbf{s}(r_{S_i}^{\beta})$ and $\beta=\mathbf{t}(r_{S_i}^{\beta})$. And then, for $\mathsf{R}\subseteq\mathsf{R}^{HS(n)}$,  we denote that $\mathbf{S}(\mathsf{R})=\{\mathbf{s}(r_{S_i}^{\beta})\mid r_{S_i}^{\beta}\in\mathsf{R}\}$ and $\mathbf{T}(\mathsf{R})=\{\mathbf{t}(r_{S_i}^{\beta})\mid r_{S_i}^{\beta}\in\mathsf{R}\}$. Particularly, we have $\mathbf{T}(\mathsf{R}^{HS(n)})=A^{HS(0)}\cup R^{HS(n)}$.

Next, we present the complete semantics of $HSAF$s, which is generalized by combining those of $SETAF$s and $BHAF$s.
\begin{defn}\label{defn24}
	A \emph{complete labelling} for a $HSAF=(\mathsf{A}^{HS(0)}, \mathsf{R}^{HS(n)})$ is a function
	$\|\cdot\| : \mathsf{A}^{HS(0)}\cup \mathsf{R}^{HS(n)}\rightarrow\{0, 1, \frac{1}{2}\}$ such that $\|\bot\|=0$, $\|r_{\{\bot\}}^{\gamma}\|=1$ (if $r_{\{\bot\}}^{\gamma}\in \mathsf{R}^{HS(n)}$) and for each $\beta\in \mathbf{T}(\mathsf{R}^{HS(n)})$
	\begin{equation*}
		\|\beta\|= 
		\begin{cases}
			1 & \quad \text{iff } \forall S_i(r_{S_i}^{\beta}\in\mathsf{R}^{HS(n)}) \exists b_{ij} \in S_i: \|b_{ij}\| = 0 \text{ or } \|r_{S_i}^{\beta}\| = 0
			\\ 0 & \quad \text{iff } \exists S_i(r_{S_i}^{\beta}\in\mathsf{R}^{HS(n)}) \forall b_{ij} \in S_i: \|b_{ij}\|= 1 \text{ and } \|r_{S_i}^{\beta}\| =1
			\\\frac{1}{2} & \quad otherwise
		\end{cases}.
	\end{equation*}
\end{defn}
The complete semantics of $HSAF$s inherits the semantic features of both $BHAF$s and $SETAF$s. For the validity of a set attacker: the set attacker is valid if and only if each member in the set attacker is valid; it is invalid if and only if there is an invalid member in the set attacker; and all other cases correspond to a semi-valid set attacker. For the validity of a set attacking action: a set attacking action is considered valid if and only if both the set attacker and its corresponding attack are valid; it is invalid if either the set attacker or the attack is invalid; all other cases result in a semi-valid set attacking action. 

We then define the stable, preferred, and grounded semantics for an $HSAF=(\mathsf{A}^{HS(0)}, \mathsf{R}^{HS(n)})$, in analogy with the $DAF$, by requiring that stable labellings be two-valued complete labellings, preferred labellings be maximal complete labellings, and the grounded labelling be the smallest complete labelling.

\begin{defn}
	A \emph{stable labelling} $\|\cdot\|^s$ for the $HSAF$ is a complete labelling $\|\cdot\|^c$ such that for any $\beta\in \mathbf{T}(\mathsf{R}^{HS(n)})$ if $\|\beta\|^s\neq1$ then $\|\beta\|^s=0$.
	
	A \emph{preferred labelling} $\|\cdot\|^p$ for the $HSAF$ is a complete labelling such that there is no complete labelling $\|\cdot\|^c$ satisfying:
	\begin{itemize}
		\item for any $\beta\in \mathbf{T}(\mathsf{R}^{HS(n)})$, if $\|\beta\|^p=1$ then $\|\beta\|^c=1$;
		\item$\exists \gamma\in \mathbf{T}(\mathsf{R}^{HS(n)})$, s.t. $\|\gamma\|^c=1$ and $\|\gamma\|^p\neq1$.
	\end{itemize}
	
	A \emph{grounded labelling} $\|\cdot\|^g$ for the $HSAF$ is a complete labelling such that for any complete labelling $\|\cdot\|^c$ and any $\beta\in \mathbf{T}(\mathsf{R}^{HS(n)})$ if $\|\beta\|^g=1$ then $\|\beta\|^c=1$.
\end{defn}

In the following, we use an example to illustrate definitions of the $HSAF$ and its semantics.
\begin{exmp}\label{ex:hs}
	Suppose that we have an argument set $A^{HS(0)}=\{a, b, c, d\}$ and attacks as follows: $R^{HS(0)}=\{(\{a,b\}, c), (\{d\}, c), (\{d\}, d)\}$ (then denoting $(\{a,b\}, c)$ by $\alpha$, denoting $(\{d\}, c)$ by $\epsilon$ and denoting $(\{d\}, d)$ by $\varepsilon$), $R^{HS(1)}=\{\alpha, \epsilon, \varepsilon, (\{c\}, \alpha)\}$ (then denoting $(\{c\}, \alpha)$ by $\beta$), and $R^{HS(2)}=\{\alpha, \epsilon,\varepsilon, \beta, (\{\alpha\}, \beta), (\{\beta\}, a)\}$ (denoting $(\{\alpha\}, \beta)$ by $\gamma$ and denoting $(\{\beta\}, a)$ by $\delta$). Then we have $pre\text{-}HSAF^{(2)}=({A}^{HS(0)}, {R}^{HS(2)})$. Since $b, \gamma, \delta, \epsilon$ and $\varepsilon$ are not attacked in the $pre\text{-}HSAF^{(2)}$, we have $HSAF^{(2)}=(\mathsf{A}^{HS(0)}, \mathsf{R}^{HS(2)})$, where $\mathsf{A}^{HS(0)}=\{a, b, c, d, \bot\}$ and $\mathsf{R}^{HS(2)}=\{\alpha, \epsilon,\varepsilon, \beta,\gamma, \delta, (\{\bot\},b), (\{\bot\},\epsilon),$ $(\{\bot\},\varepsilon), (\{\bot\},\gamma), (\{\bot\},\delta)\}$.	
	
	Next, let us find models of the $HSAF^{(2)}$ under complete semantics. By Definition \ref{defn24}, we have $\|\bot\| = 0$ and $\|(\{\bot\}, x)\| = 1$ for all $x \in \{b, \epsilon, \varepsilon, \gamma, \delta\}$. For brevity, in this example and its continued examples, these trivial cases will be omitted when referring to models of the $HSAF^{(2)}$. Since $b$, $\varepsilon$, $\gamma$, $\delta$, and $\epsilon$ all have only the set attacker ${\{\bot\}}$, it follows from Definition \ref{defn24} that $\|b\|=\|\varepsilon\|=\|\gamma\|=\|\delta\|=\|\epsilon\|=1$. Then we obtain $\|d\|=\frac{1}{2}$ immediately. Then, we discuss three cases according to values of $a$.
	\begin{itemize}
		\item Case 1, $\|a\|=1$. If $\|\alpha\|=1$, then $\|c\|=0$ and $\|\beta\|=0$. If $\|\alpha\|=\frac{1}{2}$, then $\|c\|=\frac{1}{2}$ and $\|\beta\|=\frac{1}{2}$ which contradicts $\|a\|=\|\delta\|=1$. If $\|\alpha\|=0$, then $\|c\|=\frac{1}{2}$ and $\|\beta\|=1$ which contradicts $\|\alpha\|=0$.
		\item Case 2, $\|a\|=\frac{1}{2}$. If $\|\alpha\|=1$, then $\|c\|=\frac{1}{2}$ and $\|\beta\|=0$ which contradicts $\|a\|=\frac{1}{2}$ and $\|\delta\|=1$. If $\|\alpha\|=\frac{1}{2}$, then $\|c\|=\frac{1}{2}$ and $\|\beta\|=\frac{1}{2}$. If $\|\alpha\|=0$, then $\|c\|=\frac{1}{2}$ and $\|\beta\|=1$ which contradicts $\|\alpha\|=0$.
		\item Case 3, $\|a\|=0$. Then we have $\|c\|=\frac{1}{2}$. If $\|\alpha\|=1$, then $\|\beta\|=0$ which contradicts $\|a\|=0$ and $\|\delta\|=1$. If $\|\alpha\|=\frac{1}{2}$, then $\|\beta\|=\frac{1}{2}$ which contradicts $\|a\|=0$ and $\|\delta\|=1$. If $\|\alpha\|=0$, then $\|\beta\|=1$ which contradicts $\|\alpha\|=0$ and $\|\delta\|=1$.
	\end{itemize}
	Therefore, we obtain two complete labellings of the $HSAF^{(2)}$ as follows.
	\begin{itemize}
		\item $COM_1$: $\{\|a\|=\|b\|=\|\varepsilon\|=\|\gamma\|=\|\delta\|=\|\epsilon\|=\|\alpha\|=1, \|d\|=\frac{1}{2}, \|c\|=\|\beta\|=0\}$
		\item $COM_2$: $\{\|b\|=\|\varepsilon\|=\|\gamma\|=\|\delta\|=\|\epsilon\|=1, \|a\|=\|c\|=\|d\|=\|\alpha\|=\|\beta\|=\frac{1}{2}\}$		
	\end{itemize}
	There is no stable labelling.  
	The preferred labelling is $COM_1$ and the grounded labelling is $COM_2$.
\end{exmp}
\subsection{Encoding $HSAF$s  under complete semantics}

In this subsection, we encode $HSAF$s to the $\mathcal{PL}_3^L$ and explore the model relationship. Note that for any assignment $\|\cdot\|$ in any $\mathcal{PLS}$, we always let $\|\bot\|=0$ and $\|(\{\bot\},\gamma)\|=1$ (if $(\{\bot\},\gamma)\in\mathsf{R}^{HS(n)}$ for an $HSAF=(\mathsf{A}^{HS(0)}, \mathsf{R}^{HS(n)})$). We present the normal encoding function below.

\begin{defn}
For a given $HSAF=(\mathsf{A}^{HS(0)}, \mathsf{R}^{HS(n)})$, the \emph{normal encoding function} of the $HSAF$ w.r.t. the $\mathcal{PL}_3^L$ is $ec_{HS}$ such that
\begin{equation*}
	ec_{HS}(HSAF)=\bigwedge_{\beta\in \mathbf{T}(\mathsf{R}^{HS(n)})}(\beta\leftrightarrow\bigwedge_{r_{S_i}^\beta\in \mathsf{R}^{HS(n)}}\neg(r_{S_i}^{\beta}\wedge\bigwedge_{b_{ij}\in S_i}b_{ij})).
\end{equation*}
\end{defn}
This encoding function $ec_{HS}$ translates the entire $HSAF$ structure into a single logical formula in $\mathcal{PL}_3^L$. The encoding captures the complete semantics by creating a bi-implication for each element $\beta$ in the framework: $\beta$ is acceptable if and only if for every set attack, either the attack $r_{S_i}^\beta$ is invalid or at least one element in the attacking set $S_i$ is rejected. The outermost conjunction ensures all elements and their complex attack relationships are simultaneously represented in the resulting propositional formula.
To illustrate the normal encoding, a continued example is given below.
\addtocounter{exmp}{-1}
\begin{exmp}[continued]\label{ex:hs-cont1}
	For this given $HSAF^{(2)}$, we have $\mathbf{T}(\mathsf{R}^{HS(n)})= \{a,b,c,d,\alpha,\epsilon,\varepsilon,\beta,\gamma,\delta\}$. Applying the normal encoding function, we have
	\begin{align*}
		ec_{HS}(HSAF)=&
		(a\leftrightarrow \neg(\delta\wedge \beta))
		\wedge(b\leftrightarrow \neg(\bot\wedge r_\bot^b))\\&
		\wedge(c\leftrightarrow \neg(d\wedge \epsilon)\wedge\neg(a\wedge b\wedge \alpha))\\&
		\wedge(d\leftrightarrow \neg(d\wedge\varepsilon)
		\wedge(\delta\leftrightarrow \neg(\bot\wedge r_\bot^ \delta))\\&
		\wedge(\alpha\leftrightarrow \neg(c\wedge \beta))
		\wedge(\epsilon\leftrightarrow \neg(\bot\wedge r_\bot^ \epsilon))\\&
		\wedge(\varepsilon\leftrightarrow \neg(\bot\wedge r_\bot^ \varepsilon))
		\wedge(\beta\leftrightarrow \neg(\alpha\wedge \gamma))\\&
		\wedge(\gamma\leftrightarrow \neg(\bot\wedge r_\bot^ \gamma)).
	\end{align*}
\end{exmp}

Denote that an assignment $\|\cdot\|$ is a model of an $HSAF$ under complete semantics by $\|\cdot\| \models_{\mathcal{HS}_3^C} HSAF$. Then, we present the model equivalence theorem.
\begin{thm}\label{thm-pl3}
	For an $HSAF=(\mathsf{A}^{HS(0)}, \mathsf{R}^{HS(n)})$ and an assignment $\|\cdot\|: \mathsf{A}^{HS(0)}\cup\mathsf{R}^{HS(n)}\to \{0,1,\frac{1}{2}\}$, 
	\begin{equation*}
		\|\cdot\| \models_{\mathcal{HS}_3^C} HSAF \Longleftrightarrow \|\cdot\| \models_{\mathcal{PL}_3^L} ec_{HS}(HSAF).
	\end{equation*}
\end{thm}

\begin{proof} 
	We need to check that for a given assignment $\|\cdot\|$ and $\forall \beta\in A^{BH(0)}\cup R^{BH(n)}$, any value of $\beta$ satisfies complete semantics iff
	\begin{equation}\label{eqbeta}
		\|\beta\leftrightarrow\bigwedge_{(S_i, \beta)\in \mathsf{R}^{HS(n)}}\neg(r_{S_i}^{\beta}\wedge\bigwedge_{b_{ij}\in S_i}b_{ij})\|=1
	\end{equation}
	 holds in the $\mathcal{PL}_3^L$. We need to discuss three cases.
	\begin{itemize}
		\item Case 1, $\|\beta\|=1$. 
		\\$\|\beta\|=1$ by model $\|\cdot\|$ under complete semantics
		\\$\Longleftrightarrow$ $\forall (S_i, \beta)\in \mathsf{R}^{HS(n)}$ $\exists b_{ij} \in S_i$ s.t. $\|b_{ij}\|= 0$ or $\|(S_i, \beta)\| = 0$  by complete semantics
		\\$\Longleftrightarrow$ $\forall (S_i, \beta)\in \mathsf{R}^{HS(n)}$, $\|\bigwedge_{b_{ij}\in S_i}b_{ij}\|= 0$ or $\| r_{S_i}^{\beta}\| = 0$ in the $\mathcal{PL}_3^L$
		\\$\Longleftrightarrow$ $\forall (S_i, \beta)\in \mathsf{R}^{HS(n)}$, $\|\neg\bigwedge_{b_{ij}\in S_i}b_{ij}\|= 1$ or $\|\neg r_{S_i}^{\beta}\|=1$ in the $\mathcal{PL}_3^L$
		\\$\Longleftrightarrow$ $\forall (S_i, \beta)\in \mathsf{R}^{HS(n)}$, $\|\neg r_{S_i}^{\beta}\vee \neg\bigwedge_{b_{ij}\in S_i}b_{ij}\|= 1$ in the $\mathcal{PL}_3^L$		
		\\$\Longleftrightarrow$ $\forall (S_i, \beta)\in \mathsf{R}^{HS(n)}$, $\|\neg (r_{S_i}^{\beta}\wedge\bigwedge_{b_{ij}\in S_i}b_{ij})\|= 1$ in the $\mathcal{PL}_3^L$
		\\$\Longleftrightarrow$ $\|\bigwedge_{(S_i, \beta)\in \mathsf{R}^{HS(n)}}\neg(r_{S_i}^{\beta}\wedge\bigwedge_{b_{ij}\in S_i}b_{ij})\|=1$  in the $\mathcal{PL}_3^L$
		\\$\Longleftrightarrow$ Equation \ref{eqbeta} holds in the $\mathcal{PL}_3^L$.
		
		\item Case 2, $\|\beta\|=0$. 
		\\$\|\beta\|=0$ by model $\|\cdot\|$ under complete semantics
		\\$\Longleftrightarrow$ $\exists (S_i, \beta)\in \mathsf{R}^{HS(n)}$, $\forall$ $b_{ij} \in S_i$, $\|b_{ij}\| = 1$ by complete semantics
		\\$\Longleftrightarrow$ $\exists (S_i, \beta)\in \mathsf{R}^{HS(n)}$, $\|\bigwedge_{b_{ij}\in S_i}b_{ij}\|=1$ in the $\mathcal{PL}_3^L$
		\\$\Longleftrightarrow$ $\exists (S_i, \beta)\in \mathsf{R}^{HS(n)}$, $\|\neg\bigwedge_{b_{ij}\in S_i}b_{ij}\|=0$ in the $\mathcal{PL}_3^L$
		\\$\Longleftrightarrow$ $\|\bigwedge_{(S_i, \beta)\in \mathsf{R}^{HS(n)}}(\neg \bigwedge_{b_{ij}\in S_i}b_{ij})\|= 0$ in the $\mathcal{PL}_3^L$
		\\$\Longleftrightarrow$ Equation \ref{eqbeta} holds in the $\mathcal{PL}_3^L$.
		
		\item Case 3, $\|\beta\|=\frac{1}{2}$. 
		\\$\|\beta\|=\frac{1}{2}$ by model $\|\cdot\|$ under complete semantics
		\\$\Longleftrightarrow$ $\|\beta\|\neq1$ and $\|\beta\|\neq0$ under complete semantics
		\\$\Longleftrightarrow$ $\|\bigwedge_{(S_i, \beta)\in \mathsf{R}^{HS(n)}}(\neg \bigwedge_{b_{ij}\in S_i}b_{ij})\|$ is neither 1 nor 0 in the $\mathcal{PL}_3^L$ by Case 1 and Case 2
		\\$\Longleftrightarrow$ $\|\bigwedge_{(S_i, \beta)\in \mathsf{R}^{HS(n)}}(\neg \bigwedge_{b_{ij}\in S_i}b_{ij})\|=\frac{1}{2}$ in the $\mathcal{PL}_3^L$.
		\\$\Longleftrightarrow$ Equation \ref{eqbeta} holds in the $\mathcal{PL}_3^L$.
	\end{itemize}
	From the three cases above, the value of  any $\beta$ satisfies complete semantics iff $\|\beta\leftrightarrow\bigwedge_{(S_i, \beta)\in \mathsf{R}^{HS(n)}}\neg(r_{S_i}^{\beta}\wedge\bigwedge_{b_{ij}\in S_i}b_{ij})\|=1$ in the $\mathcal{PL}_3^L$. 
	So, an assignment $\|\cdot\|$ of the $HSAF$ satisfies complete semantics iff $\|\bigwedge_{\beta\in \mathbf{T}(\mathsf{R}^{HS(n)})}(\beta\leftrightarrow\bigwedge_{(S_i, \beta)\in \mathsf{R}^{HS(n)}}\neg(r_{S_i}^{\beta}\wedge\bigwedge_{b_{ij}\in S_i}b_{ij}))\|=1$, i.e., $\|\cdot\|$ is a model of the $ec_{HS}(HSAF)$, in the $\mathcal{PL}_3^L$.
\end{proof} 
We continue to use the $HSAF^{(2)}$ in Example \ref{ex:hs} to illustrate this theorem.
\addtocounter{exmp}{-1}
\begin{exmp}[continued]\label{ex:hs-cont2}
	For the given $HSAF^{(2)}$, we take the similar processes in
	Example \ref{ex:2-cont}. Then in the $\mathcal{PL}_3^L$ we have that
	\begin{equation*}
	\|ec_{HS}(HSAF^{(2)})\|=1
	\Longleftrightarrow
		\begin{cases}
			\|b\|=\|\varepsilon\|=\|\gamma\|=\|\delta\|=\|\epsilon\|=1\\ \|d\|=\frac{1}{2}\\
			\|a\|=\|\alpha\|=1-\|\beta\|=1-\min\{\|\beta\|,\|c\|\}\\			
			\|c\|=\min\{\frac{1}{2}, 1-\|a\|\}
		\end{cases}.
	\end{equation*}
	For the above equation over the set $\{0,1, \frac{1}{2}\}$, if $\|a\|=0$, then it leads to a contradiction. If $\|a\|=\frac{1}{2}$, then we obtain a solution $\{\|b\|=\|\varepsilon\|=\|\gamma\|=\|\delta\|=\|\epsilon\|=1, \|d\|=\|a\|=\|\alpha\|=\|\beta\|=\|c\|=\frac{1}{2}\}$. If $\|a\|=1$, then we obtain the other solution $\{\|b\|=\|\varepsilon\|=\|\gamma\|=\|\delta\|=\|\epsilon\|=\|a\|=\|\alpha\|=1, \|d\|=\frac{1}{2}, \|\beta\|=\|c\|=0\}$.
	Thus, the models of $ec_{HS}(HSAF^{(2)})$ in the $\mathcal{PL}_3^L$ are the same as models of the $HSAF^{(2)}$ under complete semantics.
\end{exmp}

By using the notion $\mathcal{HSAF}$ to denote the set of all $HSAF$s, from Theorem \ref{thm-pl3}, we obtain the following corollary.
\begin{cor}
The	$ec_{HS}$ is a translation of $\mathcal{HSAF}$ associated with $\mathcal{PL}_3^L$ under complete semantics and thus complete semantics of $\mathcal{HSAF}$ is translatable.
\end{cor}

Now, using complete semantics (or equivalently, the encoding approach), we formally model the first practical argumentation scenario presented at the beginning of Section \ref{sec6}. With ${A}^{HS(0)}=\{A,B,C,P\}$, ${R}^{HS(0)}=\{(\{C\},B),$ $(\{B\},A)\}=\{r_{\{C\}}^B, r_{\{B\}}^A\}=\{\alpha, r_{\{B\}}^A\}$, and ${R}^{HS(1)}=\{\alpha, r_{\{B\}}^A, r_{\{A, \alpha\}}^P\}$, we construct the $HSAF=(\mathsf{A}^{HS(0)}, \mathsf{R}^{HS(1)})$, where $\mathsf{A}^{HS(0)}=\{A,B,C,P,\bot\}$ and $\mathsf{R}^{HS(1)}=\{\alpha, r_{\{B\}}^A, r_{\{A, \alpha\}}^P, (\{\bot\},\alpha), (\{\bot\},r_{\{B\}}^A), (\{\bot\},r_{\{A, \alpha\}}^P)\}$. The only model of this $HSAF$ under complete semantics is $\{\|A\|=\|C\|=\|\alpha\|=\|r_{\{B\}}^A\|= \|r_{\{A, \alpha\}}^P\|= \|(\{\bot\},\alpha)\|= \|(\{\bot\},r_{\{B\}}^A)\|= \|(\{\bot\},r_{\{A, \alpha\}}^P)\|=1, \|B\|$ $=\|P\|=0\}$. The result that both $P$ and $B$ are defeated aligns with our intuitive conclusion.
\subsection{Three equational semantics of $HSAF$s}
In this subsection, we emphasize three equational systems on [0,1] for $HSAF$s and then we prove the theorems of model equivalence. In the following discussion, for a given $HSAF=(\mathsf{A}^{HS(0)}, \mathsf{R}^{HS(n)})$ and $\beta\in \mathbf{T}(\mathsf{R}^{HS(n)})$, let $n_\beta= |\{S_i \mid r_{S_i}^{\beta} \in  \mathsf{R}^{HS(n)}\}|$ denote the number of set attackers of $\beta$.  For  $r_{S_i}^{\beta}\in \mathsf{R}^{HS(n)}$ and $b_{ij}\in S_i$, let $k_i= |S_i|=|\{b_{ij} \mid b_{ij} \in S_i\}|$ denote the cardinality of the set attacker $S_i$. The encoding function of the $HSAF$ w.r.t. a $\mathcal{PL}_{[0, 1]}$ is still $ec_{HS}$ such that 
\begin{equation*}
	ec_{HS}(HSAF)=\bigwedge_{\beta\in \mathbf{T}(\mathsf{R}^{HS(n)})}(\beta\leftrightarrow\bigwedge_{r_{S_i}^{\beta}\in \mathsf{R}^{HS(n)}}\neg(r_{S_i}^{\beta}\wedge\bigwedge_{b_{ij}\in S_i}b_{ij})).
\end{equation*}
Then we have three theorems of model equivalence corresponding to the three typical equational systems.
\subsubsection{The equational system $Eq_G^{HS}$}
For a $HSAF=(\mathsf{A}^{HS(0)}, \mathsf{R}^{HS(n)})$, we give the equational system $Eq_G^{HS}$ such that $\|\bot\|=0$, $\|r_{\bot}^\alpha\|=1$ (if $r_{\bot}^\alpha\in\mathsf{R}^{HS(n)}$), and the equation for each $\beta\in \mathbf{T}(\mathsf{R}^{HS(n)})$ is
\begin{equation}\label{eq8}
	\|\beta\|=\min_{i=1}^{n_\beta}\max\{1-\|r_{S_i}^{\beta}\|, 1-\|b_{i1}\|, 1-\|b_{i2}\|,\dots,1-\|b_{ik_i}\|\}.
\end{equation}
The equational semantics $Eq_G^{HS}$ for $HSAF$s evaluates argument acceptability using a min-max approach based on G\"{o}del logic. For each element $\beta$ in the framework, its truth value is determined by the most effective attack among all its set attackers. Specifically, for each set attacker $S_i$, we compute the maximum of the complement values of both the attack $r_{S_i}^\beta$ and all arguments in $S_i$, representing the strongest defensive component against that particular attack. The overall value of $\beta$ is then taken as the minimum across all these maxima, ensuring that $\beta$ is fully accepted only when all attacking actions are completely ineffective.

Denote that an assignment $\|\cdot\|$ is a model of an $HSAF$ under equational semantics $Eq_G^{HS}$ by $\|\cdot\| \models_{\mathcal{HS}_{[0,1]}^G} HSAF$. Then, we present the model equivalence theorem.
\begin{thm}\label{thm14}
	For an $HSAF=(\mathsf{A}^{HS(0)}, \mathsf{R}^{HS(n)})$ and an assignment $\|\cdot\|: \mathsf{A}^{HS(0)}\cup\mathsf{R}^{HS(n)}\to [0,1]$, 
	\begin{equation*}
		\|\cdot\| \models_{\mathcal{HS}_{[0,1]}^G} HSAF \Longleftrightarrow \|\cdot\| \models_{\mathcal{PL}_{[0,1]}^G} ec_{HS}(HSAF).
	\end{equation*}
\end{thm}

\begin{proof} 
	A model $\|\cdot\|$ of $ec_{HS}(HSAF)$ in the $\mathcal{PL}_{[0, 1]}^G$
	\\$\Longleftrightarrow $ (a solution of)
	$\|\bigwedge_{ \beta\in \mathbf{T}(\mathsf{R}^{HS(n)})}(\beta\leftrightarrow\bigwedge_{(S_i, \beta)\in \mathsf{R}^{HS(n)}}\neg(r_{S_i}^{\beta}\wedge\bigwedge_{b_{ij}\in S_i}b_{ij}))\|=1$
	\\$\Longleftrightarrow $ for each $\beta\in \mathbf{T}(\mathsf{R}^{HS(n)})$, $\|\beta\leftrightarrow\bigwedge_{(S_i, \beta)\in \mathsf{R}^{HS(n)}}\neg(r_{S_i}^{\beta}\wedge\bigwedge_{b_{ij}\in S_i}b_{ij})\|=1$
	\\$\Longleftrightarrow $ for each $\beta\in \mathbf{T}(\mathsf{R}^{HS(n)})$, $\|\beta\|=\|\bigwedge_{(S_i, \beta)\in \mathsf{R}^{HS(n)}}\neg(r_{S_i}^{\beta}\wedge\bigwedge_{b_{ij}\in S_i}b_{ij})\|$ by Lemma \ref{remark2}
	\\$\Longleftrightarrow $ for each $\beta\in \mathbf{T}(\mathsf{R}^{HS(n)})$,
	\begin{align*}
		\|\beta\|&=\min_{i=1}^{n_\beta}\|\neg(r_{S_i}^{\beta}\wedge\bigwedge_{b_{ij}\in S_i}b_{ij})\|
		\\&=\min_{i=1}^{n_\beta}(1-\|r_{S_i}^{\beta}\wedge\bigwedge_{b_{ij}\in S_i}b_{ij}\|)
		\\&=\min_{i=1}^{n_\beta}(1-\min\{\|r_{S_i}^{\beta}\|, \|b_{i1}\|, \|b_{i2}\|,\dots,\|b_{ik_i}\|\})
		\\&=\min_{i=1}^{n_\beta}\max\{1-\|r_{S_i}^{\beta}\|, 1-\|b_{i1}\|, 1-\|b_{i2}\|,\dots,1-\|b_{ik_i}\|\}
	\end{align*}
	\\$\Longleftrightarrow$ a model $\|\cdot\|$ of the $HSAF$ under equational semantics $Eq_G^{HS}$.
\end{proof} 
\begin{cor}
	The	$ec_{HS}$ is a translation of $\mathcal{HSAF}$ associated with $\mathcal{PL}_{[0,1]}^G$ under the equational semantics $Eq_G^{HS}$.
\end{cor}
\begin{cor}
The equational semantics $Eq_G^{HS}$ is a fuzzy normal encoded semantics and thus it is translatable.
\end{cor}
We continue Example \ref{ex:hs} to illustrate Theorem \ref{thm14} below.
\addtocounter{exmp}{-1}
\begin{exmp}[continued]\label{ex:hs-cont3}
	According to Equation \ref{eq8}, the system of equations for the $HSAF$ is that	 
	\begin{equation}\label{eq9}
		\begin{cases}
		\|b\|=\|\varepsilon\|=\|\gamma\|=\|\delta\|=\|\epsilon\|=1\\ \|d\|=\frac{1}{2}\\
		\|a\|=\|\alpha\|=1-\|\beta\|=1-\min\{\|\beta\|,\|c\|\}\\			
		\|c\|=\min\{\frac{1}{2}, 1-\|a\|\}
		\end{cases}.
	\end{equation}
	Thus, we have models $\{\|b\|=\|\varepsilon\|=\|\gamma\|=\|\delta\|=\|\epsilon\|=1, 
	\|d\|=\frac{1}{2},
	\|a\|=\|\alpha\|=x,
	\|\beta\|=	\|c\|=1-x\}$, where $x\in [\frac{1}{2},1]$.
	
	Next, we find models of $ec_{HS}(HSAF^{(2)})$ in the $\mathcal{PL}_{[0, 1]}^G$, i.e., all assignments that satisfy $\|ec_{HS}(HSAF^{(2)})\|=1$. In the $\mathcal{PL}_{[0, 1]}^G$ we simplify the equation $\|ec_{HS}(HSAF^{(2)})\|=1$ and then obtain an equation that is the same as Equation \ref{eq9}.
	Thus, the models of $ec_{HS}(HSAF^{(2)})$ in the $\mathcal{PL}_{[0, 1]}^G$ are the same as models of the $HSAF^{(2)}$ under equational semantics $Eq_G^{HS}$. 
\end{exmp}
\subsubsection{The equational system $Eq_P^{HS}$}
For a $HSAF=(\mathsf{A}^{HS(0)}, \mathsf{R}^{HS(n)})$, we give the equational system $Eq_P^{HS}$ such that $\|\bot\|=0$, $\|r_{\bot}^\alpha\|=1$ (if $r_{\bot}^\alpha\in\mathsf{R}^{HS(n)}$), and the equation for each $\beta\in \mathbf{T}(\mathsf{R}^{HS(n)})$ is
\begin{equation}\label{eq10}
	\|\beta\|=\prod_{i=1}^{n_\beta}(1-\|r_{S_i}^{\beta}\|\prod_{j=1}^{k_i}\|b_{ij}\|).
\end{equation}
The equational semantics $Eq_P^{HS}$ for $HSAF$s employs a multiplicative approach based on Product logic to evaluate argument acceptability. For each element $\beta$, its truth value is computed as the product of failure probabilities across all set attackers. Each term $1 - \|r_{S_i}^{\beta}\|\prod_{j=1}^{k_i}\|b_{ij}\|$ represents the probability that set attacker $S_i$ fails, calculated by subtracting the joint effectiveness of the attack and all arguments in $S_i$ from 1. The final acceptability $\|\beta\|$ equals the product of these failure probabilities, meaning $\beta$ is fully accepted only when all attacking actions completely fail, with partial acceptance decreasing multiplicatively as attacks gain strength.

Denote that an assignment $\|\cdot\|$ is a model of an $HSAF$ under equational semantics $Eq_P^{HS}$ by $\|\cdot\| \models_{\mathcal{HS}_{[0,1]}^P} HSAF$. Then, we present the model equivalence theorem.
\begin{thm}\label{thm15}
	For an $HSAF=(\mathsf{A}^{HS(0)}, \mathsf{R}^{HS(n)})$ and an assignment $\|\cdot\|: \mathsf{A}^{HS(0)}\cup\mathsf{R}^{HS(n)}\to [0,1]$, 
	\begin{equation*}
		\|\cdot\| \models_{\mathcal{HS}_{[0,1]}^P} HSAF \Longleftrightarrow \|\cdot\| \models_{\mathcal{PL}_{[0,1]}^P} ec_{HS}(HSAF).
	\end{equation*}
\end{thm}

\begin{proof} 
	A model $\|\cdot\|$ of $ec_{HS}(HSAF)$ in the $\mathcal{PL}_{[0, 1]}^P$
	\\$\Longleftrightarrow $ (a solution of)
	$\|\bigwedge_{\beta\in \mathbf{T}(\mathsf{R}^{HS(n)})}(\beta\leftrightarrow\bigwedge_{(S_i, \beta)\in \mathsf{R}^{HS(n)}}\neg(r_{S_i}^{\beta}\wedge\bigwedge_{b_{ij}\in S_i}b_{ij}))\|=1$
	\\$\Longleftrightarrow $ for each $\beta\in \mathbf{T}(\mathsf{R}^{HS(n)})$, $\|\beta\leftrightarrow\bigwedge_{(S_i, \beta)\in \mathsf{R}^{HS(n)}}\neg(r_{S_i}^{\beta}\wedge\bigwedge_{b_{ij}\in S_i}b_{ij})\|=1$
	\\$\Longleftrightarrow $ for each $\beta\in \mathbf{T}(\mathsf{R}^{HS(n)})$, $\|\beta\|=\|\bigwedge_{(S_i, \beta)\in \mathsf{R}^{HS(n)}}\neg(r_{S_i}^{\beta}\wedge\bigwedge_{b_{ij}\in S_i}b_{ij})\|$ by Lemma \ref{remark2}
	\\$\Longleftrightarrow $ for each $\beta\in \mathbf{T}(\mathsf{R}^{HS(n)})$, $\|\beta\|=\prod_{i=1}^{n_\beta}\|\neg(r_{S_i}^{\beta}\wedge\bigwedge_{b_{ij}\in S_i}b_{ij})\|=\prod_{i=1}^{n_\beta}(1-\|r_{S_i}^{\beta}\wedge\bigwedge_{b_{ij}\in S_i}b_{ij}\|)=\prod_{i=1}^{n_\beta}(1-\|r_{S_i}^{\beta}\|\prod_{j=1}^{k_i}\|b_{ij}\|)$
	\\$\Longleftrightarrow$ a model $\|\cdot\|$ of the $HSAF$ under equational semantics $Eq_P^{HS}$.
\end{proof} 
\begin{cor}
	The	$ec_{HS}$ is a translation of $\mathcal{HSAF}$ associated with $\mathcal{PL}_{[0,1]}^P$ under the equational semantics $Eq_P^{HS}$.
\end{cor}
\begin{cor}
The equational semantics $Eq_P^{HS}$ is a fuzzy normal encoded semantics and thus it is translatable.
\end{cor}
We illustrate Theorem \ref{thm15} by the continued example below.
\addtocounter{exmp}{-1}
\begin{exmp}[continued]\label{ex:hs-cont4}
	According to Equation \ref{eq10}, the system of equations for the $HSAF$ is that	 
	\begin{equation}\label{eq11}
		\begin{cases}
			\|b\|=\|\varepsilon\|=\|\gamma\|=\|\delta\|=\|\epsilon\|=1\\ \|d\|=\frac{1}{2}\\
			\|a\|=1-\|\beta\|\\	
			\|\alpha\|=1-\|\beta\|\|c\|\\
			\|\beta\|=1-\|\alpha\|\\
			\|c\|=\frac{1}{2}(1-\|a\|\|\alpha\|)
		\end{cases}.
	\end{equation}
	We have a unique solution $\{\|b\|=\|\varepsilon\|=\|\gamma\|=\|\delta\|=\|\epsilon\|=\|a\|=\|\alpha\|=1, \|\beta\|=\|c\|=0\}$ for \ref{eq11}, i.e., we have only one model under equational semantics $Eq_P^{HS}$.
	
	When we simplify the equation $\|ec_{HS}(HSAF^{(2)})\|=1$ in the $\mathcal{PL}_{[0, 1]}^P$, we obtain an equation that is the same as Equation \ref{eq11}.
	Thus, we have the same model result of the $HSAF^{(2)}$ under the equational semantics $Eq_P^{HS}$ or the fuzzy normal encoded semantics. 
\end{exmp}
\subsubsection{The equational system $Eq_L^{HS}$}
For a $HSAF=(\mathsf{A}^{HS(0)}, \mathsf{R}^{HS(n)})$, we give the equational system $Eq_L^{HS}$ such that $\|\bot\|=0$, $\|r_{\bot}^\alpha\|=1$ (if $r_{\bot}^\alpha\in\mathsf{R}^{HS(n)}$), and the equation for each $\beta\in \mathbf{T}(\mathsf{R}^{HS(n)})$ is
\begin{equation}\label{eq12}
	\|\beta\|=\begin{cases}
		0 & \quad \sum_{i=1}^{n_\beta}y_i\leqslant n_\beta-1\\
		\sum_{i=1}^{n_\beta}y_i-n_\beta+1 & \quad \sum_{i=1}^{n_\beta}y_i> n_\beta-1
	\end{cases}
\end{equation}
where
\begin{equation*}
	y_i=\begin{cases}
		1 & \quad \|r_{S_i}^{\beta}\|+\sum_{j=1}^{k_i}\|b_{ij}\|\leqslant k_i\\
		k_i+1-\|r_{S_i}^{\beta}\|-\sum_{j=1}^{k_i}\|b_{ij}\| & \quad \|r_{S_i}^{\beta}\|+\sum_{j=1}^{k_i}\|b_{ij}\|> k_i
	\end{cases}.
\end{equation*}
The equational semantics $Eq_L^{HS}$ for $HSAF$s employs a threshold-based mechanism grounded in {\L}ukasiewicz logic. For each element $\beta$, we first compute the ineffectiveness $y_i$ of each set attacker $S_i$, which equals 1 when the combined strength of the attack and its members is sufficiently low ($\leq k_i$), and decreases linearly otherwise. The final acceptability $\|\beta\|$ is then determined by comparing the total ineffectiveness $\sum y_i$ against threshold $n_\beta-1$: $\|\beta\| = 0$ if the threshold is unmet, and increases linearly with the excess otherwise. This design captures synergistic attack effects where $\beta$ gains acceptability only when a critical mass of attackers fails, enabling smooth transitions in argument status.

Denote that an assignment $\|\cdot\|$ is a model of an $HSAF$ under equational semantics $Eq_L^{HS}$ by $\|\cdot\| \models_{\mathcal{HS}_{[0,1]}^L} HSAF$. Then, we present the model equivalence theorem.
\begin{thm}\label{thm16}
	For an $HSAF=(\mathsf{A}^{HS(0)}, \mathsf{R}^{HS(n)})$ and an assignment $\|\cdot\|: \mathsf{A}^{HS(0)}\cup\mathsf{R}^{HS(n)}\to [0,1]$, 
	\begin{equation*}
		\|\cdot\| \models_{\mathcal{HS}_{[0,1]}^L} HSAF \Longleftrightarrow \|\cdot\| \models_{\mathcal{PL}_{[0,1]}^L} ec_{HS}(HSAF).
	\end{equation*}
\end{thm}

\begin{proof} 
	A model $\|\cdot\|$ of $ec_{HS}(HSAF)$ in the $\mathcal{PL}_{[0, 1]}^L$
	\\$\Longleftrightarrow $ (a solution of)
	$\|\bigwedge_{\beta\in \mathbf{T}(\mathsf{R}^{HS(n)})}(\beta\leftrightarrow\bigwedge_{(S_i, \beta)\in \mathsf{R}^{HS(n)}}\neg(r_{S_i}^{\beta}\wedge\bigwedge_{b_{ij}\in S_i}b_{ij}))\|=1$
	\\$\Longleftrightarrow $ for each $\beta\in \mathbf{T}(\mathsf{R}^{HS(n)})$, $\|\beta\leftrightarrow\bigwedge_{(S_i, \beta)\in \mathsf{R}^{HS(n)}}\neg(r_{S_i}^{\beta}\wedge\bigwedge_{b_{ij}\in S_i}b_{ij})\|=1$
	\\$\Longleftrightarrow $ for each $\beta\in \mathbf{T}(\mathsf{R}^{HS(n)})$, $\|\beta\|=\|\bigwedge_{(S_i, \beta)\in \mathsf{R}^{HS(n)}}\neg(r_{S_i}^{\beta}\wedge\bigwedge_{b_{ij}\in S_i}b_{ij})\|$ by Lemma \ref{remark2}
	\\$\Longleftrightarrow $ for each $\beta\in \mathbf{T}(\mathsf{R}^{HS(n)})$, by the formula \ref{n-t-norm},
	\begin{equation*}
		\|\beta\|=
		\begin{cases}
			0 & \quad \sum_{i=1}^{n_\beta}y_i \leq n_\beta-1\\
			\sum_{i=1}^{n_\beta}y_i-n_\beta+1 & \quad \sum_{i=1}^{n_\beta}y_i > n_\beta-1
		\end{cases}
	\end{equation*}
	where $y_i=\|\neg(r_{S_i}^{\beta}\wedge\bigwedge_{b_{ij}\in S_i}b_{ij})\|=1-\|r_{S_i}^{\beta}\wedge\bigwedge_{b_{ij}\in S_i}b_{ij}\|$.
	By the formula \ref{n-t-norm} again, we have 
	\begin{equation*}
		\|r_{S_i}^{\beta}\wedge\bigwedge_{b_{ij}\in S_i}b_{ij}\|=
		\begin{cases}
			0 & \quad \|r_{S_i}^{\beta}\|+\sum_{j=1}^{k_i}\|b_{ij}\| \leq k_i\\
			\|r_{S_i}^{\beta}\|+\sum_{j=1}^{k_i}\|b_{ij}\|-k_i & \quad \|r_{S_i}^{\beta}\|+\sum_{j=1}^{k_i}\|b_{ij}\| > k_i
		\end{cases},
	\end{equation*}
	i.e., 
	\begin{equation*}
		y_i=
		\begin{cases}
			1 & \quad \|r_{S_i}^{\beta}\|+\sum_{j=1}^{k_i}\|b_{ij}\| \leq k_i\\
			k_i+1-\|r_{S_i}^{\beta}\|-\sum_{j=1}^{k_i}\|b_{ij}\| & \quad \|r_{S_i}^{\beta}\|+\sum_{j=1}^{k_i}\|b_{ij}\| > k_i
		\end{cases}.
	\end{equation*}

The equations above are the same as the equational system $Eq_L^{HS}$. So $\|\cdot\|$ is a model of $ec_{HS}(HSAF)$ in the $\mathcal{PL}_{[0, 1]}^L$ iff it is
a model of the $HSAF$ under equational semantics $Eq_L^{HS}$, i.e., we have proven this theorem.
\end{proof} 
\begin{cor}
The	$ec_{HS}$ is a translation of $\mathcal{HSAF}$ associated with $\mathcal{PL}_{[0,1]}^L$ under the equational semantics $Eq_L^{HS}$.
\end{cor}
\begin{cor}
The equational semantics $Eq_L^{HS}$ is a fuzzy normal encoded semantics and thus it is translatable.
\end{cor}
We use a continued example to illustrate the equational semantics $Eq_L^{HS}$ and Theorem \ref{thm16}.
\addtocounter{exmp}{-1}
\begin{exmp}[continued]\label{ex:hs-cont5}
	According to Equation \ref{eq12}, the system of equations for the $HSAF$ is listed and simplified as	 
	\begin{equation}\label{eq13}
		\begin{cases}
			\|b\|=\|\varepsilon\|=\|\gamma\|=\|\delta\|=\|\epsilon\|=1\\ \|d\|=\frac{1}{2}\\
			\|a\|=\|\alpha\|=1-\|\beta\|\\	
			\|a\|=\begin{cases}
				1 &c\leqslant a\\
				1+a-c &c>a
			\end{cases}\\		
			\|c\|=\begin{cases}
				\frac{1}{2} &a\in[0, \frac{1}{2}]\\
				\frac{3}{2}-2a &a\in(\frac{1}{2}, \frac{3}{4})\\
				0 &a\in[\frac{3}{4}, 1]
			\end{cases}
		\end{cases}.
	\end{equation}
	If $c\leqslant a$ we have $\|a\|=1$ and $\|c\|=0$. If $c>a$, there is no solution. Thus, we have only one model under the equational semantics $Eq_L^{HS}$.
	
	If we simplify the equation $\|ec_{HS}(HSAF^{(2)})\|=1$ in the $\mathcal{PL}_{[0, 1]}^L$, we obtain an equation that is the same as Equation \ref{eq13}, i.e., the model under fuzzy normal encoded semantics w.r.t. $\mathcal{PL}_{[0, 1]}^L$ is the same as the model under the equational semantics $Eq_L^{HS}$.
\end{exmp}
\subsection{The formal equational approach for $HSAF$s}\label{subsec6.5}
We have presented three important equational systems for each kind of $AFSA$. Now we construct the formal equational approach for $HSAF$s following the approach in \cite{gabbay2012equational}. 
\begin{defn}\label{defn27}
	A \emph{real equation function tuple} over the real interval $[0, 1]$  is a ($k+1$)-tuple $\mathcal{H}=(h, h_1, h_2, \dots, h_i, \dots, h_k)$ such that\\
 (a) $h$ is a continuous function $h : [0, 1]^k \rightarrow [0, 1]$ with $k$ variables $\{x_1, \dots, x_k \}$ satisfying
	  \begin{itemize}
	  	\item $h(1, ... , 1) = 1$;
	  	\item $h(x_1, \dots , 0, \dots , x_k ) = 0$;
	  	\item  $h(x_1, \dots , x_k ) = h(y_1, \dots , y_k )$, where $(x_1, \dots , x_k )$ and $(y_1, \dots , y_k )$ are permutations of each other.
	  \end{itemize}
	(b) $h_i (i\in\{1,2, \dots, k\})$ is a continuous function $h_i : [0, 1]^{j_i+1} \rightarrow [0, 1]$ with $j_i+1$ variables $\{u_{i1}, u_{i2}, \dots, u_{ij_i}, u_{ij_i+1}\}$ satisfying
	\begin{itemize}
		\item $x_i=h_i(u_{i1}, u_{i2}, \dots, u_{ij_i}, u_{ij_i+1})$ for each variable $x_i$ of $h$;
		\item $h_i(1, ... , 1) = 0$;
		\item $h_i(u_{i1}, \dots, 0, \dots , u_{ij_i}, u_{ij_i+1}) = 1$;
		\item $h_i(u_{i1}, \dots, u_{ij_i}, u_{ij_i+1}) = h_i(v_{i1}, \dots, v_{ij_i}, v_{ij_i+1})$, where $(u_{i1}, \dots, u_{ij_i}, u_{ij_i+1})$ and $(v_{i1}, \dots, v_{ij_i}, v_{ij_i+1})$ are permutations of each other.
		\end{itemize}
		(c) $h$ is non-decreasing w.r.t. each variable in $\{x_1, \dots, x_k\}$, and each $h_i$ ($i \in \{1, 2, \dots, k\}$) is non-increasing w.r.t. each variable in $\{u_{i1}, u_{i2}, \dots, u_{ij_i}, u_{ij_i+1}\}$.
\end{defn}
Definition \ref{defn27} introduces a real equation function tuple $\mathcal{H}=(h, h_1, \dots, h_k)$ that provides the mathematical foundation for continuous argument evaluation. The main function $h$ aggregates $k$ input values through a continuous, symmetric function that preserves boundary conditions (yielding 1 for all-1 inputs and 0 for any zero input). Each component function $h_i$ processes the elements of a set attacker $S_i$, mapping the truth values of the attack and set members to an intermediate value $x_i$, with complementary boundary behavior: it outputs 0 when all inputs are 1 and 1 when any input is 0, while maintaining symmetry across inputs.
\begin{defn}\label{defn28}
 A \emph{general real equational system} $Eq^{HS}$ for an $HSAF = (\mathsf{A}^{HS(0)}, \mathsf{R}^{HS(n)})$ over $[0, 1]$ is defined as follows. Let $\|\bot\|=0$, $\|r_{\{\bot\}}^\alpha\|=1$ (if $r_{\{\bot\}}^\alpha\in\mathsf{R}^{HS(n)}$), and the equation for each $\beta\in \mathbf{T}(\mathsf{R}^{HS(n)})$ be
   \begin{equation*}
   	\|\beta\| = h_\beta(h_1(\|b_{11}\|,\dots,\|b_{1j_1}\|, \|r_{S_1}^\beta\|), \dots,  h_k(\|b_{k1}\|,\dots,\|b_{kj_k}\|, \|r_{S_k}^\beta\|))
   \end{equation*}
   where $(h_\beta, h_1, \dots, h_k)$ is a real equation function ($k+1$)-tuple and for each $i \in \{1, \dots, k\}$, the set $S_i = \{b_{i1}, \dots, b_{ij_i}\}$ is a set attacker of $\beta$.
\end{defn}
Definition \ref{defn28} constructs a general real equational system for $HSAF$s using the function tuple from Definition \ref{defn27}. For each element $\beta$ in the framework, its acceptability is computed by first applying component functions $h_1, \dots, h_k$ to each set attacker $S_i$ (processing both the attack $r_{S_i}^\beta$ and all elements in $S_i$), then aggregating these intermediate results through the main function $h_\beta$. This creates a unified framework where different t-norms and aggregation operators can be instantiated via specific function tuples, enabling flexible modeling of argument strength under various logical interpretations while ensuring mathematical well-foundedness through continuity and boundary preservation.

For a given $HSAF=(\mathsf{A}^{HS(0)}, \mathsf{R}^{HS(n)})$, let $Eq_{ec}^{HS}$ be an equational system which is yielded by the model of $ec_{HS}(HSAF)$ in the $\mathcal{PL}_{[0, 1]}^{\star}$, where the $\mathcal{PL}_{[0, 1]}^\star$ is a fuzzy propositional logic system equipped with continuous negation $N^\star$, continuous  t-norm $T$ and $R$-implication $I_T$. We also call the $Eq_{ec}^{HS}$ as a \emph{continuous fuzzy normal encoded} ($CFNE$) equational system. By the following theorem, we get the relationship between general real equational systems and $CFNE$ equational systems.
\begin{thm}\label{thm17}
	 For a given $HSAF$, each $CFNE$ equational system $Eq_{ec}^{HS}$ is a general real equational system $Eq^{HS}$.
\end{thm}
\begin{proof}
	Firstly, from
	\begin{equation*}
		\|ec_{HS}(HSAF)\|=\|\bigwedge_{\beta\in \mathbf{T}(\mathsf{R}^{HS(n)})}(\beta\leftrightarrow\bigwedge_{(S_i, \beta)\in \mathsf{R}^{HS(n)}}\neg(r_{S_i}^{\beta}\wedge\bigwedge_{b_{ij}\in S_i}b_{ij}))\|=1,
	\end{equation*}
 we have that for each $\beta\in \mathbf{T}(\mathsf{R}^{HS(n)})$, 
 \begin{equation*}
 	\|\beta\leftrightarrow\bigwedge_{(S_i, \beta)\in \mathsf{R}^{HS(n)}}\neg(r_{S_i}^{\beta}\wedge\bigwedge_{b_{ij}\in S_i}b_{ij})\|=1.
 \end{equation*}
	By Lemma \ref{remark2}, we have  $\|\beta\|=\|\bigwedge_{(S_i, \beta)\in \mathsf{R}^{HS(n)}}\neg(r_{S_i}^{\beta}\wedge\bigwedge_{b_{ij}\in S_i}b_{ij})\|$. Thus, by continuous negation $N^\star$ and continuous  t-norm $T$, the equational system $Eq_{ec}^{HS}$ yielded by the model of $ec_{HS}(HSAF)$ is 
	\begin{align*}
			\|\beta\|
			&=N^\star(\|r_{S_1}^{\beta}\wedge\bigwedge_{b_{1j}\in S_1}b_{1j}\|)\ast\dots\ast N^\star(\|r_{S_k}^{\beta}\wedge\bigwedge_{b_{kj}\in S_1}b_{kj}\|)
		\\&=N^\star(\|b_{11}\|\ast\dots\ast\|b_{1j_1}\|\ast\|r_{S_1}^\beta\|)\ast\dots\ast N^\star(\|b_{k1}\|\ast\dots\ast\|b_{kj_k}\|\ast\|r_{S_k}^\beta\|)\
		\\&=h'_\beta(h'_1(\|b_{11}\|,\dots,\|b_{1j_1}\|, \|r_{S_1}^\beta\|), \dots,  h_k'(\|b_{k1}\|, \dots,\|b_{kj_k}\|, \|r_{S_k}^\beta\|))
	\end{align*}
	where $h'_\beta$ is a function $[0, 1]^k \rightarrow [0, 1]$, s.t. $h'_\beta(x'_1, \dots, x'_i, \dots, x'_k)=x'_1\ast \dots\ast x'_i\ast \dots\ast x'_k$ and $x'_i=h_i'(\|b_{i1}\|, \|b_{i2}\|,\dots,\|b_{ij_i}\|, \|r_{S_i}^\beta\|)=N^\star(\|b_{i1}\|\ast \|b_{i2}\|\ast\dots\ast\|b_{ij_i}\|\ast\|r_{S_i}^\beta\|)$. Obviously, each $h_i' (i\in\{1,2, \dots, k\})$ is a function $[0, 1]^{j_i+1} \rightarrow [0, 1]$. Next let us check that $(h'_\beta, h'_1, h'_2, \dots, h'_i, \dots, h'_k)$ is a real equation function ($k+1$)-tuple:\\
	 (a) Because $N^\star$ and $\ast$ are continuous, thus $h_i'$ is continuous and then $h'_\beta$ is continuous.
	 We also have:
	  \begin{itemize}
	 	\item $h'_\beta(1, ... , 1) =1\ast ... \ast 1= 1$;
	 	\item $h'_\beta(x'_1, \dots , 0, \dots , x'_k) =x'_1\ast \dots \ast 0\ast \dots \ast x'_k= 0$;
	 	\item The commutativity of $\ast$ implies that
	 	\begin{equation*}
	 		h'_\beta(x'_1, \dots , x'_k ) = h'_\beta(y'_1, \dots , y'_k ),
	 	\end{equation*} where $(x'_1, \dots , x'_k )$ and $(y'_1, \dots , y'_k )$ are permutations of each other.
	 \end{itemize}
	 (b) $h'_i (i\in\{1,2, \dots, k\})$ satisfies that
	 \begin{itemize}
	 	\item $x'_i=h'_i(\|b_{i1}\|,\dots,\|b_{ij_i}\|, \|r_{S_i}^\beta\|)=h'_i(u_{i1}, \dots, u_{ij_i}, u_{ij_i+1})$, i.e., $h'_i$ is a function of $j_i+1$ variables where we use variables $u_{i1}, \dots, u_{ij_i}, u_{ij_i+1}$ to represent variables $\|b_{i1}\|,\dots,\|b_{ij_i}\|, \|r_{S_i}^\beta\|$ respectively;
	 	\item $h'_i(1, ... , 1) =N^\star(1\ast \dots\ast1)=N^\star(1)=0$;
	 	\item $h'_i(u_{i1}, \dots, 0, \dots , u_{ij_i}, u_{ij_i+1}) =N^\star(0)= 1$;
	 	\item According to the commutativity of $\ast$, we have $h'_i(u_{i1}, \dots, u_{ij_i}, u_{ij_i+1}) = h'_i(v_{i1}, \dots, v_{ij_i}, v_{ij_i+1})$, where $(u_{i1}, \dots, u_{ij_i}, u_{ij_i+1})$ and $(v_{i1}, \dots, v_{ij_i}, v_{ij_i+1})$ are permutations of each other.
	 \end{itemize}
	 (c) Owing to the Monotonicity of t-norms and the Antimonotonicity of negations, $h'_\beta$ is non-decreasing w.r.t. each variable in $\{x_1, \dots, x_k\}$, and each $h_i'$ ($i \in \{1, 2, \dots, k\}$) is non-increasing w.r.t. each variable in \{$\|b_{i1}\|$, $\dots$, $\|b_{ij_i}\|$, $\|r_{S_i}^\beta\|$\}.
	 
	Thus, each encoded equational system $Eq_{ec}^{HS}$ is a real equational system $Eq^{HS}$.
\end{proof}

Note that the converse of this theorem does not hold. Specifically, an example of a $DAF$ (a $pre\text{-}HSAF^{(0)}$ where each set attacker is a singleton) provided in \cite{tang2025encoding} illustrates why the converse fails. Next, we give the existence of solutions of a real equational system $Eq^{HS}$. In fact, we have a more general theorem below.

\begin{thm}
	\label{thm:brouwer-system}
	Let $\mathbf{x} = (x_1, x_2, \dots, x_n) \in [0,1]^n$, and consider a system of $n$ equations:
	\[
	x_i = f_i(\mathbf{x}) \cdot \text{for } i = 1, 2, \dots, n,
	\]
	where each function $f_i: [0,1]^n \to [0,1]$ is continuous. Then, there exists at least one solution $\mathbf{x}^* \in [0,1]^n$ such that:
	\[
	\mathbf{x}^* = \mathbf{F}(\mathbf{x}^*),
	\]
	where $\mathbf{F}(\mathbf{x}) = \big(f_1(\mathbf{x}), f_2(\mathbf{x}), \dots, f_n(\mathbf{x})\big)$.
\end{thm}
\begin{proof} 
	\textbf{Step 1: Formalization.}\\
	Rewrite the system as a fixed-point equation:
	\[
	\mathbf{F}(\mathbf{x}) = \mathbf{x},
	\]
	where $\mathbf{F}: [0,1]^n \to [0,1]^n$ is continuous. Brouwer's Fixed-Point Theorem guarantees a solution if:
	\begin{itemize}
		\item $[0,1]^n$ is a non-empty compact convex set,
		\item $\mathbf{F}$ is continuous,
		\item $\mathbf{F}$ maps $[0,1]^n$ into itself.
	\end{itemize}
	
	\textbf{Step 2: Verify compact convexity.}\\
	\begin{itemize}
		\item \textit{Compactness:} $[0,1]^n$ is closed and bounded in $\mathbb{R}^n$ (Heine-Borel Theorem).
		\item \textit{Convexity:} For any $\mathbf{x}, \mathbf{y} \in [0,1]^n$ and $\lambda \in [0,1]$:
		\[
		\lambda x_i + (1-\lambda)y_i \in [0,1] \cdot \forall i.
		\]
	\end{itemize}
	
	\textbf{Step 3: Continuity of $\mathbf{F}$.}\\
	Each $f_i$ is continuous by assumption, and continuity is preserved under Cartesian products.
	
	\textbf{Step 4: Self-mapping property.}\\
	By construction, $f_i(\mathbf{x}) \in [0,1]$ for all $i$, so $\mathbf{F}(\mathbf{x}) \in [0,1]^n$.
	
	\textbf{Step 5: Apply Brouwer's Theorem.}\\
	All conditions are satisfied; thus, there exists $\mathbf{x}^* \in [0,1]^n$ with $\mathbf{F}(\mathbf{x}^*) = \mathbf{x}^*$.
\end{proof} 

By this theorem and the continuity of each function in the ($k+1$)-tuple $\mathcal{H}$, the general real equational system $Eq^{HS}$ of an $HSAF$ meets the conditions stated in this theorem, hence the following corollary.
\begin{cor}
There exists a solution of the general real equational system $Eq^{HS}$ for a given $HSAF$.
\end{cor}
\begin{cor}
	There exists a solution of equational system $Eq_G^{HS}$ or $Eq_P^{HS}$ or $Eq_L^{HS}$ for a given $HSAF$.
\end{cor}
For any $HLAF$ (or $BHAF$ or $SETAF$), the given equational system in the corresponding section possesses a solution because given equational system meets the conditions stated in Theorem \ref{thm:brouwer-system}.
\subsection{Semantic relationships related to the complete semantics of $HSAF$s}
For each kind of $AFSA$, we have given theorems of model equivalence under complete semantics or typical equational semantics. In Subsection \ref{subsec6.5}, we have established the formal equational approach and explored the relationship between the general real equational semantics and the $CFNE$ equational semantics. In this subsection, we study relationships between the complete semantics and equational semantics and we use the consistent symbols in previous related definitions.
\subsubsection{The relationship between the complete semantics and the restricted $Eq_G^{HS}$}
We denote by $ Eq_3^{HS}$ the \emph{3-valued equational semantics} associated with the equational semantics $ Eq^{HS} $, defined by restricting the equational system $ Eq^{HS} $ to the truth-value set $\{0, 1, \frac{1}{2}\}$. Particularly, the 3-valued equational semantics associated with the $ Eq_G^{HS} $ is denoted by $ Eq_{G(3)}^{HS} $. Denote that an assignment $\|\cdot\|$ is a model of an $HSAF$ under 3-valued equational semantics $Eq_{G(3)}^{HS}$ by $\|\cdot\| \models_{\mathcal{HS}_3^G} HSAF$. Then, we present the model equivalence theorem.
\begin{thm}\label{ceqg}
	For an $HSAF=(\mathsf{A}^{HS(0)}, \mathsf{R}^{HS(n)})$ and an assignment $\|\cdot\|: \mathsf{A}^{HS(0)}\cup\mathsf{R}^{HS(n)}\to \{0, 1, \frac{1}{2}\}$, 
	\begin{equation*}
		\|\cdot\| \models_{\mathcal{HS}_3^C} HSAF \Longleftrightarrow \|\cdot\| \models_{\mathcal{HS}_3^G} HSAF.
	\end{equation*}
\end{thm}

\begin{proof}
	With respect to the truth-value set $ \{0, 1, \frac{1}{2}\} $, it is easy to check that
	\begin{itemize}
		\item 	$\forall S_i \exists b_{ij} \in S_i: \|b_{ij}\| = 0$ or  $\|(S_i, \beta)\| = 0$ iff $\min_{i=1}^{n_\beta}\max\{1-\|r_{S_i}^{\beta}\|, 1-\|b_{i1}\|, 1-\|b_{i2}\|,\dots,1-\|b_{ik_i}\|\}=1$.
		\item $\exists S_i \forall b_{ij} \in S_i: \|b_{ij}\|= 1$ and $\|(S_i, \beta)\| =1$ iff 
		$\min_{i=1}^{n_\beta}\max\{1-\|r_{S_i}^{\beta}\|, 1-\|b_{i1}\|, 1-\|b_{i2}\|,\dots,1-\|b_{ik_i}\|\}=0$.
	\end{itemize}
We need to discuss three cases for any $\beta\in \mathbf{T}(\mathsf{R}^{HS(n)})$.
\begin{itemize}
	\item Case 1, $\|\beta\|=1$.\\
	$\|\beta\|=1$ satisfies complete semantics
	\\$\Longleftrightarrow$ $\forall S_i \exists b_{ij} \in S_i: \|b_{ij}\| = 0$ or  $\|(S_i, \beta)\| = 0$
	\\$\Longleftrightarrow$ $\min_{i=1}^{n_\beta}\max\{1-\|r_{S_i}^{\beta}\|, 1-\|b_{i1}\|, 1-\|b_{i2}\|,\dots,1-\|b_{ik_i}\|\}=1$ 
	\\$\Longleftrightarrow$ $\|\beta\|=\min_{i=1}^{n_\beta}\max\{1-\|r_{S_i}^{\beta}\|, 1-\|b_{i1}\|, 1-\|b_{i2}\|,\dots,1-\|b_{ik_i}\|\}$
\\$\Longleftrightarrow$ $\|\beta\|=1$ satisfies  semantics $Eq_{G(3)}^{HS}$.

\item Case 2, $\|\beta\|=0$.\\
$\|\beta\|=0$ satisfies complete semantics
\\$\Longleftrightarrow$ $\exists S_i \forall b_{ij} \in S_i: \|b_{ij}\|= 1$ and $\|(S_i, \beta)\| =1$ 
\\$\Longleftrightarrow$ $\min_{i=1}^{n_\beta}\max\{1-\|r_{S_i}^{\beta}\|, 1-\|b_{i1}\|, 1-\|b_{i2}\|,\dots,1-\|b_{ik_i}\|\}=0$ 
\\$\Longleftrightarrow$ $\|\beta\|=\min_{i=1}^{n_\beta}\max\{1-\|r_{S_i}^{\beta}\|, 1-\|b_{i1}\|, 1-\|b_{i2}\|,\dots,1-\|b_{ik_i}\|\}$
\\$\Longleftrightarrow$ $\|\beta\|=0$ satisfies  semantics $Eq_{G(3)}^{HS}$.

\item Case 3, $\|\beta\|=\frac{1}{2}$.\\
$\|\beta\|=\frac{1}{2}$ satisfies complete semantics
\\$\Longleftrightarrow$ it is not the situation that $[\forall S_i \exists b_{ij} \in S_i: \|b_{ij}\| = 0$ or  $\|(S_i, \beta)\| = 0]$ or $[\exists S_i \forall b_{ij} \in S_i: \|b_{ij}\|= 1$ and $\|(S_i, \beta)\| =1]$
\\$\Longleftrightarrow$ $\min_{i=1}^{n_\beta}\max\{1-\|r_{S_i}^{\beta}\|, 1-\|b_{i1}\|, 1-\|b_{i2}\|,\dots,1-\|b_{ik_i}\|\}\neq1$ and $\min_{i=1}^{n_\beta}\max\{1-\|r_{S_i}^{\beta}\|, 1-\|b_{i1}\|, 1-\|b_{i2}\|,\dots,1-\|b_{ik_i}\|\}\neq0$ from Case 1 and Case 2
\\$\Longleftrightarrow$ $\min_{i=1}^{n_\beta}\max\{1-\|r_{S_i}^{\beta}\|, 1-\|b_{i1}\|, 1-\|b_{i2}\|,\dots,1-\|b_{ik_i}\|\}=\frac{1}{2}$
\\$\Longleftrightarrow$ $\|\beta\|=\min_{i=1}^{n_\beta}\max\{1-\|r_{S_i}^{\beta}\|, 1-\|b_{i1}\|, 1-\|b_{i2}\|,\dots,1-\|b_{ik_i}\|\}$
\\$\Longleftrightarrow$ $\|\beta\|=\frac{1}{2}$ satisfies  semantics $Eq_{G(3)}^{HS}$.
\end{itemize}
From the three cases above, we have proven this theorem.
\end{proof}

This theorem shows that the equational semantics $Eq_G^{HS}$ is the generalization of the complete semantics, i.e., we have the corollary below.

\begin{cor}\label{cor28}
	For an $HSAF=(\mathsf{A}^{HS(0)}, \mathsf{R}^{HS(n)})$ and an assignment $\|\cdot\|: \mathsf{A}^{HS(0)}\cup\mathsf{R}^{HS(n)}\to \{0, 1, \frac{1}{2}\}$, 
		\begin{equation*}
			\|\cdot\| \models_{\mathcal{HS}_3^C} HSAF \Longrightarrow \|\cdot\| \models_{\mathcal{HS}_{[0,1]}^G} HSAF.
		\end{equation*}
\end{cor}
\subsubsection{The relationship between the complete semantics and the encoded semantics}

We introduce the ternarization function of numerical labellings for a given $HSAF=(\mathsf{A}^{HS(0)}, \mathsf{R}^{HS(n)})$. Let $LAB$ be the set of all numerical labellings of the $HSAF$, i.e., $LAB=\{\|\cdot\|\mid\|\cdot\|: \mathsf{A}^{HS(0)}\cup \mathsf{R}^{HS(n)}\to [0,1]\}$. Let $LAB_3$ be the set of all 3-valued labellings of the $HSAF$, i.e., $LAB_3=\{\|\cdot\|\mid\|\cdot\|: \mathsf{A}^{HS(0)}\cup \mathsf{R}^{HS(n)}\to \{0,1,\frac{1}{2}\}\}$. 
\begin{defn}\label{hster}
The \emph{ternarization} of $LAB$ associated with an $HSAF$ is a total function $T_3^{HS}: LAB\rightarrow LAB_3$, $\|\cdot\|\mapsto \|\cdot\|_3$ (i.e. $T_3^{HS}(\|\cdot\|)=\|\cdot\|_3$), such that $\forall\beta\in\mathbf{T}(\mathsf{R}^{HS(n)})$
\begin{equation*}
	T_3^{HS}(\|\cdot\|)(\beta)=\|\beta\|_3=\begin{cases}
		0 & \quad \|\beta\|=0
		\\	1 & \quad \|\beta\|=1
		\\ \frac{1}{2}  & \quad \|\beta\|\in (0,1)
	\end{cases}.
\end{equation*}
\end{defn}
$T_3^{HS}(\|\cdot\|)$ (or $\|\cdot\|_3$) is called the ternarized labelling from the labelling $\|\cdot\|$. 

A t-norm $ T $ is called a \textit{zero-divisor-free} t-norm if it satisfies:
\begin{equation*}
	\forall a, b \in (0,1], \quad T(a, b) > 0.
\end{equation*}
Equivalently, $ T $ has no zero divisors, if
$\nexists \, a, b \in (0,1]$ such that $T(a, b) = 0$.

Assume that the $\mathcal{PL}_{[0,1]}$ in Definition \ref{defn13} is equipped with a zero-divisor-free t-norm $\circledcirc$, an R-implication $I_\circledcirc$ and a negation $N$. We denote it by $\mathcal{PL}_{[0,1]}^\circledcirc$. Then we have a theorem of the model relationship.
\begin{thm}\label{circled}
	For an $HSAF=(\mathsf{A}^{HS(0)}, \mathsf{R}^{HS(n)})$ and an assignment $\|\cdot\|: \mathsf{A}^{HS(0)}\cup\mathsf{R}^{HS(n)}\to [0,1]$, 
	\begin{equation*}
		\|\cdot\| \models_{\mathcal{PL}_{[0,1]}^\circledcirc} ec_{HS}(HSAF) \Longrightarrow T_3^{HS}(\|\cdot\|) \models_{\mathcal{HS}_3^C} HSAF.
	\end{equation*}
\end{thm} 

\begin{proof}
	For the given $HSAF=(\mathsf{A}^{HS(0)}, \mathsf{R}^{HS(n)})$:\\
	an assignment $\|\cdot\|$ is a model of $ec_{HS}(HSAF)$ in $\mathcal{PL}_{[0,1]}^\circledcirc$
	\\$\Longleftrightarrow$
	\begin{equation*}
		\|ec_{HS}(HSAF)\|=\|\bigwedge_{\beta\in \mathbf{T}(\mathsf{R}^{HS(n)})}(\beta\leftrightarrow\bigwedge_{(S_i, \beta)\in \mathsf{R}^{HS(n)}}\neg(r_{S_i}^{\beta}\wedge\bigwedge_{b_{ij}\in S_i}b_{ij}))\|=1,
	\end{equation*}
	\\$\Longleftrightarrow$ $\forall\beta\in \mathbf{T}(\mathsf{R}^{HS(n)})$, 
	\begin{equation*}
		\|\beta\leftrightarrow\bigwedge_{(S_i, \beta)\in \mathsf{R}^{HS(n)}}\neg(r_{S_i}^{\beta}\wedge\bigwedge_{b_{ij}\in S_i}b_{ij})\|=1,
	\end{equation*}
	\\$\Longleftrightarrow$
	by Lemma \ref{remark2}, $\forall\beta\in \mathbf{T}(\mathsf{R}^{HS(n)})$, 
	\begin{equation*}
	\|\beta\|=\|\bigwedge_{(S_i, \beta)\in \mathsf{R}^{HS(n)}}\neg(r_{S_i}^{\beta}\wedge\bigwedge_{b_{ij}\in S_i}b_{ij})\|,
	\end{equation*}
	\\$\Longleftrightarrow$ $\forall\beta\in \mathbf{T}(\mathsf{R}^{HS(n)})$,
	\begin{align*}
		\|\beta\|
		=&N(\|r_{S_1}^{\beta}\wedge\bigwedge_{b_{1j}\in S_1}b_{1j}\|)\circledcirc\dots\circledcirc N(\|r_{S_k}^{\beta}\wedge\bigwedge_{b_{kj}\in S_1}b_{kj}\|)
		\\=&N(\|b_{11}\|\circledcirc\dots\circledcirc\|b_{1j_1}\|\circledcirc\|r_{S_1}^\beta\|)\circledcirc\dots\circledcirc N(\|b_{k1}\|\circledcirc\dots\circledcirc\|b_{kj_k}\|\circledcirc\|r_{S_k}^\beta\|).
	\end{align*}
	
	 Denoting $T_3^{HS}(\|\cdot\|)$ by $\|\cdot\|_3$, we have  
	\begin{equation*}
		\|\beta\|_3=
		\begin{cases}
			0 & \quad \|\beta\|=0
			\\	1 & \quad \|\beta\|=1
			\\ \frac{1}{2}  & \quad \|\beta\|\in (0,1)
			\end{cases}.
	\end{equation*}
	Then we need to discuss three cases.
	\begin{itemize}
		\item Case 1, $\|\beta\|_3=0$.\\
		$\|\beta\|_3=\|\beta\|=0$
		\\$\Longleftrightarrow$ $\exists S_i: N(\|b_{i1}\|\circledcirc \|b_{i2}\|\circledcirc\dots\circledcirc\|b_{ij_i}\|\circledcirc\|r_{S_i}^\beta\|)=0$ 
		\\$\Longleftrightarrow$ $\exists S_i: \|b_{i1}\|\circledcirc \|b_{i2}\|\circledcirc\dots\circledcirc\|b_{ij_i}\|\circledcirc\|r_{S_i}^\beta\|=1$ 
		\\$\Longleftrightarrow$ $\exists S_i \forall b_{ij} \in S_i: \|b_{ij}\|= 1$ and $\|(S_i, \beta)\| =1$ 
		\\$\Longleftrightarrow$ $\exists S_i \forall b_{ij} \in S_i: \|b_{ij}\|_3= 1$ and $\|(S_i, \beta)\|_3 =1$ 
		\\$\Longleftrightarrow$ $\|\beta\|_3=0$ satisfies complete semantics.
		\item Case 2, $\|\beta\|_3=1$.\\
		$\|\beta\|_3=\|\beta\|=1$
		\\$\Longleftrightarrow$ $\forall S_i: N(\|b_{i1}\|\circledcirc \|b_{i2}\|\circledcirc\dots\circledcirc\|b_{ij_i}\|\circledcirc\|r_{S_i}^\beta\|)=1$ 
		\\$\Longleftrightarrow$ $\forall S_i: \|b_{i1}\|\circledcirc \|b_{i2}\|\circledcirc\dots\circledcirc\|b_{ij_i}\|\circledcirc\|r_{S_i}^\beta\|=0$ 
		\\$\Longleftrightarrow$ $\forall S_i \exists b_{ij} \in S_i: \|b_{ij}\|= 0$ or $\|(S_i, \beta)\| =0$ 		
		\\$\Longleftrightarrow$ $\forall S_i \exists b_{ij} \in S_i: \|b_{ij}\|_3 = 0$ or  $\|(S_i, \beta)\|_3 = 0$
		\\$\Longleftrightarrow$ $\|\beta\|_3=1$ satisfies complete semantics.
		\item Case 3, $\|\beta\|_3=\frac{1}{2}$.\\
		$\|\beta\|_3=\frac{1}{2}$
		\\$\Longleftrightarrow$ not the case $\|\beta\|_3=1$ or $\|\beta\|_3=0$
		\\$\Longleftrightarrow$ $\|\beta\|_3=1$ does not satisfy complete semantics and $\|\beta\|_3=0$ does not satisfy complete semantics 
		\\$\Longleftrightarrow$ $\|\beta\|_3=\frac{1}{2}$ satisfies complete semantics.
	\end{itemize}
	From the three cases, we have proven this theorem.
\end{proof}
This theorem demonstrates that the ternarization function maps every model of an $HSAF$ under the fuzzy normal encoded semantics associated with $\mathcal{PL}_{[0,1]}^\circledcirc$ to a model under the complete semantics. Since the Product t-norm is a zero-divisor-free t-norm, we have the following corollary which has been illustrated in Example \ref{ex:hs} and its continued example for $Eq_P^{HS}$.

 \begin{cor}\label{cor.p}
 		For an $HSAF=(\mathsf{A}^{HS(0)}, \mathsf{R}^{HS(n)})$ and an assignment $\|\cdot\|: \mathsf{A}^{HS(0)}\cup\mathsf{R}^{HS(n)}\to [0,1]$, 
 	\begin{equation*}
 		\|\cdot\| \models_{\mathcal{PL}_{[0,1]}^P} ec_{HS}(HSAF) \Longrightarrow T_3^{HS}(\|\cdot\|) \models_{\mathcal{HS}_3^C} HSAF,
 	\end{equation*}
 	and then, 
 	\begin{equation*}
 		\|\cdot\| \models_{\mathcal{HS}_{[0,1]}^P} HSAF \Longrightarrow T_3^{HS}(\|\cdot\|) \models_{\mathcal{HS}_3^C} HSAF.
 	\end{equation*}
\end{cor} 
It also follows from Example \ref{ex:hs} and the continued example for $Eq_P^{HS}$ that a complete labelling of an $HSAF$ may not be obtained from the ternarization of a model of $ec_{HS}(HSAF)$ in $\mathcal{PL}_{[0,1]}^P$.

Next, we explore the case where a model of an $HSAF$ under complete semantics is also a model under some fuzzy normal encoded semantics.
Let $\odot$ be a t-norm satisfying that $\frac{1}{2}\odot \frac{1}{2}=\frac{1}{2}$ and such a t-norm is called a \emph{$\frac{1}{2}$-idempotent t-norm}.
Assume that the $\mathcal{PL}_{[0,1]}$ in Definition \ref{defn13} is equipped with a $\frac{1}{2}$-idempotent t-norm $\odot$, an R-implication $I_\odot$ and a standard negation $N$, and this logic system is denoted by $\mathcal{PL}_{[0,1]}^\odot$. Then we have a theorem below.
\begin{thm}\label{idem}	
		For an $HSAF=(\mathsf{A}^{HS(0)}, \mathsf{R}^{HS(n)})$ and an assignment $\|\cdot\|: \mathsf{A}^{HS(0)}\cup\mathsf{R}^{HS(n)}\to \{0,1,\frac{1}{2}\}$, 
	\begin{equation*}
		\|\cdot\| \models_{\mathcal{HS}_3^C} HSAF \Longrightarrow \|\cdot\| \models_{\mathcal{PL}_{[0,1]}^\odot} ec_{HS}(HSAF).
	\end{equation*}
\end{thm} 

\begin{proof}
	Similar to the proof of Theorem \ref{circled}, for the $HSAF=(\mathsf{A}^{HS(0)}, \mathsf{R}^{HS(n)})$:\\
	an assignment $\|\cdot\|$ is a model of $ec_{HS}(HSAF)$ in $\mathcal{PL}_{[0,1]}^\odot$
	\\$\Longleftrightarrow$ $\forall\beta\in \mathbf{T}(\mathsf{R}^{HS(n)})$,
	\begin{equation}\label{equ-mod}
		\|\beta\|=N(\|b_{11}\|\odot\dots\odot\|b_{1j_1}\|\odot\|r_{S_1}^\beta\|)\odot\dots\odot N(\|b_{k1}\|\odot\dots\odot\|b_{kj_k}\|\odot\|r_{S_k}^\beta\|).
	\end{equation}
	
	For the given $HSAF$ and $\forall \beta\in \mathbf{T}(\mathsf{R}^{HS(n)})$, if an assignment $\|\cdot\|$ is a model of the $HSAF$ under complete semantics, then we have
		\begin{equation*}
		\|\beta\|= 
		\begin{cases}
			1 & \quad \text{if } \forall S_i \exists b_{ij} \in S_i: \|b_{ij}\| = 0 \text{ or } \|(S_i, \beta)\| = 0
			\\ 0 & \quad \text{if } \exists S_i \forall b_{ij} \in S_i: \|b_{ij}\|= 1 \text{ and } \|(S_i, \beta)\| =1
			\\\frac{1}{2} & \quad \text{otherwise}
		\end{cases}.
	\end{equation*}

	Then we need to discuss three cases.
		\begin{itemize}
		\item Case 1, $\|\beta\|=0$ under complete semantics.\\
		$\|\beta\|=0$ satisfies complete semantics
		\\$\Longleftrightarrow$ $\exists S_i \forall b_{ij} \in S_i: \|b_{ij}\|= 1$ and $\|(S_i, \beta)\|=1$ 
		\\$\Longleftrightarrow$ $\exists S_i: \|b_{i1}\|\odot \|b_{i2}\|\odot\dots\odot\|b_{ij_i}\|\odot\|r_{S_i}^\beta\|=1$ 
		\\$\Longleftrightarrow$ $\exists S_i: N(\|b_{i1}\|\odot \|b_{i2}\|\odot\dots\odot\|b_{ij_i}\|\odot\|r_{S_i}^\beta\|)=0$ 
		\\$\Longleftrightarrow$ Equation \ref{equ-mod} holds for $\beta$.
		
		\item Case 2, $\|\beta\|=1$ under complete semantics.\\
		$\|\beta\|=1$ satisfies complete semantics
		\\$\Longleftrightarrow$ $\forall S_i \exists b_{ij} \in S_i: \|b_{ij}\|= 0$ or  $\|(S_i, \beta)\|= 0$
		\\$\Longleftrightarrow$ $\forall S_i: \|b_{i1}\|\odot \|b_{i2}\|\odot\dots\odot\|b_{ij_i}\|\odot\|r_{S_i}^\beta\|=0$ 
		\\$\Longleftrightarrow$ $\forall S_i: N(\|b_{i1}\|\odot \|b_{i2}\|\odot\dots\odot\|b_{ij_i}\|\odot\|r_{S_i}^\beta\|)=1$ 
		\\$\Longleftrightarrow$ Equation \ref{equ-mod} holds for $\beta$.
	
		\item Case 3, $\|\beta\|=\frac{1}{2}$ under complete semantics.\\
		Notice that $\|b_{i1}\|\odot \|b_{i2}\|\odot\dots\odot\|b_{ij_i}\|\odot\|r_{S_i}^\beta\|\in\{0,1,\frac{1}{2}\}$. Then:\\
		$\|\beta\|=\frac{1}{2}$ satisfies complete semantics
		\\$\Longleftrightarrow$ not Case 1 and not Case 2
		\\$\Longleftrightarrow$ by Case 1 and Case 2, $\neg\exists S_i: N(\|b_{i1}\|\odot \|b_{i2}\|\odot\dots\odot\|b_{ij_i}\|\odot\|r_{S_i}^\beta\|)=0$ 
		and $\exists S_i: N(\|b_{i1}\|\odot \|b_{i2}\|\odot\dots\odot\|b_{ij_i}\|\odot\|r_{S_i}^\beta\|)\neq1$ 
		\\$\Longleftrightarrow$ $\neg\exists S_i: N(\|b_{i1}\|\odot \|b_{i2}\|\odot\dots\odot\|b_{ij_i}\|\odot\|r_{S_i}^\beta\|)=0$ 
		and $\exists S_i: N(\|b_{i1}\|\odot \|b_{i2}\|\odot\dots\odot\|b_{ij_i}\|\odot\|r_{S_i}^\beta\|)=\frac{1}{2}$ 
		\\$\Longleftrightarrow$ Equation \ref{equ-mod} holds for $\beta$.
	\end{itemize}
	From the three cases, Equation \ref{equ-mod} holds for $\forall \beta\in \mathbf{T}(\mathsf{R}^{HS(n)})$, i.e., $\|\cdot\|$ is a model of  $ec_{HS}(HSAF)$ in $\mathcal{PL}_{[0,1]}^\odot$.
\end{proof}
Since $\mathcal{PL}_{[0,1]}^G$ is a particular instance of $\mathcal{PL}_{[0,1]}^\odot$, Corollary \ref{cor28} can also be derived via this routine.

Assume that $\mathcal{PL}_{[0,1]}^\circleddash$ is a $\mathcal{PL}_{[0,1]}$ equipped with a $\frac{1}{2}$-idempotent and zero-divisor-free t-norm $\odot$, an R-implication $I_\odot$ and a standard negation $N$. The following corollary then follows from Theorems \ref{circled} and \ref{idem}.
\begin{cor}
	For a given $HSAF$, 
	\begin{equation*}
		\{\|\cdot\|\mid\|\cdot\| \models_{\mathcal{HS}_3^C} HSAF\}=\{T_3^{HS}(\|\cdot\|)\mid \|\cdot\| \models_{\mathcal{PL}_{[0,1]}^\circleddash} ec_{HS}(HSAF)\}.		
	\end{equation*}
\end{cor}
Since the $\mathcal{PL}_{[0,1]}^G$ is a $\mathcal{PL}_{[0,1]}^\circleddash$, the following corollary can be obtained.
\begin{cor}\label{cor33}
	For a given $HSAF$, 
	\begin{equation*}
		\{\|\cdot\|\mid\|\cdot\| \models_{\mathcal{HS}_3^C} HSAF\}=\{T_3^{HS}(\|\cdot\|)\mid \|\cdot\| \models_{\mathcal{PL}_{[0,1]}^G} ec_{HS}(HSAF)\}.
	\end{equation*}
\end{cor}
\begin{cor}\label{cor34}
	For a given $HSAF$, 
	\begin{equation*}
		\{\|\cdot\|\mid\|\cdot\| \models_{\mathcal{HS}_3^C} HSAF\}=\{T_3^{HS}(\|\cdot\|)\mid \|\cdot\| \models_{\mathcal{HS}_{[0,1]}^G} HSAF\}.		
	\end{equation*}
\end{cor}
Example \ref{ex:hs} and its continued example for $Eq_G^{HS}$ illustrate Corollary \ref{cor33}. Corollary \ref{cor33} provides a new method for obtaining models of an $HSAF$ under the 3-valued complete semantics. Specifically, we first solve for the models of the $HSAF$ under the equational semantics $Eq_G^{HS}$ by equation-solving software, and then apply the function $T_3^{HS}$ to these models. From Theorem \ref{ceqg} and Corollary \ref{cor34}, we derive Corollary \ref{cor35}, which states that for a given $HSAF$ and its associated equations in $Eq_G^{HS}$, the solution set obtained by first solving these equations over $[0,1]$ and then ternarizing the results is equivalent to the solution set obtained by directly solving the equations over $\{0,1,\frac{1}{2}\}$.
 \begin{cor}\label{cor35}
 	For a given $HSAF$, 
 	\begin{equation*}
 		\{\|\cdot\|\mid\|\cdot\| \models_{\mathcal{HS}_3^G} HSAF\}=\{T_3^{HS}(\|\cdot\|)\mid \|\cdot\| \models_{\mathcal{HS}_{[0,1]}^G} HSAF\}.		
 	\end{equation*}
 \end{cor}

From Corollary \ref{cor.p}, Corollary \ref{cor34}, Theorem \ref{thm14} and Theorem \ref{thm15}, we can compare models of $Eq_G^{HS}$ and $Eq_G^{HS}$ via the function $T_3^{HS}$ as follows. 
\begin{cor}	
 For a given $HSAF$, 
 \begin{equation*}
 	\{T_3^{HS}(\|\cdot\|)\mid\|\cdot\| \models_{\mathcal{HS}_{[0,1]}^P} HSAF\}\subseteq\{T_3^{HS}(\|\cdot\|)\mid\|\cdot\| \models_{\mathcal{HS}_{[0,1]}^G} HSAF\}.
 \end{equation*}
\end{cor}

\subsection{Semantic correspondence between $AFSA$s}\label{sub6.4}
Previous work has completed multi-dimensional foundational construction around $AFSA$s: Syntactically, it has been clarified that $HSAF$ serves as the most general formal framework, integrating the higher-order attack feature of $BHAF$s and the set attack feature of $SETAFs$—a 0-level $HSAF$ is equivalent to a $SETAF$; if we do not distinguish between an argument and its singleton set, an $HSAF$ where all set attackers are singletons can be reduced to a $BHAF$, and further restricting these singletons to contain only arguments allows it to be reduced to an $HLAF$. Semantically and in terms of logical encoding, complete semantics and numerical equational semantics have been defined for $HLAF$s, $BHAF$s, $SETAF$s, and $HSAFs$ respectively, and the model equivalence between these semantics and the encoded formulas in corresponding propositional logic systems (such as $\mathcal{PL}_3^L$, the fuzzy propositional logic $\mathcal{PL}_{[0,1]}$ series) has been proven, laying a foundation for unification at the semantic level.

Building on this syntactic unification and semantic-logical connections, this subsection further focuses on the ``semantic correspondence between $AFSA$s'', aiming to achieve semantic unification of various $AFSA$s through rigorous transformation mechanisms: It first presents $SETAF$-based transformations (converting $HLAF$s, $BHAF$s, and $HSAFs$ into semantically equivalent $SETAFs$), then introduces $HSAF$-based transformations (converting $HLAF$s, $BHAF$s, and $SETAFs$ to semantically equivalent $HSAFs$), and ultimately verifies the semantic consistency within the $AFSA$s family from these two distinct transformation perspectives.
\subsubsection{The semantic unification by $SETAF$-based transformations}
For an $SETAF$, a set attacker of an argument $a$ is a set of arguments.
In \cite{gabbay2009semantics}, Gabbay changes the action that argument $a$ attacks argument $b$ in an $HLAF$ to the action that argument $a$ and attack $(a,b)$ jointly attack $b$, i.e., an $HLAF$ is transformed to an $SETAF$. In this part, we present a unified $SETAF$-based transformation framework for $AFSA$s. Each transformation function converts its respective framework to a corresponding $SETAF$ by mapping both arguments and attacks (or set attacks) as arguments in the $SETAF$. Crucially, equivalence theorems establish that models under complete semantics and three specific equational semantics are perfectly preserved between each original framework and its $SETAF$ counterpart. The proofs demonstrate this semantic equivalence by showing identical encoded formulas between an $AFSA$ and its transformed $SETAF$, ensuring consistent reasoning across different representation levels while maintaining the same semantic interpretations.
\begin{defn}
	The $SETAF$-based transformation for $HLAF$s is a function $HLtoSET: \mathcal{HLAF}\to\mathcal{SETAF}$, s.t. for an $HLAF=(\mathsf{A}^{HL(0)}, \mathsf{R}^{HL(n)})$, $HLtoSET(HLAF)=SETAF^{HL}=(\mathsf A^{HL\text{-}S}, \mathsf R^{HL\text{-}S})$, where let $A^{HL\text{-}S}=A^{HL(0)}\cup\mathsf R^{HL(n)}$, $\mathsf{A}^{HL\text{-}S}=A^{HL\text{-}S}\cup \{\bot\}=\mathsf{A}^{HL(0)}\cup \mathsf{R}^{HL(n)}$, and $\mathsf R^{HL\text{-}S}=\{(\{a_i,  r_{a_i}^{\beta}\}, \beta)\mid r_{a_i}^{\beta}\in \mathsf R^{HL(n)}\}\cup\{(\bot, r_{\bot}^{\gamma})\mid r_{\bot}^{\gamma}\in \mathsf R^{HL(n)}\}$.
\end{defn}
The function $HLtoSET$ converts an $HLAF$ to an $SETAF^{HL}$ by turning all arguments and attacks into arguments in the $SETAF^{HL}$. Furthermore, it represents both the original argument attacker and its associated attack of an element in the $HLAF$ as a binary set attacker of the corresponding element in the transformed $SETAF^{HL}$. Note that $r_{a_i}^{\beta}\in \mathsf R^{HL(n)}$ implies $\beta\in A^{HL(0)}\cup R^{HL(n)}$, and $r_{\bot}^{\gamma}\in \mathsf R^{HL(n)}$ implies $r_{\bot}^{\gamma}\in\mathsf{R}^{HL(n)}\setminus R^{HL(n)}$, i.e., $r_{\bot}^{\gamma}$ is an imaginary attack in the $HLAF$. The assignment of the $HLAF$ is on $\mathsf{A}^{HL(0)}\cup \mathsf{R}^{HL(n)}$, so we need to transform all elements in $\mathsf{A}^{HL(0)}\cup \mathsf{R}^{HL(n)}$ into arguments in the $SETAF^{HL}$. However, each argument $r_{\bot}^{\gamma}\in \mathsf R^{HL(n)}$ is not attacked in the $SETAF^{HL}$, so we add an imaginary attacker in the $SETAF^{HL}$ for each $r_{\bot}^{\gamma}$, which coincides with Definition \ref{setdefn}. Thus, both the assignment of the $HLAF$ and the assignment of the $SETAF^{HL}$ are on $\mathsf{A}^{HL(0)}\cup \mathsf{R}^{HL(n)}$ and we can calculate the value of $r_{\bot}^{\gamma}$ in the $SETAF^{HL}$. 

The corresponding equivalence theorem establishes that models under complete semantics (and three specific equational semantics) are preserved between the original $HLAF$ and its transformed $SETAF^{HL}$, ensuring semantic consistency across frameworks.
\begin{thm}
	Let a transformed $SETAF^{HL}$ of an $HLAF$ be $(\mathsf A^{HL\text{-}S}, \mathsf R^{HL\text{-}S})$. 
	For an assignment $\|\cdot\|: \mathsf A^{HL\text{-}S}\to [0,1]$, 
	\begin{equation*}
		\|\cdot\| \models_{\mathcal{HL}_3^C} HLAF \Longleftrightarrow \|\cdot\| \models_{\mathcal{S}_3^C} SETAF^{HL},
	\end{equation*}
	\begin{equation*}
		\|\cdot\| \models_{\mathcal{HL}_{[0,1]}^G} HLAF \Longleftrightarrow \|\cdot\| \models_{\mathcal{S}_{[0,1]}^G} SETAF^{HL},
	\end{equation*}
\begin{equation*}
	\|\cdot\| \models_{\mathcal{HL}_{[0,1]}^P} HLAF \Longleftrightarrow \|\cdot\| \models_{\mathcal{S}_{[0,1]}^P} SETAF^{HL},
\end{equation*}
and
	\begin{equation*}
		\|\cdot\| \models_{\mathcal{HL}_{[0,1]}^L} HLAF \Longleftrightarrow \|\cdot\| \models_{\mathcal{S}_{[0,1]}^L} SETAF^{HL}.
	\end{equation*}
\end{thm}

\begin{proof}
	For an $HLAF=(\mathsf{A}^{HL(0)}, \mathsf{R}^{HL(n)})$, by $ec_{HL}$ we have
	\begin{equation*}
		ec_{HL}(HLAF)=\bigwedge_{\beta\in A^{HL(0)}\cup R^{HL(n)}}(\beta\leftrightarrow\bigwedge_{r_{a_i}^{\beta}\in \mathsf R^{HL(n)}}\neg(a_i\wedge r_{a_i}^{\beta})).
	\end{equation*}
	For the transformed $SETAF^{HL}=(\mathsf A^{HL\text{-}S}, \mathsf R^{HL\text{-}S})$ with $A^{HL\text{-}S}=A^{HL(0)}\cup\mathsf R^{HL(n)}$, $\mathsf{A}^{HL\text{-}S}=A^{HL\text{-}S}\cup \{\bot\}=\mathsf{A}^{HL(0)}\cup \mathsf{R}^{HL(n)}$, and $\mathsf R^{HL\text{-}S}=\{(\{a_i,  r_{a_i}^{\beta}\}, \beta)\mid r_{a_i}^{\beta}\in \mathsf R^{HL(n)}\}\cup\{(\bot, r_{\bot}^{\gamma})\mid r_{\bot}^{\gamma}\in \mathsf R^{HL(n)}\}$, by $ec_{S}$ we have
	\begin{align*}
		ec_{S}(SETAF^{HL})
		=&\bigwedge_{a\in A^{HL(0)}\cup\mathsf R^{HL(n)}}(a\leftrightarrow\bigwedge_{r_{S_i}^a\in \mathsf R^{HL\text{-}S}}(\neg \bigwedge_{b_{ij}\in S_i}b_{ij}))\\
		=&\bigwedge_{\beta\in A^{HL(0)}\cup R^{HL(n)}}(\beta\leftrightarrow\bigwedge_{(\{a_i,  r_{a_i}^{\beta}\}, \beta)\in \mathsf R^{HL\text{-}S}}\neg (a_i\wedge r_{a_i}^{\beta}))\\
		&\wedge\bigwedge_{r_\bot^\gamma\in\mathsf{R}^{HL(n)}\setminus R^{HL(n)}}(r_\bot^\gamma\leftrightarrow\neg\bot).
	\end{align*}
	
		In the $\mathcal{PL}_3^L$ or any $\mathcal{PL}_{[0,1]}$, we have that\\
		$\|\cdot\|$ is a model of $ec_{HL}(HLAF)$\\
	$\Longleftrightarrow$ $\|ec_{HL}(HLAF)\|=1$\\	
	$\Longleftrightarrow$ 
	\begin{equation}\label{hlset0}
		\|\bigwedge_{\beta\in A^{HL(0)}\cup R^{HL(n)}}(\beta\leftrightarrow\bigwedge_{r_{a_i}^{\beta}\in \mathsf R^{HL(n)}}\neg(a_i\wedge r_{a_i}^{\beta}))\|=1.
	\end{equation}
	In the $\mathcal{PL}_3^L$ or any $\mathcal{PL}_{[0,1]}$, we also have that\\
	 $\|\cdot\|$ is a model of $ec_{S}(SETAF^{HL})$\\
	$\Longleftrightarrow$ $\|ec_{S}(SETAF^{HL})\|=1$\\
	$\Longleftrightarrow$ 
	\begin{equation}\label{hlset1}
		\|\bigwedge_{\beta\in A^{HL(0)}\cup R^{HL(n)}}(\beta\leftrightarrow\bigwedge_{(\{a_i,  r_{a_i}^{\beta}\}, \beta)\in \mathsf R^{HL\text{-}S}}\neg (a_i\wedge r_{a_i}^{\beta}))\|=1
	\end{equation}
	 and
	 \begin{equation}\label{hlset2}
	 \|\bigwedge_{r_\bot^\gamma\in\mathsf{R}^{HL(n)}\setminus R^{HL(n)}}(r_\bot^\gamma\leftrightarrow\neg\bot)\|=1.
	 \end{equation} 
	 From Equation \ref{hlset2}, we have that $\forall r_\bot^\gamma\in\mathsf{R}^{HL(n)}$, $\|r_\bot^\gamma\|=1$, which coincides with Definition \ref{hlafdef} or related semantic definitions for the $HLAF$.
	Equation \ref{hlset0} coincides with Equation \ref{hlset1}. Thus, in the $\mathcal{PL}_3^L$ or any $\mathcal{PL}_{[0,1]}$, we have that $\|\cdot\|$ is a model of $ec_{HL}(HLAF)$ iff $\|\cdot\|$ is a model of $ec_{S}(SETAF^{HL})$. 
	
	Therefore, by Theorem \ref{hsme} and Theorem \ref{setme} (Theorem \ref{thm2} and Theorem \ref{sgme}, Theorem \ref{Thm3} and Theorem \ref{spme}, Theorem \ref{Thm4} and Theorem \ref{slme}, respectively), an assignment is a model of the $HLAF$ under complete semantics ($Eq_G^{HL}$, $Eq_P^{HL}$, $Eq_L^{HL}$, respectively) iff it is a model of the $SETAF^{HL}$ under complete semantics ($Eq_G^{S}$, $Eq_P^{S}$, $Eq_L^{S}$, respectively).
\end{proof}

\begin{defn}
	The $SETAF$-based transformation for $BHAF$s is a function $BHtoSET: \mathcal{BHAF}\to\mathcal{SETAF}$, s.t. for a $BHAF=(\mathsf{A}^{BH(0)}, \mathsf{R}^{BH(n)})$, $BHtoSET(BHAF)=SETAF^{BH}=(\mathsf A^{BH\text{-}S}, \mathsf R^{BH\text{-}S})$, where let $A^{BH\text{-}S}=A^{BH(0)}\cup\mathsf R^{BH(n)}$, $\mathsf{A}^{BH\text{-}S}=A^{BH\text{-}S}\cup \{\bot\}=\mathsf{A}^{BH(0)}\cup \mathsf{R}^{BH(n)}$, and $\mathsf R^{BH\text{-}S}=\{(\{\alpha_i,  r_{\alpha_i}^{\beta}\}, \beta)\mid r_{\alpha_i}^{\beta}\in \mathsf R^{BH(n)}\}\cup\{(\bot, r_{\bot}^{\gamma})\mid r_{\bot}^{\gamma}\in \mathsf R^{BH(n)}\}$.
\end{defn}
The function $BHtoSET$ converts a $BHAF$ to an $SETAF^{BH}$ by mapping all arguments and attacks to arguments in the $SETAF^{BH}$. It further represents both the original attacker of an element in the $BHAF$ and its associated attack as a binary set-based attacker targeting the corresponding element in the transformed $SETAF^{BH}$. The corresponding equivalence theorem establishes that models under complete semantics (as well as three specific equational semantics) are preserved between the original $BHAF$ and its transformed $SETAF^{BH}$, thereby guaranteeing semantic consistency across the two frameworks.
\begin{thm}
Let a transformed $SETAF^{BH}$ of an $BHAF$ be $(\mathsf A^{BH\text{-}S}, \mathsf R^{BH\text{-}S})$. 
	For an assignment $\|\cdot\|: \mathsf A^{BH\text{-}S}\to [0,1]$, 
	\begin{equation*}
		\|\cdot\| \models_{\mathcal{BH}_3^C} BHAF \Longleftrightarrow \|\cdot\| \models_{\mathcal{S}_3^C} SETAF^{BH},
	\end{equation*}
	\begin{equation*}
		\|\cdot\| \models_{\mathcal{BH}_{[0,1]}^G} BHAF \Longleftrightarrow \|\cdot\| \models_{\mathcal{S}_{[0,1]}^G} SETAF^{BH},
	\end{equation*}
	\begin{equation*}
		\|\cdot\| \models_{\mathcal{BH}_{[0,1]}^P} BHAF \Longleftrightarrow \|\cdot\| \models_{\mathcal{S}_{[0,1]}^P} SETAF^{BH},
	\end{equation*}
	and
	\begin{equation*}
		\|\cdot\| \models_{\mathcal{BH}_{[0,1]}^L} BHAF \Longleftrightarrow \|\cdot\| \models_{\mathcal{S}_{[0,1]}^L} SETAF^{BH}.
	\end{equation*}
\end{thm}

\begin{proof}
	For a $BHAF=(\mathsf{A}^{BH(0)}, \mathsf{R}^{BH(n)})$, by $ec_{BH}$ we have
	\begin{equation*}
		ec_{BH}(BHAF)=\bigwedge_{\beta\in A^{BH(0)}\cup R^{BH(n)}}(\beta\leftrightarrow\bigwedge_{r_{\alpha_i}^{\beta}\in \mathsf{R}^{BH(n)}}\neg(\alpha_i\wedge r_{\alpha_i}^{\beta})).
	\end{equation*}
	For the transformed $SETAF^{BH}=(\mathsf A^{BH\text{-}S}, \mathsf R^{BH\text{-}S})$ with $A^{BH\text{-}S}=A^{BH(0)}\cup\mathsf R^{BH(n)}$, $\mathsf{A}^{BH\text{-}S}=A^{BH\text{-}S}\cup \{\bot\}=\mathsf{A}^{BH(0)}\cup \mathsf{R}^{BH(n)}$, and $\mathsf R^{BH\text{-}S}=\{(\{\alpha_i,  r_{\alpha_i}^{\beta}\}, \beta)\mid r_{\alpha_i}^{\beta}\in \mathsf R^{BH(n)}\}\cup\{(\bot, r_{\bot}^{\gamma})\mid r_{\bot}^{\gamma}\in \mathsf R^{BH(n)}\}$, by $ec_{S}$ we have
	\begin{align*}
		ec_{S}(SETAF^{BH})
		=&\bigwedge_{a\in A^{BH(0)}\cup\mathsf R^{BH(n)}}(a\leftrightarrow\bigwedge_{r_{S_i}^a\in \mathsf R^{BH\text{-}S}}(\neg \bigwedge_{b_{ij}\in S_i}b_{ij}))\\
		=&\bigwedge_{\beta\in A^{BH(0)}\cup R^{BH(n)}}(\beta\leftrightarrow\bigwedge_{(\{\alpha_i,  r_{\alpha_i}^{\beta}\}, \beta)\in \mathsf R^{BH\text{-}S}}\neg (\alpha_i\wedge r_{\alpha_i}^{\beta}))\\
		&\wedge\bigwedge_{r_\bot^\gamma\in\mathsf{R}^{BH(n)}\setminus R^{BH(n)}}(r_\bot^\gamma\leftrightarrow\neg\bot).
	\end{align*}
	
	In the $\mathcal{PL}_3^L$ or any $\mathcal{PL}_{[0,1]}$, we have that\\
	$\|\cdot\|$ is a model of $ec_{BH}(BHAF)$\\
	$\Longleftrightarrow$ $\|ec_{BH}(BHAF)\|=1$\\	
	$\Longleftrightarrow$ 
	\begin{equation}\label{bhset0}
		\|\bigwedge_{\beta\in A^{BH(0)}\cup R^{BH(n)}}(\beta\leftrightarrow\bigwedge_{r_{\alpha_i}^{\beta}\in \mathsf R^{BH(n)}}\neg(\alpha_i\wedge r_{\alpha_i}^{\beta}))\|=1.
	\end{equation}
	In the $\mathcal{PL}_3^L$ or any $\mathcal{PL}_{[0,1]}$, we also have that\\
	$\|\cdot\|$ is a model of $ec_{S}(SETAF^{BH})$\\
	$\Longleftrightarrow$ $\|ec_{S}(SETAF^{BH})\|=1$\\
	$\Longleftrightarrow$ 
	\begin{equation}\label{bhset1}
		\|\bigwedge_{\beta\in A^{BH(0)}\cup R^{BH(n)}}(\beta\leftrightarrow\bigwedge_{(\{\alpha_i,  r_{\alpha_i}^{\beta}\}, \beta)\in \mathsf R^{BH\text{-}S}}\neg (\alpha_i\wedge r_{\alpha_i}^{\beta}))\|=1
	\end{equation}
	and
	\begin{equation}\label{bhset2}
		\|\bigwedge_{r_\bot^\gamma\in\mathsf{R}^{BH(n)}\setminus R^{BH(n)}}(r_\bot^\gamma\leftrightarrow\neg\bot)\|=1.
	\end{equation} 
	From Equation \ref{bhset2}, we have that $\forall r_\bot^\gamma\in\mathsf{R}^{BH(n)}$, $\|r_\bot^\gamma\|=1$, which coincides with Definition \ref{bhdefn} or related semantic definitions for the $BHAF$.
	Equation \ref{bhset0} coincides with Equation \ref{bhset1}. Thus, in the $\mathcal{PL}_3^L$ or any $\mathcal{PL}_{[0,1]}$, we have that $\|\cdot\|$ is a model of $ec_{BH}(BHAF)$ iff $\|\cdot\|$ is a model of $ec_{S}(SETAF^{BH})$. 
	
	 Therefore, by Theorem \ref{bhcs} and Theorem \ref{setme} (Theorem \ref{bheg} and Theorem \ref{sgme}, Theorem \ref{bhep} and Theorem \ref{spme}, Theorem \ref{bhel} and Theorem \ref{slme}, respectively), an assignment is a model of the $BHAF$ under complete semantics ($Eq_G^{BH}$, $Eq_P^{BH}$, $Eq_L^{BH}$, respectively) iff it is a model of the $SETAF^{BH}$ under complete semantics ($Eq_G^{S}$, $Eq_P^{S}$, $Eq_L^{S}$, respectively).
\end{proof}
\begin{defn}
	The $SETAF$-based transformation for $HSAF$s is a function $HStoSET: \mathcal{HSAF}\to\mathcal{SETAF}$, s.t. for an $HSAF=(\mathsf{A}^{HS(0)}, \mathsf{R}^{HS(n)})$, $HStoSET(HSAF)=SETAF^{HS}=(\mathsf A^{HS\text{-}S}, \mathsf R^{HS\text{-}S})$, where let $A^{HS\text{-}S}=A^{HS(0)}\cup\mathsf R^{HS(n)}$, $\mathsf{A}^{HS\text{-}S}=A^{HS\text{-}S}\cup \{\bot\}=\mathsf{A}^{HS(0)}\cup \mathsf{R}^{HS(n)}$, and $\mathsf R^{HS\text{-}S}=\{(S_i\cup \{r_{S_i}^{\beta}\}, \beta)\mid r_{S_i}^{\beta}\in \mathsf R^{HS(n)}\}\cup\{(\{\bot\}, r_{\{\bot\}}^{\gamma})\mid r_{\{\bot\}}^{\gamma}\in \mathsf R^{HS(n)}\}$.
\end{defn}
The function $HStoSET$ converts an $HSAF$ to an $SETAF^{HS}$ by mapping all arguments and set attacks to arguments in the $SETAF^{HS}$. Furthermore, it aggregates both the set attacker and its associated set attack of an element in the $HSAF$ as a multiple set attacker of the corresponding element in the transformed $SETAF^{HS}$. The corresponding equivalence theorem establishes that models under complete semantics—as well as under three specific equational semantics—are preserved between the original $HSAF$ and its transformed counterpart, thereby ensuring semantic consistency across both frameworks.	
\begin{thm}
	Let a transformed $SETAF^{HS}$ of an $HSAF$ be $(\mathsf A^{HS\text{-}S}, \mathsf R^{HS\text{-}S})$. 
	For an assignment $\|\cdot\|: \mathsf A^{HS\text{-}S}\to [0,1]$, 
	\begin{equation*}
		\|\cdot\| \models_{\mathcal{HS}_3^C} HSAF \Longleftrightarrow \|\cdot\| \models_{\mathcal{S}_3^C} SETAF^{HS},
	\end{equation*}
	\begin{equation*}
		\|\cdot\| \models_{\mathcal{HS}_{[0,1]}^G} HSAF \Longleftrightarrow \|\cdot\| \models_{\mathcal{S}_{[0,1]}^G} SETAF^{HS},
	\end{equation*}
	\begin{equation*}
		\|\cdot\| \models_{\mathcal{HS}_{[0,1]}^P} HSAF \Longleftrightarrow \|\cdot\| \models_{\mathcal{S}_{[0,1]}^P} SETAF^{HS},
	\end{equation*}
	and
	\begin{equation*}
		\|\cdot\| \models_{\mathcal{HS}_{[0,1]}^L} HSAF \Longleftrightarrow \|\cdot\| \models_{\mathcal{S}_{[0,1]}^L} SETAF^{HS}.
	\end{equation*}
\end{thm}
\begin{proof}
	For an $HSAF=(\mathsf{A}^{HS(0)}, \mathsf{R}^{HS(n)})$, by $ec_{HS}$ we have
	\begin{equation*}
		ec_{HS}(HSAF)=\bigwedge_{\beta\in \mathbf{T}(\mathsf{R}^{HS(n)})}(\beta\leftrightarrow\bigwedge_{(S_i, \beta)\in \mathsf{R}^{HS(n)}}\neg(r_{S_i}^{\beta}\wedge\bigwedge_{b_{ij}\in S_i}b_{ij})).
	\end{equation*}
	For the transformed $SETAF^{HS}=(\mathsf A^{HS\text{-}S}, \mathsf R^{HS\text{-}S})$ with $A^{HS\text{-}S}=A^{HS(0)}\cup\mathsf R^{HS(n)}$, $\mathsf{A}^{HS\text{-}S}=A^{HS\text{-}S}\cup \{\bot\}=\mathsf{A}^{HS(0)}\cup \mathsf{R}^{HS(n)}$, and $\mathsf R^{HS\text{-}S}=\{(S_i\cup \{r_{S_i}^{\beta}\}, \beta)\mid r_{S_i}^{\beta}\in \mathsf R^{HS(n)}\}\cup\{(\{\bot\}, r_{\{\bot\}}^{\gamma})\mid r_{\{\bot\}}^{\gamma}\in \mathsf R^{HS(n)}\}$, by $ec_{S}$ we have
	\begin{align*}
		ec_{S}(SETAF^{HS})
		=&\bigwedge_{a\in A^{HS(0)}\cup\mathsf R^{HS(n)}}(a\leftrightarrow\bigwedge_{r_{T_i}^a\in \mathsf R^{HS\text{-}S}}(\neg \bigwedge_{b_{ij}\in T_i}b_{ij}))\\
		=&\bigwedge_{\beta\in A^{HS(0)}\cup R^{HS(n)}}(\beta\leftrightarrow\bigwedge_{(S_i, \beta)\in \mathsf{R}^{HS(n)}}\neg(r_{S_i}^{\beta}\wedge\bigwedge_{b_{ij}\in S_i}b_{ij}))\\
		&\wedge\bigwedge_{r_{\{\bot\}}^\gamma\in\mathsf{R}^{HS(n)}\setminus R^{HS(n)}}(r_{\{\bot\}}^\gamma\leftrightarrow\neg\bot).
	\end{align*}
	
	In the $\mathcal{PL}_3^L$ or any $\mathcal{PL}_{[0,1]}$, we have that\\
	$\|\cdot\|$ is a model of $ec_{HS}(HSAF)$\\
	$\Longleftrightarrow$ $\|ec_{HS}(HSAF)\|=1$\\	
	$\Longleftrightarrow$ 
	\begin{equation}\label{hsset0}
		\|\bigwedge_{\beta\in \mathbf{T}(\mathsf{R}^{HS(n)})}(\beta\leftrightarrow\bigwedge_{(S_i, \beta)\in \mathsf{R}^{HS(n)}}\neg(r_{S_i}^{\beta}\wedge\bigwedge_{b_{ij}\in S_i}b_{ij}))\|=1.
	\end{equation}
	In the $\mathcal{PL}_3^L$ or any $\mathcal{PL}_{[0,1]}$, we also have that\\
	$\|\cdot\|$ is a model of $ec_{S}(SETAF^{HS})$\\
	$\Longleftrightarrow$ $\|ec_{S}(SETAF^{HS})\|=1$\\
	$\Longleftrightarrow$ 
	\begin{equation}\label{hsset1}
		\|\bigwedge_{\beta\in A^{HS(0)}\cup R^{HS(n)}}(\beta\leftrightarrow\bigwedge_{(S_i, \beta)\in \mathsf{R}^{HS(n)}}\neg(r_{S_i}^{\beta}\wedge\bigwedge_{b_{ij}\in S_i}b_{ij}))\|=1
	\end{equation}
	and
	\begin{equation}\label{hsset2}
		\|\bigwedge_{r_{\{\bot\}}^\gamma\in\mathsf{R}^{HS(n)}\setminus R^{HS(n)}}(r_{\{\bot\}}^\gamma\leftrightarrow\neg\bot)\|=1.
	\end{equation} 
	From Equation \ref{hsset2}, we have that $\forall r_\bot^\gamma\in\mathsf{R}^{HS(n)}$, $\|r_\bot^\gamma\|=1$, which coincides with Definition \ref{hsafd} or related semantic definitions for the $HSAF$.
	Equation \ref{hsset0} coincides with Equation \ref{hsset1}. Thus, in the $\mathcal{PL}_3^L$ or any $\mathcal{PL}_{[0,1]}$, we have that $\|\cdot\|$ is a model of $ec_{HS}(HSAF)$ iff $\|\cdot\|$ is a model of $ec_{S}(SETAF^{HS})$. 
	
	 Therefore, by Theorem \ref{thm-pl3} and Theorem \ref{setme} (Theorem \ref{thm14} and Theorem \ref{sgme}, Theorem \ref{thm15} and Theorem \ref{spme}, Theorem \ref{thm16} and Theorem \ref{slme}, respectively), an assignment is a model of the $HSAF$ under complete semantics ($Eq_G^{HS}$, $Eq_P^{HS}$, $Eq_L^{HS}$, respectively) iff it is a model of the $SETAF^{HS}$ under complete semantics ($Eq_G^{S}$, $Eq_P^{S}$, $Eq_L^{S}$, respectively).
\end{proof}
These imply that any $HLAF$, $BHAF$, or $HSAF$ can be transformed into a semantically equivalent $SETAF$, i.e., one that preserves the models of the original framework under corresponding semantics. Hence, we achieve a unification by the encoding method from the aspect of semantic equivalence. This is also the reason why we call any of $HLAF$s, $BHAF$s, $SETAF$s or $HSAF$s \emph{argumentation frameworks with set attackers} ($AFSA$s).
\subsubsection{The semantic unification by $HSAF$-based transformations}
We have achieved the semantic unification of $AFSAs$ through $SETAF$-based transformations. Previously, the understanding that the semantics of $HSAFs$ extend those of $BHAFs$ and $SETAFs$ was derived mainly from intuitive analysis, without formalized transformation rules or rigorous specifications of their relationships. As a general framework integrating higher-order attacks (inherited from $BHAFs$) and set attacks (inherited from $SETAFs$), $HSAFs$ inherently have the potential to act as a core carrier for semantically unifying diverse $AFSAs$. To further refine and formalize the semantic connections among various $AFSAs$, this part focuses on developing strict transformation mechanisms to convert $HLAFs$, $BHAFs$, and $SETAFs$ into semantically equivalent $HSAFs$ respectively. This formal transformation ensures that the complete semantics and equational semantics of the original frameworks are fully retained, ultimately establishing rigorous model equivalence relationships between each type of $AFSA$ and $HSAFs$—filling the gap of relying solely on intuition in prior analyses of semantic extensions.
\begin{defn}
	The $HSAF$-based transformation for $HLAF$s is a function $HLtoHS: \mathcal{HLAF}\to\mathcal{HSAF}$, s.t. for an $HLAF=(\mathsf{A}^{HL(0)}, \mathsf{R}^{HL(n)})$, $HLtoHS(HLAF)=HSAF^{HL}=(\mathsf A^{HL\text{-}HS}, \mathsf R^{HL\text{-}HS})$, where let $\mathsf A^{HL\text{-}HS}=\mathsf A^{HL(0)}$ and $\mathsf R^{HL\text{-}HS}=\{r_{\{a_i\}}^{\beta}\mid r_{a_i}^{\beta}\in \mathsf R^{HL(n)}\}$.
\end{defn}
Without distinguishing an element and its singleton set, we have $R^{HL\text{-}HS}= \mathsf R^{HL(n)}$. In this sense, we have the theorem below.
\begin{thm}\label{hlcg}
Let a transformed $HSAF^{HL}$ of an $HLAF$ be $(\mathsf A^{HL\text{-}HS}, \mathsf R^{HL\text{-}HS})$. 
For an assignment $\|\cdot\|: \mathsf A^{HL\text{-}HS}\cup\mathsf R^{HL\text{-}HS}\to [0,1]$, 
\begin{equation*}
	\|\cdot\| \models_{\mathcal{HL}_3^C} HLAF \Longleftrightarrow \|\cdot\| \models_{\mathcal{HS}_3^C} HSAF^{HL},
\end{equation*}
\begin{equation*}
	\|\cdot\| \models_{\mathcal{HL}_{[0,1]}^G} HLAF \Longleftrightarrow \|\cdot\| \models_{\mathcal{HS}_{[0,1]}^G} HSAF^{HL},
\end{equation*}
\begin{equation*}
	\|\cdot\| \models_{\mathcal{HL}_{[0,1]}^P} HLAF \Longleftrightarrow \|\cdot\| \models_{\mathcal{HS}_{[0,1]}^P} HSAF^{HL},
\end{equation*}
and
\begin{equation*}
	\|\cdot\| \models_{\mathcal{HL}_{[0,1]}^L} HLAF \Longleftrightarrow \|\cdot\| \models_{\mathcal{HS}_{[0,1]}^L} HSAF^{HL}.
\end{equation*}
\end{thm}

\begin{proof}
	For the given $HLAF=(\mathsf{A}^{HL(0)}, \mathsf{R}^{HL(n)})$, by $ec_{HL}$ we have
	\begin{equation*}
		ec_{HL}(HLAF)=\bigwedge_{\beta\in A^{HL(0)}\cup R^{HL(n)}}(\beta\leftrightarrow\bigwedge_{r_{a_i}^{\beta}\in \mathsf R^{HL(n)}}\neg(a_i\wedge r_{a_i}^{\beta})).
	\end{equation*}
	For the transformed $HSAF^{HL}=(\mathsf A^{HL\text{-}HS}, \mathsf R^{HL\text{-}HS})$ with $\mathsf A^{HL\text{-}HS}=\mathsf A^{HL(0)}$ and $\mathsf R^{HL\text{-}HS}=\{r_{\{a_i\}}^{\beta}\mid r_{a_i}^{\beta}\in \mathsf R^{HL(n)}\}$, in the $\mathcal{PL}_3^L$ or any $\mathcal{PL}_{[0,1]}$s by $ec_{HS}$ we have
	\begin{align*}
		ec_{HS}(HSAF^{HL})&=\bigwedge_{\beta\in \mathbf{T}(\mathsf R^{HL\text{-}HS})}(\beta\leftrightarrow\bigwedge_{(S_i, \beta)\in \mathsf R^{HL\text{-}HS}}\neg(r_{S_i}^{\beta}\wedge\bigwedge_{b_{ij}\in S_i}b_{ij}))\\
		&=\bigwedge_{\beta\in A^{HL(0)}\cup R^{HL(n)}}(\beta\leftrightarrow\bigwedge_{r_{\{a_i\}}^{\beta}\in \mathsf R^{HL\text{-}HS}}\neg(r_{\{a_i\}}^{\beta}\wedge a_i)).
	\end{align*}	
	Therefore, the encoded formula of the $HLAF$ is the same as the encoded formula of the transformed $HSAF^{HL}$ without distinguishing the representations of an element and its singleton set (i.e., $\{a_i\}$ and $a_i$ in the encoded formulas above). Therefore, by Theorem \ref{hsme} and Theorem \ref{thm-pl3} (Theorem \ref{thm2} and Theorem \ref{thm14}, Theorem \ref{Thm3} and Theorem \ref{thm15}, Theorem \ref{Thm4} and Theorem \ref{thm16}, respectively), an assignment is a model of the $HLAF$ under complete semantics ($Eq_G^{HL}$, $Eq_P^{HL}$, $Eq_L^{HL}$, respectively) iff it is a model of the $HSAF^{HL}$ under complete semantics ($Eq_G^{HS}$, $Eq_P^{HS}$, $Eq_L^{HS}$, respectively).
\end{proof}

\begin{defn}
	The $HSAF$-based transformation for $BHAF$s is a function $BHtoHS: \mathcal{BHAF}\to\mathcal{HSAF}$, s.t. for a $BHAF=(\mathsf{A}^{BH(0)}, \mathsf{R}^{BH(n)})$, $BHtoHS(BHAF)=HSAF^{BH}=(\mathsf A^{BH\text{-}HS}, \mathsf R^{BH\text{-}HS})$, where let $\mathsf A^{BH\text{-}HS}=\mathsf{A}^{BH(0)}$ and $\mathsf R^{BH\text{-}HS}=\{r_{\{\alpha_i\}}^{\beta}\mid r_{\alpha_i}^{\beta}\in \mathsf{R}^{BH(n)}\}$.
\end{defn}
Without distinguishing an element and its singleton set, we have $R^{BH\text{-}HS}= \mathsf R^{BH(n)}$. In this sense, we have the theorem below.
\begin{thm}\label{bhcg}
Let a transformed $HSAF^{BH}$ of an $BHAF$ be $(\mathsf A^{BH\text{-}HS}, \mathsf R^{BH\text{-}HS})$. 
For an assignment $\|\cdot\|: \mathsf A^{BH\text{-}HS}\cup\mathsf R^{BH\text{-}HS}\to [0,1]$, 
\begin{equation*}
	\|\cdot\| \models_{\mathcal{BH}_3^C} BHAF \Longleftrightarrow \|\cdot\| \models_{\mathcal{HS}_3^C} HSAF^{BH},
\end{equation*}
\begin{equation*}
	\|\cdot\| \models_{\mathcal{BH}_{[0,1]}^G} BHAF \Longleftrightarrow \|\cdot\| \models_{\mathcal{HS}_{[0,1]}^G} HSAF^{BH},
\end{equation*}
\begin{equation*}
	\|\cdot\| \models_{\mathcal{BH}_{[0,1]}^P} BHAF \Longleftrightarrow \|\cdot\| \models_{\mathcal{HS}_{[0,1]}^P} HSAF^{BH},
\end{equation*}
and
\begin{equation*}
	\|\cdot\| \models_{\mathcal{BH}_{[0,1]}^L} BHAF \Longleftrightarrow \|\cdot\| \models_{\mathcal{HS}_{[0,1]}^L} HSAF^{BH}.
\end{equation*}
\end{thm}
\begin{proof}
	For the given $BHAF=(\mathsf{A}^{BH(0)}, \mathsf{R}^{BH(n)})$, by $ec_{BH}$ we have
	\begin{equation*}
		ec_{BH}(BHAF)=\bigwedge_{\beta\in A^{BH(0)}\cup R^{BH(n)}}(\beta\leftrightarrow\bigwedge_{r_{\alpha_i}^{\beta}\in \mathsf{R}^{BH(n)}}\neg(\alpha_i\wedge r_{\alpha_i}^{\beta})).
	\end{equation*}
	For the transformed $HSAF^{BH}=(\mathsf A^{BH\text{-}HS}, \mathsf R^{BH\text{-}HS})$ with $\mathsf A^{BH\text{-}HS}=\mathsf{A}^{BH(0)}$ and $\mathsf R^{BH\text{-}HS}=\{r_{\{\alpha_i\}}^{\beta}\mid r_{\alpha_i}^{\beta}\in \mathsf{R}^{BH(n)}\}$, by $ec_{HS}$ we have
	\begin{align*}
		ec_{HS}(HSAF^{BH})&=\bigwedge_{\beta\in \mathbf{T}(\mathsf R^{BH\text{-}HS})}(\beta\leftrightarrow\bigwedge_{(S_i, \beta)\in \mathsf R^{BH\text{-}HS}}\neg(r_{S_i}^{\beta}\wedge\bigwedge_{b_{ij}\in S_i}b_{ij}))\\
		&=\bigwedge_{\beta\in \mathbf{T}(\mathsf R^{BH\text{-}HS})}(\beta\leftrightarrow\bigwedge_{r_{\{\alpha_i\}}^{\beta}\in \mathsf R^{BH\text{-}HS}}\neg(r_{\{\alpha_i\}}^{\beta}\wedge \alpha_i)).
	\end{align*}
	Therefore, the encoded formula of the $BHAF$ is the same as the encoded formula of the transformed $HSAF^{BH}$ without distinguishing the representations of an element and its singleton set (i.e., $\{\alpha_i\}$ and $\alpha_i$ in the encoded formulas above). Therefore, by Theorem \ref{bhcs} and Theorem \ref{thm-pl3} (Theorem \ref{bheg} and Theorem \ref{thm14}, Theorem \ref{bhep} and Theorem \ref{thm15}, Theorem \ref{bhel} and Theorem \ref{thm16}, respectively), an assignment is a model of the $BHAF$ under complete semantics ($Eq_G^{BH}$, $Eq_P^{BH}$, $Eq_L^{BH}$, respectively) iff it is a model of the $HSAF^{BH}$ under complete semantics ($Eq_G^{HS}$, $Eq_P^{HS}$, $Eq_L^{HS}$, respectively).
\end{proof}

\begin{defn}\label{seths}
	The $HSAF$-based transformation for $SETAF$s is a function $StoHS: \mathcal{SETAF}\to\mathcal{HSAF}$, s.t. for a given $SETAF=(\mathsf A^{S}, \mathsf R^S)$, $StoHS(SETAF)=HSAF^{S}=(\mathsf A^{S\text{-}HS}, \mathsf R^{S\text{-}HS})$, where let $\mathsf A^{S\text{-}HS}=\mathsf A^{S}$ and $\mathsf R^{S\text{-}HS}=\mathsf R^S\cup(\{\{\bot\}\}\times R^S)$.
\end{defn}
Denote an assignment on $\mathsf A^{S\text{-}HS}\cup\mathsf R^{S\text{-}HS}$ as $\|\cdot\|_\mathsf{AR}$. Denote the restriction of $\|\cdot\|_\mathsf{AR}$ on $\mathsf A^{S}$ as $\|\cdot\|_{\mathsf A}$, i.e., $\forall a\in\mathsf A^{S}$, $\|a\|_{\mathsf A}=\|a\|_\mathsf{AR}$. Then we have the theorem below.
\begin{thm}\label{setcg}
Let a transformed $HSAF^{S}$ of an $SETAF$ be $(\mathsf A^{S\text{-}HS}, \mathsf R^{S\text{-}HS})$. 
For an assignment $\|\cdot\|_\mathsf{AR}: \mathsf A^{S\text{-}HS}\cup\mathsf R^{S\text{-}HS}\to [0,1]$, 
\begin{equation*}
	\|\cdot\|_\mathsf{A} \models_{\mathcal{S}_3^C} SETAF \Longleftrightarrow \|\cdot\|_\mathsf{AR} \models_{\mathcal{HS}_3^C} HSAF^{S},
\end{equation*}
\begin{equation*}
	\|\cdot\|_\mathsf{A} \models_{\mathcal{S}_{[0,1]}^G} SETAF \Longleftrightarrow \|\cdot\|_\mathsf{AR} \models_{\mathcal{HS}_{[0,1]}^G} HSAF^{S},
\end{equation*}
\begin{equation*}
	\|\cdot\|_\mathsf{A} \models_{\mathcal{S}_{[0,1]}^P} SETAF \Longleftrightarrow \|\cdot\|_\mathsf{AR} \models_{\mathcal{HS}_{[0,1]}^P} HSAF^{S},
\end{equation*}
and
\begin{equation*}
	\|\cdot\|_\mathsf{A} \models_{\mathcal{S}_{[0,1]}^L} SETAF \Longleftrightarrow \|\cdot\|_\mathsf{AR} \models_{\mathcal{HS}_{[0,1]}^L} HSAF^{S}.
\end{equation*}
\end{thm}

\begin{proof}
	For the given $SETAF=(\mathsf A^{S}, \mathsf R^S)$, by $ec_{S}$ we have
	\begin{equation*}
		ec_{S}(SETAF)=\bigwedge_{a\in A^S}(a\leftrightarrow\bigwedge_{r_{S_i}^a\in \mathsf R^S}(\neg \bigwedge_{b_{ij}\in S_i}b_{ij})).
	\end{equation*}
	For the transformed $HSAF^{S}=(\mathsf A^{S\text{-}HS}, \mathsf R^{S\text{-}HS})$ with $\mathsf A^{S\text{-}HS}=\mathsf A^{S}$ and $\mathsf R^{S\text{-}HS}=\mathsf R^S\cup(\{\{\bot\}\}\times R^S)$, in the $\mathcal{PL}_3^L$ or any $\mathcal{PL}_{[0,1]}$ by $ec_{HS}$ we have
	\begin{align*}
		ec_{HS}(HSAF^{S})=&\bigwedge_{\beta\in \mathbf{T}(\mathsf R^{S\text{-}HS})}(\beta\leftrightarrow\bigwedge_{(S_i, \beta)\in \mathsf R^{S\text{-}HS}}\neg(r_{S_i}^{\beta}\wedge\bigwedge_{b_{ij}\in S_i}b_{ij}))\\
		=&\bigwedge_{\beta\in A^{S}\cup R^{S}}(\beta\leftrightarrow\bigwedge_{r_{S_i}^{\beta}\in \mathsf R^{S\text{-}HS}}\neg(r_{S_i}^{\beta}\wedge\bigwedge_{b_{ij}\in S_i}b_{ij}))\\
		=&\bigwedge_{a\in A^{S}}(a\leftrightarrow\bigwedge_{r_{S_i}^{a}\in \mathsf R^{S}}\neg(r_{S_i}^{a}\wedge\bigwedge_{b_{ij}\in S_i}b_{ij}))\wedge\bigwedge_{\beta\in R^{S}}(\beta\leftrightarrow\neg(r_{\{\bot\}}^{\beta}\wedge \bot))\\
		=&\bigwedge_{a\in A^{S}}(a\leftrightarrow\bigwedge_{r_{S_i}^{a}\in \mathsf R^{S}}\neg(r_{S_i}^{a}\wedge\bigwedge_{b_{ij}\in S_i}b_{ij}))\wedge\bigwedge_{\beta\in R^{S}}(\beta\leftrightarrow\neg\bot).
	\end{align*}
	
	In the $\mathcal{PL}_3^L$ or any $\mathcal{PL}_{[0,1]}$, we have that\\
	$\|ec_{S}(SETAF)\|_\mathsf A=1$\\
	$\Longleftrightarrow$ $\|\bigwedge_{a\in A^S}(a\leftrightarrow\bigwedge_{r_{S_i}^a\in \mathsf R^S}(\neg \bigwedge_{b_{ij}\in S_i}b_{ij}))\|_\mathsf A=1$\\
	and\\
	 $\|ec_{HS}(HSAF^{S})\|_\mathsf{AS}=1$\\
	$\Longleftrightarrow$ $\|\bigwedge_{a\in A^{S}}(a\leftrightarrow\bigwedge_{r_{S_i}^{a}\in \mathsf R^{S}}\neg(r_{S_i}^{a}\wedge\bigwedge_{b_{ij}\in S_i}b_{ij}))\wedge\bigwedge_{\beta\in R^{S}}(\beta\leftrightarrow\neg\bot)\|_\mathsf{AS}=1$\\
	$\Longleftrightarrow$ $\|\bigwedge_{a\in A^{S}}(a\leftrightarrow\bigwedge_{r_{S_i}^{a}\in \mathsf R^{S}}\neg(r_{S_i}^{a}\wedge\bigwedge_{b_{ij}\in S_i}b_{ij}))\|_\mathsf{AS}=\|\bigwedge_{\beta\in R^{S}}(\beta\leftrightarrow\neg\bot)\|_\mathsf{AS}=1$\\
	$\Longleftrightarrow$ $\|\bigwedge_{a\in A^{S}}(a\leftrightarrow\bigwedge_{r_{S_i}^{a}\in \mathsf R^{S}}\neg(r_{S_i}^{a}\wedge\bigwedge_{b_{ij}\in S_i}b_{ij}))\|_\mathsf{AS}=1$ and $\forall \beta\in R^{S}:$ $\|\beta\|_\mathsf{AS}=1$.\\
	Since for all $\beta\in R^{S}$ we have $\|\beta\|_\mathsf{AS}=1$, and if $(\{\bot\}, \beta)\in \mathsf{R}^{HS(n)}$ then we have $\|(\{\bot\}, \beta)\|_\mathsf{AS}=1$ (from Definition \ref{hsafd}), we obtain that $\|r_{S_i}^{a}\|_\mathsf{AS}=1$ ($\forall r_{S_i}^{a}\in \mathsf R^{S}$).
	Thus, $\|\bigwedge_{a\in A^{S}}(a\leftrightarrow\bigwedge_{r_{S_i}^{a}\in \mathsf R^{S}}\neg(r_{S_i}^{a}\wedge\bigwedge_{b_{ij}\in S_i}b_{ij}))\|_\mathsf{AS}=\|\bigwedge_{a\in A^{S}}(a\leftrightarrow\bigwedge_{r_{S_i}^{a}\in \mathsf R^{S}}\neg(\bigwedge_{b_{ij}\in S_i}b_{ij}))\|_\mathsf{AS}=1$. Therefore, we have that\\ $\|ec_{HS}(HSAF^{S})\|_\mathsf{AS}=1$\\
	$\Longleftrightarrow$ $\|\bigwedge_{a\in A^{S}}(a\leftrightarrow\bigwedge_{r_{S_i}^{a}\in \mathsf R^{S}}\neg(\bigwedge_{b_{ij}\in S_i}b_{ij}))\|_\mathsf{AS}=1$ and $\forall \beta\in R^{S}:$ $\|\beta\|_\mathsf{AS}=1$.
	
	Thus, in the $\mathcal{PL}_3^L$ or any $\mathcal{PL}_{[0,1]}$, $\|\cdot\|_\mathsf A$ is a model of the $ec_{S}(SETAF)$ iff $\|\cdot\|_\mathsf{AR}$ (s.t. $\forall \beta\in R^{S}:$ $\|\beta\|_\mathsf{AS}=1$) is a model of the $ec_{HS}(HSAF^{S})$. Therefore, by Theorem \ref{setme} and Theorem \ref{thm-pl3} (Theorem \ref{sgme} and Theorem \ref{thm14}, Theorem \ref{spme} and Theorem \ref{thm15}, Theorem \ref{slme} and Theorem \ref{thm16}, respectively), the theorem is proven.
\end{proof}
Similar to Definition \ref{hster}, we directly denote the ternarization of labellings associated with an $HLAF$ (or a $BHAF$ or an $SETAF$, respectively) as $T_3^{HL}$ (or $T_3^{BH}$ or $T_3^{S}$, respectively). Then, by Corollary \ref{cor33} and Theorem \ref{hlcg} (or Theorem \ref{bhcg} or Theorem \ref{setcg}, respectively), we have corollaries as follows.
\begin{cor}
		For a given $HLAF$, 
	\begin{equation*}
		\{\|\cdot\|\mid\|\cdot\| \models_{\mathcal{HL}_3^C} HLAF\}=\{T_3^{HL}(\|\cdot\|)\mid \|\cdot\| \models_{\mathcal{PL}_{[0,1]}^G} ec_{HL}(HLAF)\}.		
	\end{equation*}
\end{cor}
\begin{cor}
	For a given $BHAF$, 
\begin{equation*}
	\{\|\cdot\|\mid\|\cdot\| \models_{\mathcal{BH}_3^C} BHAF\}=\{T_3^{BH}(\|\cdot\|)\mid \|\cdot\| \models_{\mathcal{PL}_{[0,1]}^G} ec_{BH}(BHAF)\}.		
\end{equation*}
\end{cor}
\begin{cor}
		For a given $SETAF$, 
	\begin{equation*}
		\{\|\cdot\|\mid\|\cdot\| \models_{\mathcal{BH}_3^C} SETAF\}=\{T_3^{S}(\|\cdot\|)\mid \|\cdot\| \models_{\mathcal{PL}_{[0,1]}^G} ec_{S}(SETAF)\}.		
	\end{equation*}
\end{cor}
Theorem \ref{hlcg} indicates that models of a given $HLAF$ under its complete semantics (corresponding to $Eq_G^{HL}$, $Eq_P^{HL}$, or $Eq_L^{HL}$) are equivalent to models of its transformed $HSAF^{HL}$ under $HSAF$'s complete semantics (corresponding to $Eq_G^{HS}$, $Eq_P^{HS}$, or $Eq_L^{HS}$). Theorem \ref{bhcg} shows that models of a given $BHAF$ under its complete semantics (corresponding to $Eq_G^{BH}$, $Eq_P^{BH}$, or $Eq_L^{BH}$) are equivalent to models of its transformed $HSAF^{BH}$ under $HSAF$'s complete semantics (corresponding to $Eq_G^{HS}$, $Eq_P^{HS}$, or $Eq_L^{HS}$). Theorem \ref{seths} reveals that models of a given $SETAF$ under its complete semantics (corresponding to $Eq_G^{S}$, $Eq_P^{S}$, or $Eq_L^{S}$) are equivalent to restricted models of its transformed $HSAF^{S}$ under $HSAF$'s complete semantics (corresponding to $Eq_G^{HS}$, $Eq_P^{HS}$, or $Eq_L^{HS}$). Thus, the complete and three equational semantics of $HLAF$s, $BHAF$s, and $SETAF$s can be equivalently expressed via the complete and three equational semantics of their respective transformed $HSAF$s. Based on these results, it is natural to define the general real equational system of $HLAF$s, $BHAF$s, and $SETAF$s—denoted by $Eq^{HL}$, $Eq^{BH}$, and $Eq^{S}$ respectively—based on the general real equational system of their respective transformed $HSAF$s. Thus, we present the corresponding definitions separately below, with notation specific to each definition.
\begin{itemize}
	\item A \emph{general real equational system} $Eq^{HL}$ for an $HLAF = (\mathsf{A}^{HL(0)}, \mathsf{R}^{HL(n)})$ over $[0, 1]$ is defined as follows. Let $\|\bot\|=0$, $\|r_{\bot}^a\|=1$ (if $r_{\bot}^a\in\mathsf{R}^{HL(n)}$), and the equation for each $\beta\in {A}^{HL(0)}\cup{R}^{HL(n)}$ be
	\begin{equation*}
		\|\beta\| = h_\beta(h_1(\|a_1\|, \|r_{a_1}^\beta\|), \dots,  h_k(\|a_k\|, \|r_{a_k}^\beta\|)),
	\end{equation*}
	where $(h_\beta, h_1, \dots, h_k)$ is a real equation function ($k+1$)-tuple and for each $i \in \{1, \dots, k\}$, $(a_i, \beta)\in\mathsf{R}^{HL(n)}$.
	
	\item A \emph{general real equational system} $Eq^{BH}$ for a $BHAF = (\mathsf{A}^{BH(0)}, \mathsf{R}^{BH(n)})$ over $[0, 1]$ is defined as follows. Let $\|\bot\|=0$, $\|r_{\bot}^\alpha\|=1$ (if $r_{\bot}^\alpha\in\mathsf{R}^{BH(n)}$), and the equation for each $\beta\in {A}^{BH(0)}\cup{R}^{BH(n)}$ be
	\begin{equation*}
		\|\beta\| = h_\beta(h_1(\|\alpha_1\|, \|r_{\alpha_1}^\beta\|), \dots,  h_k(\|\alpha_k\|, \|r_{\alpha_k}^\beta\|)),
	\end{equation*}
	where $(h_\beta, h_1, \dots, h_k)$ is a real equation function ($k+1$)-tuple and for each $i \in \{1, \dots, k\}$, $(\alpha_i, \beta)\in\mathsf{R}^{BH(n)}$.
	
	\item A \emph{general real equational system} $Eq^{S}$ for an $SETAF = (\mathsf{A}^{S}, \mathsf{R}^{S})$ over $[0, 1]$ is defined as follows. Let $\|\bot\|=0$, $\|r_{\{\bot\}}^c\|=1$ (if $r_{\{\bot\}}^c\in\mathsf{R}^{S}$), and the equation for each $a\in {A}^{S}$ be
	\begin{equation*}
		\|a\| = h_a(h_1(\|b_{11}\|,\dots,\|b_{1j_1}\|, 1), \dots,  h_k(\|b_{k1}\|,\dots,\|b_{kj_k}\|, 1)),
	\end{equation*}
	where $(h_a, h_1, \dots, h_k)$ is a real equation function ($k+1$)-tuple and for each $i \in \{1, \dots, k\}$, the set $S_i = \{b_{i1}, \dots, b_{ij_i}\}$ is a set attacker of $a$.
\end{itemize}
\section{Related work}
The logical encoding of advanced structures in argumentation frameworks has been an active area of research, extending beyond the basic graph-based model to capture more complex phenomena. In this section, we provide a discussion and comparison between our work and other related works, with a particular focus on encoding methods for collective or higher-order structures of $AF$s. We examine how different logical systems have been employed to represent these structures and highlight the unique advantages of our propositional logic-based approach in terms of semantic unity and computational efficiency.

In \cite{gabbay2015attack}, the author encodes $AF$s by extending classical logic with strong negation (CN system) and its modal extension CNN, formalizing attacks via a custom strong negation operator $N$ (interpreted through a two-world Kripke model in CNN) and constructing a theory of separate non-logical axioms to map argument labellings. While they handle $SETAF$s and $HLAF$s in CNN, their approach relies on novel logical operators and a discrete two-world structure. In contrast, our work uses {\L}ukasiewicz’s 3-valued and fuzzy propositional logic, encodes the entire framework as a single integrated formula (avoiding scattered axioms), addresses more general $HSAF$s, enables seamless extension from 3-valued to fuzzy semantics, and fits lightweight solvers—filling gaps in general set-based higher-order handling while prioritizing simplicity and computational efficiency. 

In \cite{dvovrak2023expressiveness}, the authors encode $SETAF$s by expressing them as acceptance conditions in abstract dialectical frameworks ($ADF$s). The logical framework used for the encoded $SETAF$ formulas is $\mathcal{PL}_2$. The authors also implement three-valued semantics for $SETAF$s, which is equivalent to the three-valued semantics of the corresponding $ADF$s, by utilizing $\mathcal{PL}_2$ as the underlying computational engine. This is achieved through a characteristic operator based on partial evaluations of acceptance conditions. In our paper, we encode $SETAF$s into $\mathcal{PL}_3^L$ and $\mathcal{PL}_{[0,1]}$. Therefore, the approaches in the two papers are based on different logic systems. Due to the encoding methodology associated with $\mathcal{PL}_3^L$, we do not require any transformation such as a characteristic operator. The model equivalence relationship is straightforward, and the logical essence of the complete semantics is clearly revealed through $\mathcal{PL}_3^L$.

Related work on logical encoding of enriched argumentation frameworks primarily relies on first-order logic to capture complex interactions such as recursive attacks \cite{cayrol2017logical}, higher-order relations \cite{claudette2018logical}, evidential (or necessary) supports \cite{cayrol2020logical, lagasquie2021evidential, lagasquie2021necessary, lagasquie2023handling}, and generic cases \cite{besnard2023generic}. These first-order logic-based approaches prioritize generality, using quantifiers and predicate symbols to formalize the structure and semantics of argumentation frameworks, enabling the expression of universal and existential constraints over arguments and interactions. However, they often introduce syntactic and computational complexity, making direct integration with efficient automated reasoning tools (e.g., SAT solvers) challenging. In contrast, our work differs fundamentally by adopting propositional logic systems as the encoding target for $AFSA$s. By avoiding quantifiers and leveraging propositional connectives and numerical truth-values, our encoding establishes direct model equivalence between $AFSA$s and propositional formulas. This design not only simplifies the syntactic structure of encoded formulas but also aligns with the practical needs of automated reasoning—propositional logic systems are inherently compatible with lightweight solvers, facilitating efficient computation of extensions and equational semantics. Furthermore, our approach also achieves a unique advantage that the same formal encoding formula enables natural extension from 3-valued to fuzzy semantics. This propositional-centric approach bridges the gap between theoretical argumentation semantics and practical reasoning tools, complementing the generality of first-order logic-based work with enhanced computational feasibility.
\section{Conclusion}

This paper advances the study of $AF$s by establishing a formal connection between $AFSA$s and $\mathcal{PLS}$s. Through encoding methods and model equivalence theorems, we provide logical foundations for both complete semantics and numerical equational semantics of $AFSA$s, bridging the gap between argumentation theory and logical formalisms. The main contributions of this paper are as follows:
\begin{itemize}
	\item \textbf{Syntactic Unification}: This work achieves a syntactic unification and generalization of argumentation frameworks. We first generalize the syntactic structure from $HLAF$s to $BHAF$s, and then integrate the structures of $BHAF$s and $SETAF$s to establish the $HSAF$, which features a more general architecture. This progression realizes both the generalization and unification of syntactic structures across different framework types.
	\item \textbf{Semantic Encoding}: For each $AFSA$, we define complete semantics by generalizing existing definitions and translate these semantics by $\mathcal{PL}_3^L$. Additionally, we propose numerical equational semantics based on $\mathcal{PL}_{[0,1]}$s, providing quantitative methods to evaluate argument acceptability. We also introduce transformations from any $AFSA$ to an $SETAF$ or an $HSAF$.
	\item \textbf{Model Equivalence}: We prove model equivalence theorems showing that assignments (labellings) of an $AFSA$ under a given semantics are exactly the models of the encoded formula in the corresponding $\mathcal{PLS}$. This establishes that the semantics of $AFSA$s can be logically characterized by propositional logic, providing a rigorous foundation for their analysis.
	\item \textbf{Establishment of the Formal Equational Approach}: We introduce a formal equational approach for $HSAF$s. This approach provides a comprehensive numerical equational framework for evaluating the acceptability of arguments. By leveraging Brouwer’s Fixed-Point Theorem, we show that the proposed equational systems for $AFSA$s admit at least one solution, ensuring the feasibility of numerical computations for argument evaluation.
\end{itemize}

This paper offers both theoretical insights and practical tools for reasoning with conflicting information. The model equivalence results and equational semantics established here serve as a robust basis for future developments in argumentation-based applications. 

Building on the preliminary practical application guidelines presented in this work, future research will further deepen the understanding of the practical implications of the proposed equational systems and provide refined guidelines for their selection and deployment in real-world scenarios—with potential extensions to probabilistic interpretations. Additionally, we will extend the proposed encoding methodology to argumentation frameworks incorporating support relations. Another promising avenue lies in transforming the proposed equational semantics into gradual semantics, which offers substantial potential for broader applicability.

\section*{CRediT authorship contribution statement}
\textbf{Shuai Tang:} Conceptualization, Methodology, Investigation, Software, Formal analysis, Writing – original draft, Writing – review \& editing, Validation.
\textbf{Jiachao Wu:} Methodology, Investigation, Software, Validation, Writing – review \& editing, Project administration.
\textbf{Ning Zhou:} Software, Investigation, Writing – review \& editing.

\section*{Declaration of competing interest}
The authors declare that they have no known competing financial interests or personal relationships that could have appeared to influence the work reported in this paper.
\section*{Data availability}
No data was used for the research described in the article.

\section*{Funding sources}
This research did not receive any specific grant from funding agencies in the public, commercial, or not-for-profit sectors.

\section*{Acknowledgement}
We sincerely thank the editors and anonymous reviewers for their valuable comments and constructive suggestions. 




\bibliographystyle{elsarticle-num} 
\bibliography{ref}



\end{document}